\documentclass[a4paper,12pt,twoside]{article}%
\def\ueberschrift{Vectors of Higher Rank on a Hadamard Manifold with Compact Quotient}
\usepackage{ifthen}
\usepackage{makeidx}
\usepackage{calc}
\usepackage[T1]{fontenc}
\usepackage{graphicx}\usepackage{psfrag}
\usepackage{float}

\usepackage[latin1]{inputenc}  
\usepackage[german,english,british]{babel}
\usepackage{amsmath}
\usepackage{amsfonts} 
\usepackage{amscd}\usepackage{amssymb}\usepackage{latexsym}
\usepackage{array} \usepackage{delarray}
\usepackage{theorem}

\newlength{\meinbuffer}
\renewcommand{\emph}{\textbf}
\newcommand{\drunter}[2]{\underset{#1}{\underbrace{#2}}}
\theoremstyle{break}
\newtheorem{Theo}{Theorem}[section] 
\newenvironment{theo}[3] 
          {\ifthenelse{\equal{#1}{}}{\begin{Theo}\label{#2}\index{#3}}
          {\begin{Theo}[#1]\label{#2}\index{#3}}}
        {\end{Theo}} 
\newenvironment{theop}[2]
        {\emph{Proof of Theorem~\ref{#1}:}\textbf{#2}}
        {\hfill$\Box$\\ \vspace{\parsep}}
\newtheorem{Prop}[Theo]{Proposition} 
\newenvironment{prop}[3] 
        {\ifthenelse{\equal{#1}{}}{\begin{Prop}\label{#2}\index{#3}}
          {\begin{Prop}[#1]\label{#2}\index{#3}}}
        {\end{Prop}} 
\newenvironment{propp}[2]
        {\emph{Proof of Proposition~\ref{#1}:}\textbf{#2}}
        {\hfill$\Box$\\ \vspace{\parsep}}

\newtheorem{Lemm}[Theo]{Lemma} 
\newenvironment{lemm}[3] 
        {\ifthenelse{\equal{#1}{}}{\begin{Lemm}\label{#2}\index{#3}}
          {\begin{Lemm}[#1]\label{#2}\index{#3}}}
        {\end{Lemm}} 
\newenvironment{lemmp}[2]
        {\emph{Proof of Lemma~\ref{#1}:}\textbf{#2}}
        {\hfill$\Box$\\ \vspace{\parsep}}

\newtheorem{Coro}[Theo]{Corollary} 
\newenvironment{coro}[3] 
        {\ifthenelse{\equal{#1}{}}{\begin{Coro}\label{#2}\index{#3}}
          {\begin{Coro}[#1]\label{#2}\index{#3}}}
        {\end{Coro}} 
\newenvironment{corop}[2]
        {\emph{Proof of Corollary~\ref{#1}:}\textbf{#2}}
        {\hfill$\Box$\\ \vspace{\parsep}}

\newtheorem{Defi}{Definition}[section] 
\newenvironment{defi}[1] 
        {\ifthenelse{\equal{#1}{}}{\begin{Defi}}
          {\begin{Defi}[#1]}}
          {\end{Defi}}

\newtheorem{Rema}{Remark}[section]
\newenvironment{rema}[3]
        {\ifthenelse{\equal{#1}{}}{\begin{Rema}\label{#2}\index{#3}}
          {\begin{Rema}[#1]\label{#2}\index{#3}}}
          {\end{Rema}}
\newcommand{\U}{{\cal U}}
\renewcommand{\frak}{\mathfrak}
\renewcommand{\flat}{\operatorname{flat}}
\newcommand{\bars}[1]{\bar{\bar#1}}

\newcommand{\inv}{^{-1}}
\newcommand{\horo}{{\cal H}}
\newcommand{\sphere}{{\cal S}}
\newcommand{\W}{{\cal W}}
\newcommand{\R}{{\cal R}}

\newcommand{\bbR}{$\mathbb{R}$ } 
\newcommand{\mathbbR}{\mathbb{R}}
\newcommand{\mathbbN}{\mathbb{N}}
\newcommand{\Td}{\operatorname{Td}}
\newcommand{\Hd}{\operatorname{Hd}}\newcommand{\diam}{\operatorname{diam}}
\newcommand{\grad}{\operatorname{grad}}
 
\newcommand{\rank}{\operatorname{rank}} \renewcommand{\Im}{\operatorname{Im}}
\newcommand{\id}{\operatorname{id}}

\newcommand{\Arcosh}{\operatorname{Arcosh}}
\newcommand{\Ric}{\operatorname{Ric}}
\newcommand{\spann}{\operatorname{span}}
\newcommand{\sdim}{\operatorname{s-dim}}
\newcommand{\Gam}{$\Gamma$}
\newcommand{\neuesubsection}[1]
     {\ifthenelse{\equal{#1}{}}
         {}
         {\markright{\small\thesubsection{ }#1}
         }
     }
\newcommand{\neuesection}[1]
     {\ifthenelse{\equal{#1}{}}
         {}
         {\markboth{\small\thesection.{ }#1}
                   {\small\thesection.{ }#1}
          \thispagestyle{empty}
         }
     }
\newcommand{\mysection}[2]{\cleardoublepage
                           \section{#1}
                           \ifthenelse{\equal{#2}{}}
                                      {\neuesection{#1}}
                                      {\neuesection{#2}}
                          }
\newcommand{\mysubsection}[2]{\subsection{#1}
                              \ifthenelse{\equal{#2}{}}
                                         {\neuesubsection{#1}}
                                         {\neuesubsection{#2}}
                          }
\makeindex

\begin{document}
\pagestyle{myheadings}     
\pagenumbering{roman}

\selectlanguage{german}\thispagestyle{empty}
\begin{center}
{\bf
  {\Large\ueberschrift \\ \vfill}$ $ \\ \vfill
  Dissertation %\\ \vfill
  zur          %\\ \vfill
  Erlangung  der   \\
                 naturwissenschaftlichen Doktorwürde  %\\
  (Dr.\ sc.\ nat.)\\ \vfill
  vorgelegt der\\ \vfill
  Mathematisch-naturwissenschaftlichen Fakultät\\     % \vfill
  der %\\                                               % \vfill
  Universität Zürich\\ \vfill
  von\\ \vfill
  Bert Reinold\\                                       \vfill
  aus\\                                                \vfill
  Altnau (TG)\\ und\\
  Köln-Porz, Nordrhein-Westfalen\\Bundesrepublik Deutschland\\ \vfill
$ $\\ \vfill

  Begutachtet von\\ \vfill
  Prof. Dr. Viktor Schroeder\\
  Prof. Dr. Keith Burns\\ \vfill
  Zürich 2003

}

\end{center}
\pagebreak
\thispagestyle{empty}
$ $
\vfill

Die vorliegende Arbeit wurde von der Mathematisch - naturwissenschaftlichen
Fakultät der Universität Zürich auf Antrag von Prof. Dr. Thomas Kappeler
und Prof. Dr. Viktor
Schroeder als Dissertation angenommen.

\cleardoublepage

\setcounter{page}{1}
\pagenumbering{roman}

\selectlanguage{british}
\addcontentsline{toc}{section}{\protect\numberline{}Abstract \& Zusammenfassung}
\section*{Abstract}
Consider a closed, smooth manifold $M$ of nonpositive sectional
curvature. Write $p:UM\rightarrow M$ for the unit tangent bundle over $M$
and let $\R_>\subset UM$ denote the subset consisting of all vectors of
higher rank. This subset is closed and invariant under the geodesic flow
$\phi$ on $UM$. We will define the structured dimension $\sdim R_>$ which,
essentially, is the dimension of the set $p(\R_>)$ of base points of
$\R_>$.\\

The main result holds for manifolds with $\sdim \R_> <{\dim M}/2$:
For every $\epsilon>0$ there is an $\epsilon$-dense, flow invariant, closed
subset $Z_\epsilon\subset UM\backslash\R_>$ such that $p(Z_\epsilon)=M$.
For every point in $M$ this means that through this point there is a
complete geodesic for which the velocity vector field avoids a neighbourhood of $\R_>$.

\vfill
\selectlanguage{german}
\section*{Zusammenfassung}

Gegeben sei eine geschlossene, glatte Mannigfaltigkeit $M$ nichtpositiver
Schnittkrümmung. Das Einheitstangentialbündel sei mit $p:UM\rightarrow M$
bezeichnet und die Teilmenge aller Vektoren höheren Ranges mit
$\R_>\subset UM$. Diese Teilmenge ist abgeschlossen und invariant unter dem
geodätischen Fluss $\phi$ auf $UM$. Wir definieren die Strukturdimension
$\sdim \R_>$ von $\R_>$, die, im Wesentlichen, die Dimension der
Fußpunktmenge $p(\R_>)$ misst.\\

Das Hauptergebnis gilt unter der Bedingung, dass $\sdim\R_><(\dim M)/2$
gilt: Für jedes $\epsilon>0$ gibt es eine $\epsilon$-dichte,
flussinvariante, abgeschlossene Teilmenge $Z_\epsilon\subset
UM\backslash\R_>$, für die gilt $p(Z_\epsilon)=M$. Dies bedeutet, dass es
durch jeden Punkt in $M$ eine vollständige Geodäte gibt, deren Geschwindigkeitsfeld eine
Umgebung von $\R_>$ vermeidet.

\vfill
\pagebreak
\addcontentsline{toc}{section}{\protect\numberline{}Danksagung}
\section*{Danksagung}
{}
{
Notwendige und hinreichende Bedingung für das Entstehen dieser Arbeit war
die exzellente Betreuung durch meinen Doktorvater Prof.~Viktor                
Schroeder. Ich möchte ihm danken für viele interessante Gespräche im Laufe    
der letzten Jahre und dafür, dass er sich stets Zeit für mich nahm,       
wenn ich eine Frage hatte und er nicht gerade für ein halbes Jahr in           
Chicago weilte.                 \\                                             
                                                         
Meinen Zimmerkollegen Thomas Foertsch und Michael             
Scherrer, sowie allen anderen Institutsmitgliedern danke ich für das           
angenehme Arbeitsklima am Institut für Mathematik. \\                          
                                                                          
Durch diverse Veranstaltungen haben Prof.~Ruth Kellerhals und Prof.~Bruno      
Colbois es in den letzten Jahren verstanden, für einen regelmässigen           
Kontakt zwischen den mathematischen Instituten vieler schweizerischer          
Universitäten zu sorgen. Insbesondere                                          
die Möglichkeit, andere Doktoranden aus meinem Gebiet kennenzulernen und       
mit ihnen gemeinsam neue Dinge zu erarbeiten, hat mir oft neue Motivation      
für meine Arbeit gegeben.                            \\                        
                                                                               
Meinen Eltern danke ich für das Vertrauen in meine Fähigkeiten, welches        
sie nie verloren haben, und für ihre fortwährende Hoffnung, dass meine         
Handschrift eines Tages leserlich sein wird.           \\                      
                                                                               
Zuletzt möchte ich meiner Freundin Anke danken, dafür, dass sie diese          
Arbeit fehlergelesen hat, für ihre Unterstützung und für ihre schier           
unerschöpfliche Geduld insbesondere in den letzten Wochen vor dem              
Abgabetermin.                                            \\                    
                                                                               
Danke schön                                               \\                   
                                                                               
\hspace{1cm}Bert Reinold                                     \\                  
                                                                               
Zürich, im Juli 2003        

}

\selectlanguage{british}

\newpage
{\thispagestyle{empty}
\tableofcontents
\addcontentsline{toc}{section}{\protect\numberline{}Table of Contents}
}
\cleardoublepage

\setcounter{page}{1}
\pagenumbering{arabic}
\mysection{Preliminaries}{}\label{preliminaries}

Before going into detail the introduction in Subsection~\ref{introduction.tex}
will give a rough outline of what this thesis is concerned with. In
Subsection~\ref{structure.tex} the structure of the text is explained section by
section. Then Subsection~\ref{grundlagen.tex} will introduce the main
notation used throughout the thesis. The theory of Riemannian manifolds of
nonpositive curvature as far as considered well known is summed up in
Subsection~\ref{nichtpositiv.tex} for reference.

\mysubsection{Introduction}{}\label{introduction.tex}

Given a complete Riemannian Manifold $M$ there is a natural action of $\mathbbR$
on the unit tangent bundle $UM$, which is called the geodesic flow and is
denoted by $$\phi:\mathbbR\times UM\rightarrow UM.$$
If $\gamma_v$ denotes the unique geodesic defined by the initial condition
$\dot\gamma_v(0)=v$ then $\phi_t(v)$ is the tangent vector to $\gamma_v$ at
time $t$
$$\phi_t(v)=\frac d{ds}\big|_{s=t}\gamma_v(s).$$

\paragraph{Flow Invariant Subsets of $UM$}$ $

Suppose $\gamma_v$ is a simply closed geodesic of length $l$. Then the
orbit of $v$ under the geodesic flow is isometric to a circle of perimeter
$l$. This is the simplest example of a flow invariant subset of $UM$.\\
Another simple example of a flow invariant subset of $UM$ is the unit
tangent bundle $UN$ of a \emph{totally geodesic}\index{totally geodesic
  submanifold}\index{submanifold!totally geodesic}%
\footnote{A submanifold $N\subset M$ is called totally geodesic if every
  geodesic in $N$ is a geodesic in $M$.}
submanifold $N\subset M$.\\
What other examples of flow invariant subsets of $UM$ can we think
of?\\
In the case where $M$ is a compact manifold of strictly negative sectional
curvature, there are many interesting results. E.\ g.\ the geodesic flow is
\emph{ergodic}\index{ergodic flow}\index{flow!ergodic}%
\footnote{ergodic means that any flow invariant measurable subset of $UM$ has
  either full measure or zero measure in $UM$.}
and for almost every geodesic $\gamma$ in $M$ the set $\dot\gamma(\mathbbR)$ is dense
in $UM$. There are many references for these
results, see for example \cite{Pesin1981}.

\paragraph{From Negative to Nonpositive Curvature}$ $

Neither ergodicity nor the existence of a geodesic dense in $UM$ holds for
the flat torus. But what about those manifolds that are neither flat nor
strictly negatively curved?\\
The manifolds of nonpositive curvature are classified by their rank. The
rank of a geodesic is defined to be the maximal number of linearly
independent parallel Jacobi fields along that geodesic. For a manifold the
rank is the smallest rank of its geodesics.\\
By a result of Werner Ballmann\index{Ballmann, Werner} \cite[Appendix 1]{BGS} closed manifolds of higher rank
(i.\ e.\ rank greater than one)
are either products or locally symmetric spaces.\\
However, manifolds of rank one are similar to manifolds of negative curvature in
many aspects. We will use the fact that
geodesics in the vicinity of rank one vectors show a behaviour similar to
the behaviour of geodesics in hyperbolic space rather than in Euclidean
space (compare Proposition~\ref{rank 1 hyper}).

\paragraph{Objective}$ $

The aim of this thesis is to generalize an unpublished result of Keith
Burns\index{Burns, Keith} and Mark Pollicott\index{Pollicott, Mark} for manifolds of constant negative curvature to the
rank one case.\\
For a compact manifold $M$ of constant negative curvature and a given
vector $v_0\in UM$, there is a subset $Z\subset UM$ with the following
properties:
\begin{align}
  \underset{\,\qquad t\in\mathbb R\qquad\,}{\forall}\quad\phi_t(Z)&=Z
\tag{flow invariant}
\index{flow invariant subset}\index{subset!flow invariant}\\
  \underset{v_0\in\W\subset UM \text{ open}}{\exists}\quad \W\cap Z&=\emptyset
\tag{not dense}
\index{not dense subset}\index{subset!not dense}\\
  p(Z)&=M.
\tag{full}
\index{full subset}\index{subset!full}
\end{align}
The construction works in a complete (not necessarily compact) Riemannian manifold $M$ of dimension at
least three and curvature bounded away from zero as Viktor Schroeder showed
in~\cite{Schroeder2001}.\\
For a compact manifold $M$ of nonpositive
curvature Sergei Buyalo\index{Buyalo, Sergei} and Viktor Schroeder\index{Schroeder, Viktor} proved in ~\cite{whoknows}
that this set $Z$ can be constructed for every vector $v_0\in UM$ of rank one.
To construct $Z$ they
proved the following proposition:

\begin{prop}{{\cite{whoknows}}}{general}{}
  Let $M$ be a manifold of nonpositive curvature and dimension at least
  three and let $v_0\in UM$ be a vector
  of rank one. Then for any point $o\in M$ we can find an $\epsilon>0$ and a
  geodesic $\gamma_o$ passing through $o$ with the property that the
  distance between $\dot\gamma_o(t)$ and $v_0$ is greater than $\epsilon$
  for all $t\in \mathbb R$.
\end{prop}

If we take a closer look at the proof we notice that $\epsilon$ can be
chosen to depend continuously on $o$ and thus if $M$ is compact we can
choose $\epsilon$ globally for all $o\in M$. Thus if we define
$$Z:=\underset{o\in M}{\bigcup}\dot\gamma_o(\mathbb R)$$
and $\W$ denotes
the $\epsilon$-neighbourhood of $v_0$ in $UM$, we have the desired properties.

So for every rank one vector $v_0$ there is a flow invariant, full subset
of $UM$ the complement of which contains a neighbourhood of $v_0$.\\
In fact, $Z$ can be chosen to be $\epsilon'$-dense for some
$\epsilon'=\epsilon'(\epsilon)>0$.\\

The main result of this thesis is a generalization of this result: Under
certain conditions on the set of higher rank vectors there is a full, flow
invariant $\epsilon'$-dense subset of $UM$ avoiding an
$\epsilon$-neighbourhood of all vectors of higher rank.\\
But first we need to understand the set of vectors of higher rank.

\paragraph{Higher Rank Vectors in a Rank One Manifold}$ $

Given a manifold $M$ of rank one the vectors of higher rank ($\rank(v)>1$)
form a closed subset, say $\R_>$. If $M$ is real analytic this subset is
subanalytic. Projecting $\R_>$ down to $M$ gives another
subanalytic set, say $Y$.
In Section~\ref{subanalytic.tex} we will describe subanalytic subsets in
more detail. For the moment it suffices to think of them as being
stratified%
\footnote{i.\ e.\ they can be written as locally finite union of
  submanifolds, compare Definition~\ref{strata}}
by a locally finite family of real analytic submanifolds.\\

So for a real analytic manifold of rank one the vectors of higher rank are
all tangent to a locally finite union of submanifolds of $M$. This fact
motivates the definition of the \emph{structured dimension} of a subset of
$UM$ in Definition~\ref{def s-dimension}. For flow invariant subsets, like
$\R_>$ this is essentially the dimension of the set of base points:

\paragraph{Definition (s-Dimension)}\hfill{(compare Remark~\ref{sdimflowinv})}\\

\begin{it}
Suppose $\R\subset UM$ is a flow invariant subset of the unit tangent
bundle.
\begin{itemize}
\item An {\bf s-support} of $\R$ is a locally finite union $N=\bigcup N_i$
  of closed submanifolds of $M$ such that $p(\R)\subset N$. The {\bf
    dimension of $N$} is defined to be $\max\dim N_i$.
\item The {\bf structured dimension} $\sdim \R_>$ of $\R_>$ is the minimal
  dimension of an s-support of $\R$.
\end{itemize} 
Even though there
might be different s-supports for $\R$ and different stratifications for
the same s-support, the structured dimension of $\R$ is well defined.
\end{it}\\

For any subset $\R\subset UM$ we have 
$$0\leq \sdim(\R)\leq\dim(M).$$
Obviously nontrivial flow invariant subsets of $UM$ contain at least one
geodesic and hence have structured dimension at least one. For a closed totally
geodesic submanifold $N\subset M$ the structured dimension of $UN$ is just
the dimension of $N$: $\sdim(UN)=\dim(N)$ and $N$ is an s-support of $UN$. Any full subset of $UM$
(i.\ e.\ it covers all of $M$ under the base point projection) has structured
dimension $\dim(M)$.\\
So the simplest example of a flow invariant subset of s-dimension one is
the set of tangent vectors to a a finite collection of simply closed geodesic. 

\paragraph{Main Theorem}$ $

Now we have the means to state the main theorem which will be proved in Section~\ref{topologie.tex}.
\paragraph{Theorem~\ref{THE Theorem}}$ $\\
{\it
Let $M$ be a compact manifold of nonpositive curvature. Suppose the
s-dimension of the set of vectors of higher rank $\R_>$ is bounded by 
$$\sdim(\R_>)<\frac{\dim M}2.$$ 
Then for every
  $\epsilon>0$ there is a closed, flow invariant, full, $\epsilon$-dense
  subset $Z_\epsilon$ of the unit tangent bundle $UM$ consisting only of
  vectors of rank one.
}\\[\baselineskip]
Notice that $UM\backslash Z_\epsilon$ will be a neighbourhood of
$\R_>$. Flow-invariant means that $Z_\epsilon$ consists of velocity fields
of geodesics. Since $Z_\epsilon$ is full, we can find a geodesic through
every point in $M$, in fact our proof will show that initial vectors of
such geodesics are even $\epsilon$-dense in $U_oM$ for any given point.\\
 
\newpage
\thispagestyle{empty}
\mysubsection{Structure of this Thesis}   {}{}
\label{structure.tex}
Here is a short survey on the sections of this thesis, on their purpose and
the main results.
%1
\paragraph{Section~\ref{preliminaries}:}
This section consists of an introduction, a summary of the thesis and
introduces the main notation.

\paragraph{Section~\ref{vertical.tex2}:}
To define the rank of geodesics and manifolds we need to better understand
the tangent bundle $TTM$ of the tangent bundle $TM$. Elements of this object
correspond to Jacobi fields in the manifold. Working with $TTM$ might seem
a bit complicated at first glance but once we have established the
connection to Jacobi fields the desired results pop up easily. 
\paragraph{Section~\ref{subanalytic.tex}:}
By now we know that the vectors of lowest rank form an open subset of $UM$.
For real analytic manifolds we can describe its closed complement more
precisely. In fact, all the sets of vectors of constant rank form
subanalytic subsets of $UM$. This motivates the definition of the structured
dimension and is applied in Section~\ref{flats.tex}.
\paragraph{Section~\ref{sphaere.tex}:}
This section introduces the notation for spheres and horospheres in a
Hadamard manifold $X$. If the opposite is not explicitely stated we will always
consider spheres and horospheres as objects in $UX$ rather than in $X$, by
identifying outward orthogonal vectors with their base points. Both, in
$X$ and $UX$ the horospheres can be understood as limits of spheres of
growing radius. This convergence is uniform on compact subsets.\\
\\
Summing up all we know about higher rank vectors on real analytic manifolds
of rank one we get a result that stands a little apart from the other results
of this thesis. But it is interesting in its own right. The idea for the
proof arose from a discussion with Sergej Buyalo and Viktor Schroeder.

\paragraph{Corollary~\ref{coro2}}$ $\\
{\it
  Let $X$ be an analytic rank one Hadamard manifold with compact quotient.
  Then for any horosphere or sphere in $UX$ the subset of rank one vectors
  is dense. 
}

\paragraph{Section~\ref{hyperbolic.tex}:}
Here we show that in the vicinity of rank one vectors geodesics have some
widening property that reminds one of geodesics in hyperbolic space.
The proof goes back to Sergei Buyalo and Viktor Schroeder
in~\cite{whoknows}. We improve the result slightly, which will be necessary
in the following.

\paragraph{Section~\ref{rank1.tex}:}
In this section we describe the construction by Sergei Buyalo and Viktor
Schroeder in~\cite{whoknows}\index{Schroeder, Viktor}\index{Buyalo, Sergei}
as far as it is
necessary for the understanding of our construction in the higher rank case.

\paragraph{Section~\ref{ausweichen.tex}:}
Suppose $X$ is a Hadamard manifold with compact quotient $M=X/\Gamma$ and $R\subset X$ a
$\Gamma$-compact submanifold. Fix a
sphere $S$ of large radius, say $L$. Given an outward orthogonal vector to
$S$, we 
want to find a close vector on $S$ which is comparably far away from
$UR$.\\
Obviously we have to impose some dimensional restriction on $R$. But if
$\dim R$ is small enough (when compared to $\dim X$) we can even find a
continuous deformation $\Psi_L$ of $UX\backslash UR$ into
$UX\backslash\W_{\epsilon}(UR)$ which sends every vector to an
vector that is outward normal to the same sphere of radius $L$. This deformation
moves every vector by less than $C$ to a vector that is farther than
$\epsilon$ away from $UR$. \\
The result generalizes to the case where $R$ is stratified by
submanifolds of restricted dimension. This is an important step towards
avoiding the set of higher rank vectors, if its structured dimension is
small.

\paragraph{Section~\ref{rankr5.tex}:}
Putting together what we learned in the other sections we get the first
result. For any given point there are many geodesic rays starting in that
point and avoiding a neighbourhood of all vectors of higher rank%
\footnote{Here $\R_>$ denotes the subset of $UM$ consisting of all vectors
  of higher rank.}:
\paragraph*{Proposition~\ref{summary}}$ $\\
\begin{it}
Let $X$ denote a rank one Hadamard manifold with compact quotient
$M=X/\Gamma$. Suppose $\sdim(\R_>)<\dim X-1$ then there are constants $\epsilon$ and $c$ such that for any $\eta<\epsilon$
there is an $\eta'<\eta$ for which the following holds: 
\begin{quote} 
For any compact manifold $Y$ with 
$$\dim Y<\dim X -\sdim \R_>$$ and any continuous map
$v_0:Y\rightarrow U_oX$ from $Y$ to the unit tangent sphere at a point
$o\in X$ we can find a continuous map $v_\infty:Y\rightarrow U_oX$ which is
$c\eta$-close to $v_0$ and satisfies $$d({\R_>},\phi_{\mathbbR_+}(v_{\infty}))\geq\eta'.$$
\end{quote}
\end{it}

So on a surface we have many geodesic rays avoiding all vectors of higher
rank if there are only finitely many geodesics of higher rank which are all
closed. There are many examples of surfaces for which this is true, as the
construction by Werner Ballmann, Misha Brin and Keith Burns in~\cite{BBB}
shows.\index{Ballmann, Werner}\index{Brin, Misha}\index{Burns, Keith}

\paragraph{Section~\ref{topologie.tex}:}
So now we know that every point is the initial point of many geodesic rays
which do not come close to any vector of higher rank. But we are interested
in flow invariant subsets of $UM$ and therefore we have somehow to find two
geodesic rays adding up to give a full geodesic. This can be done using a
topological argument on $U_oX$, the unit tangent sphere at our starting
point $o$. This gives a stronger restriction on the structured dimension of
$\R_>$:
\paragraph{Theorem~\ref{maintheorem}}$ $\\
{\it
Let $X$ denote a rank one Hadamard manifold with compact quotient $M$,
and suppose that
$$\sdim(\R_>)<\frac{\dim X}2.$$
Then there are constants $\epsilon$ and $c$ such that for any
$\eta<\epsilon$ there is an $\eta'<\eta$ with the following property:\\
For every vector $v_0\in
UX$ there is a $c\eta$-close vector $v$ with the same base point such that
$\gamma_v$ avoids an $\eta'$-neighbourhood of all vectors of higher rank.
}\\

An immediate consequence is our main result,

\paragraph{Theorem~\ref{THE Theorem}}$ $\\
{\it
Let $M$ be a compact manifold of nonpositive curvature. Suppose the
s-dimension of the set of vectors of higher rank $\R_>$ is bounded by 
$$\sdim(\R_>)<\frac{\dim M}2.$$ Then for every
  $\epsilon>0$ there is a closed, flow invariant, full, $\epsilon$-dense
  subset $Z_\epsilon$ of the unit tangent bundle $UM$ consisting only of
  vectors of rank one.
}\\

\paragraph{Appendix
:}
This appendix gives a short introduction to the theory of manifolds and
Riemannian manifolds without proofs. Many of the terms are used throughout
this thesis and some of the notation might not be standard.
 
\newpage
\thispagestyle{empty}
\mysubsection{Notation}{}\label{grundlagen.tex}
In this section we introduce the notation as used throughout the text. The
terminology should be standard. However, should problems arise, please
refer to the appendix on smooth manifolds or to the index to find
definitions or further explanations. A general reference is the book by 
Takashi Sakai~\cite{Sakai1996}.

\paragraph{Manifold}$ $\\
$M$\index{M@$M$}\index{$m$@$M$} denotes a compact, nonpositively curved Riemannian manifold
of dimension $n$. Usually we will expect $M$ to be smooth ($C^\infty$), but
sometimes it will be useful to consider the real analytic case. If $M$ is
an analytic manifold this will be said explicitly.

\paragraph{Tangent Bundle}$ $\\
By $TM$ %\index{TM@$TM$} 
we denote the
{tangent bundle} 
over
the manifold $M$. The unit tangent bundle will be denoted by $UM$. The base
point projection is the bundle map $p:TM\rightarrow M$, respectively
$p:UM\rightarrow M$. $TM$ is a $2n$-dimensional Riemannian manifold with
$(2n-1)$-dimensional submanifold $UM$, both of the same differentiability as
$M$.\\

We write $\langle.,.\rangle_p$ for the
scalar product on $T_pM$ given by the  
Riemannian metric and might omit the indexed $p$ where there is no danger
of ambiguity.\\

$\langle.,.\rangle_p$ defines the distance function $d(.,.)$ on $M$. A
geodesic $\gamma:\mathbbR\rightarrow M$ will be a geodesic with respect to
the metric $d$. For $v\in TM$ the geodesic $\gamma_v$ is the unique
geodesic with $\gamma_v(0)=p(v)$ and $\dot\gamma_v(0)=v$. Most of the time
we will consider unit speed geodesics only, i.\ e.\ $\gamma_v$ with $v\in
UM$.
The geodesic flow at time $t$ will be denoted by $\phi_t$. It can be
understood as an action of $\mathbbR$ on either $TM$ or $UM$.

\paragraph{Universal Covering}$ $\\
The universal covering of $M$ will be denoted by $\pi:X\rightarrow M$.\\

In
general $X$ will denote a Hadamard manifold with compact quotient, i.\ e.\
a complete, simply connected manifold of nonpositive curvature which admits
an action of a group $\Gamma$ such that $M=X/\Gamma$ is a compact manifold
covered by the projection $\pi:X\rightarrow X/\Gamma$. Notice that in this
case  $\Gamma$\index{Gamma@$\Gamma$}\index{$gamma$@$\Gamma$} denotes the group of \emph{deck
  transformations}\index{deck transformations}\index{group!deck
  transformation} of $\pi:X\rightarrow M$, i.\ e.\ the subgroup of isometries $\sigma$ of $X$
satisfying $\pi\circ\sigma=\pi$.\\

We write $TX$, respectively $UX$ for the (unit) tangent bundle over $X$ and
write $p$ for the base point projection. The covering map $\pi$ induces a
covering map $d\pi:UX\rightarrow UM$ by the following diagram:

$$\begin{CD}
UX@>{d\pi}>>UM\\
@V{p}VV @VV{p}V\\
X@>{\pi}>>M.
\end{CD}$$
$X$ carries a Riemannian structure induced by the structure on $M$ and
denoted by
$$\langle u,v\rangle_x:=\langle
d\pi(u),d\pi(v)\rangle_{\pi(x)}\;\;\left(=\pi^*\langle u,v\rangle\right).$$

\paragraph{$\Gamma$-Invariance}$ $\\
Throughout this text we will be concerned with structures of $X$ and $UX$
that arise from structures on $M$ and $UM$. For example if you take the
preimage under $\pi:X\rightarrow M$ of a compact subset of $M$, it is not
necessarily compact itself. But still it will inherit some nice properties
from the underlying compact set. We use the following terms.\\

A subset $C$ of $X$ is called
\emph{$\Gamma$-invariant}\index{Gamma-invariant@$\Gamma$-invariant subset}
\index{subset!$\Gamma$-invariant}\index{invariant subset} if $\Gamma C=C$, i.\ e.\ $\sigma(p)\in C$ for
all $\sigma\in \Gamma$ and $p\in C$. A subset $C\subset UX$ is called
\emph{$\Gamma$-invariant}, if for all $v\in C$ and all $\sigma\in \Gamma$
it holds
$d\sigma(v)\in C$.%
\footnote{Notice that $UM=UX/\tilde \Gamma$ for
  $\tilde\Gamma:=\{d\sigma\;|\;\sigma\in\Gamma\}$ if $M=X/\Gamma$. So saying
  that $C\subset UX$ is $\Gamma$-invariant as a subset of the unit tangent
  bundle of $X$ is equivalent to saying that
  it is $\tilde\Gamma$-invariant as a subset of the manifold $UX$.}
For any element $x$ of $X$, respectively $UX$ the \emph{$\Gamma$-orbit} \index{Gamma-orbit@$\Gamma$-orbit}\index{orbit} of $x$ is
the smallest $\Gamma$-invariant subset of $X$, respectively $UX$ containing
$x$. A map is called \emph{$\Gamma$-compatible}\index{Gamma-compatible@$\Gamma$-compatible}\index{compatible} (or compatible with the
compact quotient structure) if it is constant on orbits. A subset $C$ of
$X$, respectively $UX$, is called
\emph{$\Gamma$-compact}\index{Gamma-compact@$\Gamma$-compact}
\index{compact!$\Gamma$-compact} if it is $\Gamma$-invariant and has
compact image under the projection $\pi: X\rightarrow M$.\\
Remember that in a first countable space compactness coincides with
sequential compactness. The following lemma is an equivalent to sequential
compactness\index{sequential compactness}\index{compact!sequentially}
for $\Gamma$-compact sets. 

\begin{lemm}{$\Gamma$-Compactness}{proc}{}
\index{Gamma-Compactness@$\Gamma$-Compactness}
Suppose $X$ is a Riemannian manifold with compact Riemannian quotient
$M=X/\Gamma$. Suppose 
$C$ is a \Gam-compact subset of $X$ and $c:X\rightarrow Q$ a
\Gam-compatible map.\\
Then for any sequence $(y_i)_{i\in\mathbbN}$ of points in $C$ such that $c(y_i)\rightarrow
c_0$ there is a sequence $(y'_i)_{i\in\mathbbN}$ of points in
$C$ such that
\begin{enumerate}
\item $y'_i\rightarrow y\in C$ and
\item $(\pi (y'_i))_{i\in\mathbbN}$ is a subsequence of
  $(\pi (y_i))_{i\in\mathbbN}$ and hence $c(y'_i)\rightarrow c_0$.
\end{enumerate}
Note that $c(y)=c_0$ if $c$ is continuous.
\end{lemm}
\begin{lemmp}{proc}{}
Since $\pi (C)$ is compact we can find a converging subsequence
$\pi (y_{\iota(i)})\rightarrow z\in \pi (C)$. Fix $y\in \pi \inv(z)$ and a countable
basis of neighbourhoods $\{U_k\}$ of $y$ in $X$. The images $\pi (U_k)$ are
neighbourhoods of $z$ in $M$. There are therefore
infinitely many $\pi (y_i)$ contained in any of the $\pi (U_k)$. Fix a sequence $z_i=\pi (y_{j(i)})\in \pi (U_i)$ where $y_{j(i)}$ is a
subsequence of $y_i$ as follows. Define $j(0):=0$ and recursively
$$j(i+1):=\min\{j>j(i)\;|\;\pi (x_j)\in \pi (U_{i+1})\}\quad\text{ and }\quad z_{i+1}:=\pi (x_{j(i+1)}).$$
Now choose $y'_i\in \pi \inv(z_i)\cap U_i$ to get the desired sequence.
\end{lemmp}

\paragraph{Neighbourhoods}$ $\\
There is a canonical Riemannian metric, the Sasaki metric, on the manifolds $TM$, $TX$, $UM$ and
$UX$ which will be introduced in Section~\ref{vertical.tex}. So we have
scalar products $\langle.,.\rangle$ on the tangent spaces and we can talk
about distances $d(u,v)$ between tangent vectors in e.\ g.\ $UM$ or $TX$ and
hence about neighbourhoods in these spaces.
In general we will use the same notations for vector bundles as for
the original manifolds, since they are manifolds, too. With one
exception, though:\\
Sometimes we will have to talk about neighbourhoods of points and
neighbourhoods of vectors in their respective spaces. To distinguish these
we use the following notation. Write 
$$\U_\epsilon(p):=\{q\in X\;|\;d(p,q)<\epsilon\}\index{$u epsilon(p)$@$\U_\epsilon(p)$}\index{U epsilon(p)@$\U_\epsilon(p)$}$$
for the \emph{$\epsilon$-neighbourhood of a point $p\in X$} and \index{neighbourhood}
$$\W_\epsilon(v):=\{u\in UX\;|\; d(u,v)<\epsilon\}\index{W
  epsilon(v)@$\W_\epsilon(v)$}\index{$w epsilon(v)$@$\W_\epsilon(v)$}$$
for the \emph{$\epsilon$-neighbourhood of a vector $v\in UX$}. We will see that
for two vectors $u,v$
$$d(p(u),p(v))\leq d(u,v)\quad\text{and}\quad p(\W_\epsilon(v))=\U_\epsilon(p(v)).$$

\paragraph{Spheres and Horospheres}$ $\\
A sphere in $X$ can be defined by its radius $r$ and one outward pointing
vector $v\in UX$. We will do this and write $S_v(r)$ for this sphere in
$X$. In fact we will be more interested in the set of outward pointing
vectors to this sphere. This is a subset of $UX$ which we will call the
sphere of radius $r$ in $UX$ defined by $v$ and we will denote by
$\sphere_v(r)$. \\
 
A similar notation will be used for horospheres. $H_v(r)$ will denote a
horosphere in $X$ while  $\horo_v(r)$, the set of outward pointing vectors of
$H_v(r)$, will be called a horosphere in $UX$.

\paragraph{Rank}$ $\\
For every vector $v\in UX$ and every geodesic in $X$ we define its rank
$\rank(v)\in\{1,2,\ldots,n\}$. The subset of $UX$ consisting of all vectors
of rank one will be denoted by $\R_1$. A vector $v\in UX$ with $\rank(v)>1$
will be called a vector of higher rank. The set of all vectors of higher
rank will be denoted by $\R_>$.%

\newpage
\thispagestyle{empty}
\mysubsection{Manifolds of Nonpositive Sectional Curvature}{Nonpositive
Curvature}
\label{nichtpositiv.tex}
In this section we state results from~\cite{BBE} and~\cite{BGS} for a
manifold $X$ of nonpositive curvature. Several useful structures and
notions are introduced which we 
will use throughout the text.

\paragraph{Convexity}$ $\\
Recall that a function $f:\mathbbR\rightarrow\mathbbR$ is
\emph{convex}\index{convex function}\index{function!convex} if
 $f(t)\leq f(a)+\frac{t-a}{b-a}(f(b)-f(a))$ holds for all $a<t<b$. If $f$ is
smooth then this is equivalent to $f''\geq 0$.
A map $f:M\rightarrow\mathbbR$ from a smooth manifold with linear
connection into the reals is called \emph{convex}\index{convex map}\index{map!convex} if it is convex for every
restriction to a geodesic $f\circ\gamma:\mathbbR\rightarrow\mathbbR$. 

\begin{lemm}{}{convexlemma}{}
Suppose $M$ is a Riemannian manifold of nonpositive sectional
curvature. Let $\gamma$ denote a geodesic and $J$ a Jacobi field along
$\gamma$. Then the function
$$t\longrightarrow\|J_t\|$$
is convex.\\
If $M$ is furthermore simply connected then the following maps
are convex, too,
\begin{align*}
t&\longrightarrow d(\gamma(t),\sigma(t))\\
(p,q)&\longrightarrow d(p,q)\\
p&\longrightarrow d(p,N).
\end{align*}
where $\sigma$ is another geodesic and $N$ is a \emph{totally
  geodesic}\index{totally geodesic submanifold}
\index{submanifold!totally geodesic}% 
\footnote{I.\ e.\ all geodesics in $N$ are geodesics in $M$}
submanifold of $M$. 
\end{lemm}

It is easy to see that $\|J\|$ is convex:
\begin{align*}
\|J\|''&=\sqrt{\langle J,J\rangle}''=\frac12\left(\frac{\langle
    J',J\rangle}{\sqrt{\langle J,J\rangle}}\right)'\\
       &=\frac1{2\|J\|^3}\left(\langle J'',J\rangle\|J\|^2
       +\|J'\|^2\|J\|^2-\frac12 \langle J',J\rangle^2\right)\\
       &=\frac1{2\|J\|^3}\Big(-\drunter{\sim K(J,\dot\gamma)\leq0}{\langle
         R(J,\dot\gamma)\dot\gamma,J\rangle}\|J\|^2 +\frac12\|J'\|^2\|J\|^2
       +\frac12\drunter{\geq0}{\left(\|J'\|^2\|J\|^2-\langle
           J',J\rangle^2\right)}\Big)\\
       &\geq 0
\end{align*} 

The means to prove the second part of the lemma is the \emph{second variation formula}\index{second
variation formula}\index{variation formula!second}. Suppose $h_s(t)$ is a
variation of the geodesic $\gamma(t)=h_0(t)$ and denote by
$X_t=\frac\partial{\partial s}\big|_{s=0}h_s(t)$ the variational
vector field along $\gamma$. Then for the length function $L(s):=L(h_s)$
we have 
\begin{align*}
L''(0)=&\frac1{\|\dot\gamma\|}\int_a^b\left(\left\langle\frac\nabla{dt}X^\perp,\frac\nabla{dt}X^\perp\right\rangle-
  \left\langle R(X^\perp,\dot\gamma)\dot\gamma,X^\perp\right\rangle\right)\;dt\\
&-\frac1{\|\dot\gamma\|}\left.\left\langle\frac\nabla{\partial s}\frac{\partial
    h}{\partial s},\dot\gamma\right\rangle\right|_{(0,a)}
 +\frac1{\|\dot\gamma\|}\left.\left\langle\frac\nabla{\partial s}\frac{\partial
    h}{\partial s},\dot\gamma\right\rangle\right|_{(0,b)}
\end{align*}
where $X^\perp:=X-\frac{\left\langle
    X,\dot\gamma\right\rangle}{\|\dot\gamma\|^2}\dot\gamma$\index{X
  ortho@$X^\perp$}\index{$x ortho$@$X^\perp$} and $[a,b]$ is
the range of $t$. 

\paragraph{Toponogov's Theorem}$ $\\
An important tool when dealing with manifolds of bounded curvature are
comparisons with model spaces of constant
curvature.\\

Consider three points $p,x,y\in X$ and write $\gamma$, resp.\ $\sigma$, for
the shortest geodesics joining $p$ to $x$, resp.\ $y$, parametrized such
that $\gamma(0)=p=\sigma(0)$ and $\gamma(1)=x$ and $\sigma(1)=y$.
By Lemma~\ref{convexlemma} the function $t\mapsto d(\gamma(t),\sigma(t))$
is convex and hence
$d(\gamma(s),\sigma(s))\leq sd(x,y)$ holds for all $s\in[0,1]$.
In Euclidean space equality holds and we conclude 
\begin{lemm}{}{}{}
Suppose $p,x,y\in X$ and $p',x',y'\in\mathbbR^2$ are triples of points in a
manifold $X$ of nonpositive curvature and in $\mathbbR^2$ respectively,
such that $d(p,x)=d(p',x')$, $d(p,y)=d(p',y')$ and $d(x,y)=d(x',y')$. Let
$\gamma$, resp.\ $\sigma$, denote a shortest geodesics joining $p$ to $x$,
resp.\ $y$ (such that $\gamma(1)=x$ and $\sigma(1)=y$), and let $\gamma'$, 
 resp.\ $\sigma'$, denote shortest geodesics joining $p'$ to $x'$,
resp.\ $y'$ (such that $\gamma(1)=x'$ and $\sigma(1)=y'$).\\
Then for all $s\in[0,1]$ 
$$d(\gamma(s),\sigma(s))\leq d(\gamma'(s),\sigma'(s)).$$
\end{lemm}

In fact, if equality holds the triangles $\Delta(p,x,y)$ and
$\Delta(p',x',y')$ are isometric. This is the rigidity result contained in

\begin{theo}{Toponogov's Comparison Theorem}{toponogov}{}
Let $X$ denote a manifold of nonpositive curvature. Suppose $p,x,y\in X$
and $p',x',y'\in\mathbbR^2$ are given points, $\gamma$, resp.\ $\sigma$,
are shortest geodesics from $p$ to $x$, resp.\ $y$, and $\gamma'$, resp.\
$\sigma'$, are shortest geodesics from $p'$ to $x'$, resp.\ $y'$, such that
$$d(p,x)=d(p',x')\qquad d(p,y)=d(p',y')\qquad L(\gamma)=L(\gamma')\qquad L(\sigma)=L(\sigma').$$
Then
\begin{itemize}
\item if $\measuredangle_p(x,y)=\measuredangle_{p'}(x',y')$ then
  $d(x,y)\geq d(x',y')$,
\item if $d(x,y)= d(x',y')$ then $\measuredangle_p(x,y)\leq
  \measuredangle_{p'}(x',y')$
\end{itemize}
and if equality holds in either case then the points $p,x,y$ are contained
in a totally geodesic flat triangle in $X$.
\end{theo}

\begin{rema}{}{}{}
\begin{itemize}
\item For any triple of side lengths satisfying the triangle inequality we
  can find a triangle in $\mathbbR^2$ with the same side length. Similarly
  for any given lengths of two sides and a given angle between these sides
  we can find a triangle in $\mathbbR^2$ with the same constellation.
  Therefore for any given points $p,x,y\in X$ we can find a comparison
  triangle in $\mathbbR^2$ such that Toponogov's Theorem applies.
\item An easy consequence of Toponogov's Theorem is that for any triangle
  in $X$ the sum over the inner angles is less or equal to $\pi$ and
  equality holds only for totally geodesic flat triangles.
\end{itemize}
\end{rema}

\paragraph{Hadamard Manifold}$ $\\
By the \emph{Theorem of
  Hadamard/Cartan}\index{Hadamard/Cartan}\index{Cartan/Hadamard}
  \index{Theorem!Hadamard/Cartan} for a nonpositively curved,
complete manifold the exponential map is a covering map when restricted to
any tangent space. So if we fix any $p\in M$ the map
$\exp_p:T_pM\rightarrow M$ is a covering map. Since the tangent space is
diffeomorphic to $\mathbbR^n$ which is simply connected, we conclude that
the universal covering of a manifold of nonpositive curvature is always
diffeomorphic to $\mathbbR^n$. A simply connected manifold of nonpositive
curvature is called a \emph{Hadamard manifold}.\index{Hadamard
  manifold}\index{manifold!Hadamard}\\

Now let $X$ be a Hadamard manifold%
\footnote{It would be sufficient, if $X$ was simply connected without
  conjugate points. But since we will only work with nonpositive curvature,
  we will only encounter Hadamard manifolds.}
of dimension $n$.
Recall that by Lem\-ma~\ref{convexlemma} the distance function on $X$
is convex.
Two vectors $v,w\in UX$ are said to be
\emph{asymptotic}\index{asymptotic vectors}\index{vector!asymptotic} if $d(\gamma_v(t),\gamma_w(t))$ is bounded for
$t\in\mathbbR_+$. The quotient of $UX$ under this equivalence relation is
called the \emph{boundary at infinity}\index{boundary at infinity}{}, the
\emph{boundary sphere}\index{boundary sphere}{} or the \emph{sphere at
  infinity}\index{sphere!at infinity}. We denote it by
$\partial X$\index{del X@$\partial X$}\index{$del X$@$\partial X$} and write
$\gamma_v(\infty)$\index{gamma v(infinity)@$\gamma_v(\infty)$}\index{$gamma
  v(infinity)$@$\gamma_v(\infty)$} for the equivalence class of $v$.
For any $p\in X$ there are two bijections of $U_pX$ onto $\partial X$ by
$v\mapsto\gamma_v(\infty)$ and $v\mapsto\gamma_{-v}(\infty)
=:\gamma_v(-\infty)$, respectively. These induce the topology of the sphere
on the
boundary. Sometimes it is necessary to consider the union of $X$ and its
boundary sphere. We write $X(\infty):=X\cup \partial X$\index{X infinity@$X(\infty)$}\index{$x(infinity)$@$X(\infty)$}. There is a natural
topology%
\footnote{the cone topology, see~\cite{EberleinHamenstaedtSchroeder1993}\index{topology!cone}\index{cone topology}} on $X(\infty)$ such that the topological
subspaces $X$ and $\partial X$ are equipped with their original topology and
$X$ is a dense, open subset of $X(\infty)$.

\paragraph{Tits Metric}$ $\\
For a Hadamard manifold $X$ there are two complete metrics of interest on $\partial
X$.\\
The \emph{angle metric} $\measuredangle$\index{angle
  metric}\index{metric!angle}\index{$angle(.,.)$@$\measuredangle(.,.)$} is defined by
$$\measuredangle(\xi,\zeta):=\underset{p\in
  X}\sup\;\measuredangle_p(\xi,\zeta)$$ 
where $\measuredangle_p(\xi,\zeta)$ denotes the angle between two geodesic
rays $\gamma$ and $\sigma$ starting in $p$ such that $\gamma(\infty)=\xi$
and $\sigma(\infty)=\zeta$.\\

The \emph{asymptotic growth rate} $l$\index{metric!asymptotic growth rate}
\index{asymptotic growth rate}\index{growth rate}\index{$l(.,.)$}\index{l(.,.)@$l(.,.)$} is
another metric on $\partial X$ defined as follows:\\
Fix $x\in X$. For $\xi,\zeta\in \partial X$ take two unit speed geodesics
$\gamma$ and $\sigma$ starting in $x$ with $\gamma(\infty)=\xi$ and
$\sigma(\infty)=\zeta$. Now define
$$l(\xi,\zeta):=\underset{t\rightarrow\infty}\lim\frac{d(\gamma(t),\sigma(t))}t.$$
This metric is related to the angle metric by the equation
$$l(\xi,\zeta)=2\sin\left(\frac{\measuredangle(\xi,\zeta)}2\right).$$
The \emph{Tits metric} $\Td$ \index{Tits
  metric}\index{metric!Tits}\index{Td@$\Td$}\index{$td$@$\Td$} is the inner
(pseudo)-metric on $\partial X$ with respect to either of the above
metrics, i.\ e.\
$$\Td(\xi,\zeta):=\inf\left\{L(c)\;\big|\;c\text{ continuous curve in
    }\partial X \text{ joining }\xi\text{ to }\zeta\right\}\in[0,\infty].$$

The Tits metric indicates whether two points on the boundary can be joined
by a geodesic in the manifold.

\begin{lemm}{{\protect\cite[p.46]{BGS}}}{}{}
\begin{itemize}
\item If $\Td(\xi,\zeta)>\pi$ then there is a geodesic $\gamma$ in $X$
  joining $\xi$ to $\zeta$ (i.\ e. 
  $\gamma(-\infty)=\xi$ and $\gamma(\infty)=\zeta$).
\item If $\gamma$ is a geodesic in $X$ then 
  $\Td(\gamma(-\infty),\gamma(\infty))\geq\pi$ and equality implies that 
  $\gamma$ bounds a flat half plane.
\end{itemize}
\end{lemm}

\paragraph{Horospheres}$ $\\
Now we can define the \emph{stable}\index{stable bundle}\index{bundle!stable} and
\emph{unstable bundle}\index{unstable bundle}\index{bundle!unstable} over $X$
\begin{align*}
  W^u(v)&:=\{w\in UX\,|\,\gamma_v(-\infty)=\gamma_w(-\infty)\}
  \index{W u (v)@$W^u(v)$}
  \index{$w u (v)$@$W^u(v)$}
\tag{unstable}\\
  W^s(v)&:=\{w\in UX\,|\,\gamma_v(\infty)=\gamma_w(\infty)\}
  \index{W s(v)@$W^s(v)$}
  \index{$w s(v)$@$W^s(v)$}
\tag{stable}
\end{align*}
which are subbundles of $UX$. The partition of $UX$ into these bundles is
called 
\index{stable foliation}\index{foliation!stable}
\index{unstable foliation}\index{foliation!unstable}
the \emph{stable foliation}, respectively \emph{unstable foliation}.
For every $v\in UX$ the leaves $W^u(v)$ and $W^s(v)$ are $n$-dimensional
$C^2$-submanifolds of $UX$. Notice that the leaves are invariant under the
geodesic flow and that $W^u(v)=W^s(-v)$.

We can further foliate the leaves of the unstable foliation by
means of the Busemann functions:\\
Define the \emph{generalized Busemann function}\index{generalized Busemann
  function}\index{Busemann function!generalized}\index{function!generalized
  Busemann}
\begin{align*}
  b:UX\times X&\longrightarrow\mathbbR\\
  (u,p)&\longmapsto
  \underset{t\rightarrow\infty}\lim\,d(\gamma_{-u}(t),p)-t.
\end{align*}
This function is continuous and for every $v\in UX$ the \emph{Busemann
  function}\index{Busemann function}\index{function!Busemann}
$b_v:=b(v,.)$\index{b v@$b_v$}\index{$b v$@$b_v$} is of class $C^2$. For two vectors $u,v\in
UX$ the 
difference $b_v-b_u$ is constant if and only if the vectors $-u$ and $-v$
are asymptotic, i.\ e.\ if $W^u(v)=W^u(u)$. In this case
$$b_v(p(\phi_t(u)))=b_v(\gamma_u(t))=b_v(p(u))+t.$$
For any $v\in UX$, the
gradient field $\grad b_v:X\rightarrow UX$ is a $C^1$-dif\-feo\-morph\-ism onto
$W^u(v)$. Classically $H_v:=b_v\inv(0)$ is called the
\emph{horosphere}{}\index{horosphere!classical terminology} centered at $\gamma_v(-\infty)$
through $p(v)$ or the horosphere determined by $v$.  However we prefer to
think of the \emph{horosphere}\index{horosphere!in $UX$} as an object in $UX$ and write
$\horo_v:=\grad b_v(b_v\inv(0))$;
  \index{h v@$\horo_v$} \index{$h v$@$\horo_v$}
traditionally the foliation of $UX$ into these
horospheres in $UX$ is called the \emph{strong
unstable foliation}. If we need to distinguish we will refer to
$H_v:=p(\horo_v)=b_v\inv(0)$
  \index{H v@$\horo_v$}  \index{$h v$@$\horo_v$}
as the \emph{horosphere in
  $X$}\index{horosphere!in $X$}.  The
\emph{strong unstable}\index{foliation!strong unstable}\index{strong unstable foliation} and
\emph{strong stable foliation}\index{foliation!strong stable}\index{strong stable foliation} are
continuous foliations of the unstable and stable foliation, respectively. The leaves
are $(n-1)$-dimensional $C^1$-manifolds in $UX$. Notice the behaviour under
the geodesic flow
$$\phi_t(\horo_v)=\horo_{\phi_t(v)}.$$
\index{geodesic flow}

\paragraph{Parallel Vectors}{}$ $\\
%\subsection{Parallel Vectors}
%\neuesubsection{Parallel Vectors}
Among the asymptotic vectors the
parallel vectors play a special r\^ole.  Two vectors $v,u\in UX$ are said
to be \emph{parallel}\index{vector!parallel}\index{parallel vectors}, if they are asymptotic and
$-u$ and
$-v$ are asymptotic, too.\\
\\
\begin{lemm}{}{lemma1}{}
  For any $u,v\in UX$ the conditions
\begin{enumerate}
\item $u$ and $v$ are parallel vectors
\item $d(\gamma_u(t),\gamma_v(t))$ is constant in \bbR\\
  (i.\ e.\ $\gamma_u$ and $\gamma_v$ are \emph{parallel
    geodesics}\index{parallel geodesics}\index{geodesics!parallel})
\item $d(\gamma_u(t),\gamma_v(t))$ is bounded in \bbR
\item $u\in W^u(v)\cap W^s(v)$
\item $W^u(v)=W^u(u)$ and $W^s(v)=W^s(u)$
\item $\gamma_u$ and $\gamma_v$ span a flat strip\\
  (i.\ e.\ the convex hull of the geodesics is a totally geodesic, flat
  submanifold of $X$)
\end{enumerate}
are equivalent.\\
Furthermore the following holds for parallel vectors $u, v$ with $u\neq\dot\gamma_v(\mathbbR)$
\begin{enumerate}
\item $\gamma_u$ and the geodesic joining $p(u)$ to $p(v)$ span a flat
  plane
\item the restriction of $\exp_{p(u)}$ to the subspace
  $\langle\{u,\exp_{p(u)}\inv(p(v))\}\rangle$ of $\,UX$ is an isometry onto
  the image.
\end{enumerate}
\end{lemm}
Being parallel is an equivalence relation: Suppose
$d(\gamma_v(t),\gamma_{u_i}(t))=\delta_i$ for all $t\in\mathbb R$. Then
$d(\gamma_{u_1}(t),\gamma_{u_2}(t))\leq\delta_1+\delta_2$ is bounded and
hence
constant, too.\\

We write $\bars v$ \index{v@$\bars v$}\index{$v $@$\bars v$} for the
equivalence class of $v$, 
i.\ e.\ the set of all vectors parallel to $v$. 
$$\bars v:=W^u(v)\cap W^s(v)$$

\paragraph{Stable and Unstable Jacobi Fields}$ $\\
A Jacobi vector field $J$ is called \emph{stable}, respectively
\emph{unstable},\index{unstable Jacobi field}\index{stable Jacobi
  field}\index{Jacobi field!stable}\index{Jacobi field!unstable} 
if $\|J(t)\|$ is bounded for $t>0$, respectively $t<0$. The Jacobi field
$J$ is called \emph{parallel}\index{parallel Jacobi field}\index{Jacobi
  field!parallel} if $\|J(t)\|$ is bounded for all $t\in\mathbbR$. For every geodesic
$\gamma_v$ and any vector $w\in T_{p(v)}X$ there is a unique stable Jacobi
field $J^s$ along $\gamma_v$ with $J^s(0)=w$ and a unique stable Jacobi
field $J^u$ along $\gamma_v$ with $J^u(0)=w$. A Jacobi field is obviously
parallel, if and only if it is stable and unstable.\\
As the name suggests, stable and unstable Jacobi fields are closely related to the stable and
unstable foliation. There is a canonical identification $\xi\mapsto J_\xi$
of elements of $TTX$ with Jacobi fields. We will describe this in
Section~\ref{infflats.tex} in more detail. It is then easy to show that
$J_\xi$ is a stable Jacobi field along $\gamma_v$ if and only if $\xi$ is an element of the
tangent bundle $TW^s(v)$ of the stable leaf of $v$. An analogue result
holds for unstable and parallel Jacobi fields.

\mysection{Rank and Flatness}{}\label{vertical.tex2}

The rank of a geodesic is the maximal number of linearly independent
parallel Jacobi fields along that geodesic. In this section we will
identify Jacobi fields with 
elements of $TTX$. This way it is easy to see that the rank is
semicontinous. If we think of parallel Jacobi fields as `infinitesimal
flats' the finite equivalent is the flatness, describing how many
independent flats there are along a geodesic.\\
As will be shown in Subsection~\ref{hyperbolic.tex} rank one vectors show
hyperbolic behaviour. In Subsection~\ref{flats.tex} we will see that in the
real analytic setting higher rank vectors might be integrated to find
flats.\\

Keep in mind that we will always work with complete manifolds that are smooth
or even real analytic.

\mysubsection{Horizontal and Vertical Structure}{}\label{vertical.tex}
If $X$ is any $C^\infty$-manifold, then the tangent bundle $TX$ and the
unit tangent bundle $UX$ are $C^\infty$-manifolds in a natural way.
Consider $UX$ as
a submanifold of $TX$.\\
We have the base point projections $p:TX\longrightarrow X$ and
$p':TTX\longrightarrow TX$\index{p'@$p':TTX\longrightarrow
  TX$}\index{$p':TTX\longrightarrow TX$}. Consider the differential
$dp$\index{d p@$dp$}\index{$d p$} of the
projection (Defined by $dp(\xi):=\frac{d}{dt}\big|_{t=0}(p\circ\sigma(t))$
where $\sigma$ is a curve in $TX$ with $\dot\sigma(0)=\xi$).
$$
\begin{CD}
  TTX@>{dp}>>TX\\ @V{p'}VV@VVpV\\ TX@>>{p}>X
\end{CD}
$$

For any $v\in TX$ we get a linear map $d_vp:T_vTX\rightarrow T_{p(v)}X$
with $n$-dimensional kernel $V_v:=\ker d_vp=T_vT_{p(v)}X\subset T_uTX$. We
can restrict $p'$ to the union of these $V_v$\index{V v@$V_v$}\index{$v v$@$V_v$} to get the
\emph{vertical bundle}\index{bundle!vertical}\index{vertical bundle} $V$ over $TX$. If we have a
Riemannian metric $g$ on $X$, we can define a complementary bundle $H$, the
\emph{horizontal bundle}\index{bundle!horizontal}\index{horizontal bundle} (cf. \cite[p. 53
II,4.1]{Sakai1996}):\\ For any vector $v\in TX$ and any geodesic
$c:I\rightarrow X$ with $c(0)=p(v)$ let $_c\|{^t_0}v$ describe the vector
in $T_{c(t)}X$ one gets by parallel transport of $v$ along $c|_{[0,t]}$.
Define the map
\begin{align*}
  h_v^{-1}:T_{p(v)}X&\longrightarrow T_vTX\\ u&\longmapsto
  \frac{\partial}{\partial t}\Big|_{t=0}(_{\gamma_u}\|^t_0v)
\end{align*}
which is an isomorphism onto the image $H_v:=\Im(h_v^{-1})$. Thus we have
$$
  H_v=\Big\{\dot X(0)\;\Big|\;\begin{array}{c}X:I\rightarrow TX \text{ parallel along a geodesic
    }\\\gamma:I\rightarrow X \text{ with } \gamma(0)=p(v) \text{ and
    }X(0)=v\end{array}\Big\}\index{H v@$H_v$}\index{$h v$@$H_v$}
$$
and in the case of a Hadamard manifold we can identify the elements of $H$ with the parallel vector fields on
$X$.
This is a one-to-one relation.\\

It is an easy calculation to show that $dp\circ
h_u=\id\big|_{T_{p(u)}X}$. Hence $T_uTX=H_u\oplus V_u$ and
$(dp\big|_{H_u})\inv=h_u\inv$. We interpret $dp$ as projection onto the
horizontal bundle and therefore write $h$ for $dp$. The projection onto the 
vertical bundle is the natural homomorphism
$v_u:V_u=T_uT_{p(u)}X\rightarrow T_{p(u)}X$ defined by
$v_u\inv(w):=\frac{d}{dt}\big|_{t=0}(u+tw)$ which is extended to the whole
$TTX$ by $v\big|_{H}=0$.\\
Another way to define $v(\xi)$ for $\xi\in TTX$ is by
$$v(\xi):=\frac\nabla{dt}(\gamma'_\xi(t))$$
where $\gamma'_\xi$ is a curve%
\footnote{We will see that $TX$ is a Riemannian manifold, hence we could
  take a geodesic in $TX$ here.}
in $TX$ with $\dot\gamma'_\xi(0)=\xi$, and $\frac\nabla{dt}$ is the
covariant derivative along the curve $p\circ\gamma'_\xi$.\\

If we restrict ourselves to $UX$ we can define the vertical and horizontal
bundle in an analogous way. The resulting vertical bundle is an
$(n-1)$-dimensional vector bundle over $UX$, the resulting horizontal
bundle is
the restriction of $H$ to $\underset{v\in UX}{\bigcup}H_v$.\\

Using the projections $h$ and $v$ we can define a Riemannian metric, the
\emph{Sasaki metric}\index{metric!Sasaki}\index{Sasaki metric} on $TX$ 
via the scalar product which is determined by the facts that $H$ and $V$ are 
orthogonal bundles and the projections $h_u\big|_{H_u}$ and $v_u\big|_{V_u}$ are isomorphisms: 
$$\langle\xi,\zeta\rangle:=\langle h(\xi),h(\zeta)\rangle+\langle
v(\xi),v(\zeta)\rangle\index{$<.,.>$@$\langle.,.\rangle$}$$
The resulting metric on $TX$ or $UX$ respectively will be denoted by
$d(.,.)$\index{d(.,.)@$d(.,.)$}\index{$d(.,.)$}.\\           
Consider a geodesic $\gamma$ in $X$. $\dot\gamma$ is a parallel vector
field along $\gamma$. $\dot\gamma$ is hence a geodesic in $TX$ all of whose 
tangent vectors are horizontal. As a consequence 
$d(\dot\gamma(0),\dot\gamma(1))=d(\gamma(0),\gamma(1))$ holds for any shortest geodesic and for any
two vectors $u,v\in TX$ we have $d(u,v)\geq d(p(u),p(v))$.\\
On the other hand, any curve in $U_{p(u)}X$ has only vertical tangent
vectors. It may therefore be considered as a curve in $S^{n-1}$. We
conclude that $d(u,v)\leq \pi$ if $u$ and $v$ have common base
point. Furthermore in this case $d(u,v)$ is exactly the angle between the
geodesics $\gamma_u$ and $\gamma_v$.\\
We sum up these results for the Riemannian metric $d$ on $UX$:

\begin{prop}{}{keinbeweis}{}
If $\gamma\big|_{[0,t]}$ is any shortest geodesic in $X$ and $u\in
U_{\gamma(0)}X$ any unit vector, then
$$d(u,\dot\gamma(t))\leq
d(\gamma(0),\gamma(t))+\measuredangle\left(\dot\gamma(0),u\right)$$
More generally for parallelly transported vectors we get ($w\in U_{\gamma(t)}X$)
$$d(u,w)\leq d(\gamma(0),\gamma(t))+d(_\gamma\big\|_0^tu,w)=d(p(u),p(w))+\measuredangle(_\gamma\big\|_0^tu,w).$$
\end{prop}
If the connecting geodesic $\gamma$ between any two points is unique
(e.\ g.\ if $X$ is a Hadamard space), we can define another `metric' on $UX$ by the second estimation:
$$d^1(u,w):=d(p(u),p(w))+\measuredangle(_\gamma\big\|_0^tu,w)
\index{d 1(.,.)@$d^1(.,.)$}\index{$d 1(.,.)$@$d^1(.,.)$}$$
In this case, $d^1$ is continuous as can be seen by rewriting
$$\measuredangle(_\gamma\big\|_0^tu,w)=\arccos\langle_{\gamma_{\exp\inv(p(w))}}\|_0^1u,w\rangle.$$

$d^1$ is an upper estimation of $d$. It is not a metric itself, since it
needs not satisfy the triangle inequality. Still it is an appropriate means
for estimating local distances, since the sets 
$${\cal W}^1_\epsilon(w):=\{u\in UX\;|\; d^1(u,w)<\epsilon\}$$
constitute a basis of the topology defined by $d$. I.e.
$$
\underset{u\in UX}\forall\quad
\underset{\delta}\forall\quad
\underset{\epsilon}\exists\quad {\cal W}_\epsilon(u)\subset {\cal
  W}^1_\delta(u)\subset{\cal W}_\delta(u) 
$$
where ${\cal W}_\delta(u)$ denotes the $\delta$-neighbourhood of $u$ in
$UX$ with respect to the metric $d$.
If $X$ has a compact quotient, $\epsilon$ may be chosen globally
$$
\underset{\delta}\forall\quad
\underset{\epsilon}\exists\quad
\underset{u\in UX}\forall
\quad {\cal W}_\epsilon(u)\subset {\cal
  W}^1_\delta(u)\subset{\cal W}_\delta(u). 
$$
Thus we may work with $d^1$ instead of $d$ whenever this is more convenient.

\begin{lemm}{}{equidist}{}
Suppose $X$ is a Riemannian manifold diffeomorphic to $\mathbbR^n$ which
has a
compact quotient. Then there is a constant $\nu\in]0,1]$ and a
neighbourhood $\cal N$ of the diagonal in $X\times X$ such that
$$(u,v)\in{\cal N}\qquad\Longrightarrow\qquad \nu
d^1(u,v)<d(u,v)<d^1(u,v).$$
\end{lemm}

\mysubsection{Flats and Jacobi Fields}{}\label{infflats.tex}
An important task when considering manifolds of nonpositive curvature is to
find behaviour that is similar to the Euclidean or hyperbolic case. In this
section we want to motivate the idea that locally flat behaviour can be
assigned to parallel Jacobi fields. We introduce flatness and rank of a
geodesic to measure the `flat' behaviour in the vicinity of the geodesic.
\subsection*{Flat Strips}
A \emph{flat}\index{flat} or a \emph{plane}\index{plane} in a manifold $M$ is an isometric
embedding of Euclidean $\mathbbR^2$ into the manifold.
A \emph{flat strip}\index{strip!flat}\index{flat strip} is the isometric embedding of the
open subset $]a,b[\times\mathbbR$ of Euclidean plane $\mathbbR^2$ into a
manifold $X$. 
Obviously, whenever we have an embedded flat or flat strip we have some
geodesics that behave like geodesics in Euclidean space. To be more
precise, a flat strip is in fact a geodesic variation
$h:I\times\mathbbR\rightarrow M$ consisting of parallel geodesics
$h_s$. Here parallel means that for any two $s_1,s_2$ the distance
$d(h_{s_1}(t),h_{s_2}(t))$ is constant in $t$.\\
In this case the variational vector field $X$ of $h$ is a parallel Jacobi field
along $h_0$, it hence satisfies
\begin{align}
X'&=0\tag{parallel}\\
R(X,\dot\gamma)\dot\gamma&=0\tag{Jacobi equation}
\end{align}
Notice that in the real analytic case it suffices if these differential
equations are satisfied for a small intervall. Furthermore, in the real
analytic case, every flat strip is part of a flat.\\
So parallel Jacobi fields may arise as variational vector fields of
geodesic variations. There are, however, parallel Jacobi fields which are
not variational vector fields of flat strips. We consider these
\emph{infinitesimal flat strips}.\index{infinitesimal flat
  strip}\index{flat strip!infinitesimal}\index{strip!infinitesimal flat}\\

We will define for every geodesic a number that measures how much flatness
and infinitesimal flatness we encounter in its vicinity. To this end we identify Jacobi fields with
element of the tangent bundle $TUM$ of $UM$. 

\subsection*{Jacobi Fields and $TTM$}
There is a natural identification of elements $\xi\in TTM$ with a
Jacobi field $J_\xi\in\frak J(\gamma_{p'(\xi)})$ satisfying $J_\xi(0)=h(\xi)$
and $J'_\xi(0)=v(\xi)$ where $h$ and $v$ are the horizontal respectively
vertical projection $T_vUM\rightarrow T_{p(v)}M$. The map
$J:TTM\rightarrow\frak J$ is a bijection onto the set of all Jacobi fields
along geodesics in $M$. Every vector in $\xi\in TTM$ is uniquely identified
by a point $p:=p\circ p'(\xi)\in M$, and the three vectors
$p'(\xi),h(\xi),v(\xi)\in T_pM$. In fact the following submanifolds of
$TTM$ can be defined by the respective equations:
$$
\begin{array}{|c|c|c|c|c|}
\hline
\text{Manifold}&\multicolumn{3}{c|}{\text{Conditions}}&\text{Dimension}\\\hline
TTM&&&&4n\\\hline
T_uTM&p'(\xi)=u&&&2n\\\hline
\underset{u\in UM}\bigcup T_uTM&\|p'(\xi)\|=1&&&4n-1\\\hline
TUM&\|p'(\xi)\|=1&v(\xi)\perp p'(\xi)&&4n-2\\\hline
T_uUM&p'(\xi)=u&v(\xi)\perp p'(\xi)&&2n-1\\\hline
H&&v(\xi)=0&&3n\\\hline
H_u&p'(\xi)=u&v(\xi)=0&&n\\\hline
\underset{u\in UM}\bigcup H_u&\|p'(\xi)\|=1&v(\xi)=0&&3n-1\\\hline
V&&&h(\xi)=0&3n\\\hline
V_u&p'(\xi)=u&&h(\xi)=0&n\\\hline
\end{array}$$
Notice that $v(\xi)=0$ implies $v(\xi)\perp p'(\xi)$ and hence $H_u\subset
T_uUM$ hold for $u\in UM$s. Now consider the image of these submanifolds
of $TTM$ under the identification with Jacobi fields we just
explained. Clearly $\|p'(\xi)\|=1$ means that $J_\xi$ is a Jacobi field
along a unit speed geodesic. $p'(\xi)=u$ means that the geodesic is
$\gamma_u$. If $v(\xi)\perp p'(\xi)$, then $J'_\xi\perp \gamma_{p'(\xi)}$
for all times.%
\footnote{Here we use that $\langle J',\dot\gamma\rangle$ is
  constant for any Jacobi field $J$ along a geodesic $\gamma$.} 
Furthermore $v(\xi)=0$ means $J'_\xi(0)=0$ and $h(\xi)=0$ means $J_\xi(0)=0$.\\
We will only consider unit speed geodesics $\gamma_u$ ($u\in UM$). The unit
speed geodesic variations of $\gamma_u$ are represented by $H_u$.

\subsection*{Flatness and Rank}
\begin{figure}[h] \caption{Definition of $\frak f$}\label{deffrakf}
\psfrag{TpX}{$T_{p}X$}
\psfrag{X}{$X$}
\psfrag{p}{$p$}
\psfrag{v}{$v$}
\psfrag{hxi}{$h(\xi)$}
\psfrag{thxi}{$th(\xi)$}
\psfrag{tv}{$tv$}
\psfrag{ghxit}{$\gamma_{h(\xi)}(t)$}
\psfrag{gvt}{$\gamma_v(t)$}
\psfrag{gv}{$\gamma_v$}
\psfrag{ghxi}{$\gamma_{h(\xi)}$}
\psfrag{text1}{$v=p'(\xi)$}
\psfrag{text2}{$p=p(v)$}
\includegraphics[width=\textwidth]{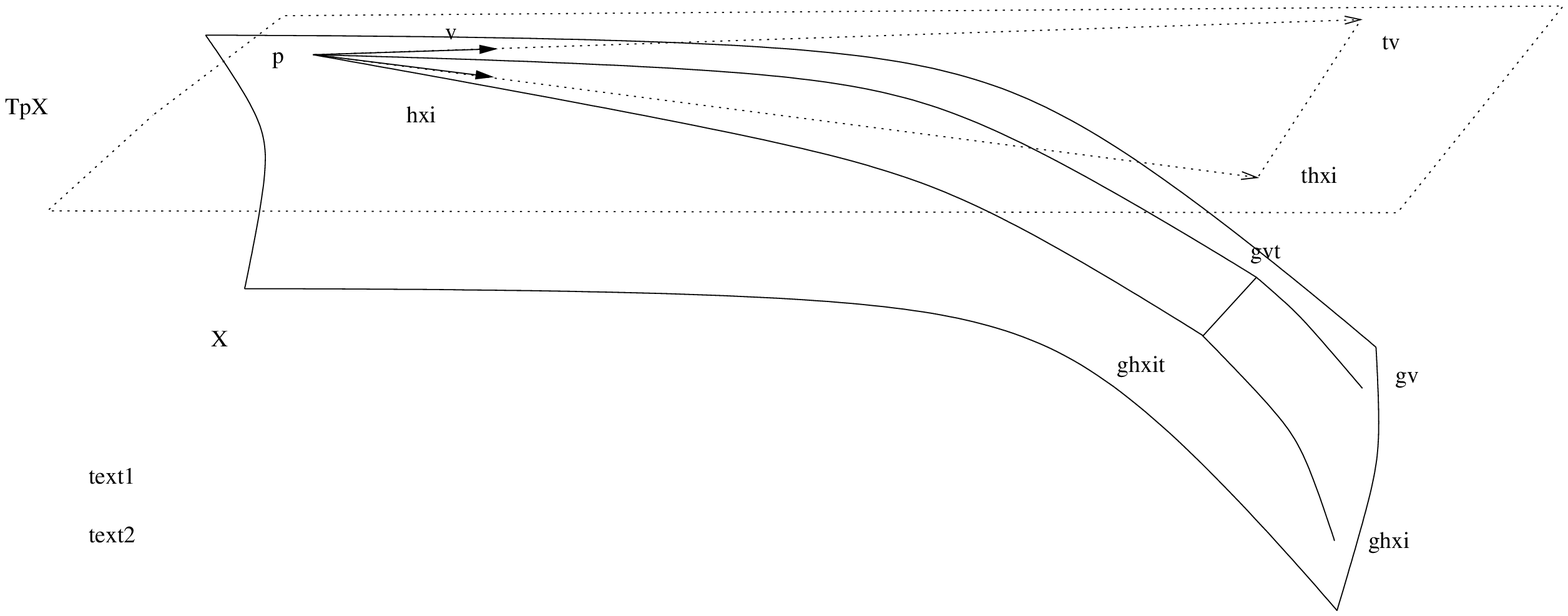}
\end{figure}

For every unit speed vector we want to define the rank and the flatness. 
The flatness is the dimension of vectors parallel to $\gamma_u$. The rank
is the dimension of parallel Jacobi fields along $\gamma_u$. To see that
both these values are semicontinuous we need to work on $\tilde
H:=\underset{u\in UM}\bigcup H_u$.\\
We define two smooth functions to measure flatness and infinitesimal
flatness:
\index{f(v)@$\frak f(v)$}
\index{$f(v)$@$\frak f(v)$}
\index{$j(v)$@$\frak j(v)$}
\index{j(v)@$\frak j(v)$}
$$\begin{array}{cccc}\label{Defj}\label{Deff}
\frak f:&\mathbbR\times\tilde H&\longrightarrow&\mathbbR\\
  &(t,\xi)&\longmapsto&d(\gamma_{p'(\xi)}(t),\gamma_{h(\xi)}(t))^2-t^2\|p'(\xi)-h(\xi)\|^2\\\text{ }\\
\frak j:&\mathbbR\times\tilde H&\longrightarrow&\mathbbR\\
   &(t,\xi)&\longmapsto&\|J_\xi(t)\|^2-\|J_\xi(0)\|^2
\end{array}$$
$\frak f$ measures the difference between corresponding sides of geodesic triangles
in $T_{p(v)}X$ and $X$ where the triangle in $X$ is the image of the
triangle in $T_{p(v)}X$ under the exponential map. This is illustrated in Figure~\ref{deffrakf}. $\frak j$ measures
whether the length of the Jacobi field defined by $\xi\in H_v$ is constant
along $\gamma_v$.\\

We will need the following properties of $\frak f$ and $\frak j$ on nonpositively
curved manifolds.
 
\begin{lemm}{}{frakf}{}
Suppose $X$ is a Hadamard manifold. Then

\begin{enumerate}
\item $\frak f\geq 0$ and $\frak f(t,\xi)=0$ for $t=0$ or if $p'(\xi)$ and $h(\xi)$ are
  collinear.
\item For $0< s\leq t$ the implication $\frak f(t,\xi)=0\Longrightarrow
  \frak f(s,\xi)=0$ holds.
\item For $t\leq s< 0$ the implication $\frak f(t,\xi)=0\Longrightarrow
  \frak f(s,\xi)=0$ holds.
\item The triangle $\Delta(p\circ
  p'(\xi),\gamma_{p'(\xi)}(t),\gamma_{h(\xi)}(t))$ is
  flat\footnote{isometric to a triangle in $\mathbbR^2$} if and only if
  $\frak f(t,\xi)=0$. This triangle is degenerated only if $p'(\xi)$ and $h(\xi)$
  are collinear or if $t=0$. 
\end{enumerate}
Furthermore if $X$ is analytic, then so it $\frak f$ and
\begin{enumerate}\setcounter{enumi}{4}
\item Any nondegenerated flat triangle is contained in a flat. Hence
  $\exp:\spann(p'(\xi),h(\xi))\rightarrow X$ is an isometry if and only if
  $p'(\xi)$ and $h(\xi)$ are not collinear and $\frak f(t,\xi)=0$ for some $t\neq
  0$.
\end{enumerate}
\end{lemm}

\begin{lemm}{}{frakj}{}
Suppose $M$ is a manifold of nonpositive curvature. Then
\begin{enumerate}
\item $\frak j(.,\xi)$ is a convex nonnegative function. $\frak j(t,\xi)=0$ for $t=0$ and
  whenever $p'(\xi)$ and $h(\xi)$ are collinear.
\item For $0< s\leq t$ the implication $\frak j(t,\xi)=0\Longrightarrow
  \frak j(s,\xi)=0$ holds.
\item For $t\leq s< 0$ the implication $\frak j(t,\xi)=0\Longrightarrow
  \frak j(s,\xi)=0$ holds.
\item If $\frak j(t,\xi)=0$ then $J_\xi\big|_{[0,t]}$ is a parallel Jacobi field.
\end{enumerate}
If $M$ is an analytic manifold then $\frak j$ is analytic, too, and 
\begin{enumerate}\setcounter{enumi}4
\item $J_\xi$ is a parallel Jacobi field if $\frak j(t,\xi)=0$ for some $t\neq
  0$.
\end{enumerate}
\end{lemm}

\begin{lemmp}{frakf}{}
Compare for example \cite[Prop. 1.4.1]{Eberlein1996} to see that $\frak
f\geq 0$. From the proof it can be seen that $\frak f(t,\xi)=0$ if and only
if $\frak f(s,\xi)=0$ for all $0\leq s\leq t$ and in this case the
exponential map is an isometry of the triangle
$\Delta(0,tp'(\xi),th(\xi))\subset T_{p\circ p'(\xi)}X$ onto the triangle
$\Delta(p\circ p'(\xi),\gamma_{p'(\xi)}(t),\gamma_{h(\xi)}(t))\subset
X$ which is a totally geodesic flat submanifold of $X$. In the real
analytic case this flat triangle must be contained in a flat, namely the
image under the exponential map of $\spann(p'(\xi),h(\xi))$. This is the
convex hull of the two geodesics $\gamma_{p'(\xi)}$ and $\gamma_{h(\xi)}$.
\end{lemmp} 

\begin{lemmp}{frakj}{}
Fix $\xi$ and write $\frak j$ for $t\rightarrow\frak j_\xi(t)$ and $J$ for
$J_\xi$. An easy calculation shows that
\begin{align*} 
\frak j'&=2\langle J',J\rangle\\
\frak j''&=2\|J'\|^2-2\langle R(J,\dot\gamma)\dot\gamma,J\rangle.
\end{align*}
Recall that $\langle
R(J,\dot\gamma)\dot\gamma,J\rangle=K(J,\dot\gamma)(\|J\|^2\|\dot\gamma\|^2-\langle
J,\dot\gamma\rangle^2)$ and therefore this term has the same sign as the curvature,
namely it is nonpositive and therefore $\frak j$ is convex, $\frak j(0)=0$ and $\frak j'(0)=0$ since $J'(0)=0$.
We conclude that $\frak j(t)\geq0$ and $\frak j'(t)\geq0$ for
$t\geq0$. Now suppose $\frak j(t)=0$ then $\frak j\big|_{[0,t]}\equiv 0$ and
$\frak j'\big|_{[0,t]}\equiv 0$. This implies $\frak j''\big|_{[0,t]}\equiv 0$
and hence $\|J'\|^2\equiv\langle R(J,\dot\gamma)\dot\gamma,J\rangle\leq
0$. We conclude that in this case $\|J'\|\big|_{[0,t]}\equiv 0$ which means
that $J\big|_{[0,t]}$ is parallel. An analogous argument works for $t<0$.
\end{lemmp}

\begin{lemm}{}{UVRe}{}
\begin{enumerate}

\item If $X$ is a Hadamard manifold then the set
$$\operatorname{Flat}_v:=\{\xi\in H_v\;|\;\frak f(t,\xi)=0\text{ for all
  }t\in\mathbbR\}$$
\index{Flat v@$\operatorname{Flat}_v$}
\index{$flat v$@$\operatorname{Flat}_v$}
is a subspace of $H_v$ for all $v\in UX$. Therefore the
\emph{flatness}\index{flatness}\index{$flat(v)$@$\flat(v)$}\index{flat(v)@$\flat(v)$}, defined by
\begin{align*}
\flat\;:\;UX&\longrightarrow \{1,2,\ldots ,n\}\\
v&\longrightarrow\dim\{\xi\in H_v\;|\;\frak f(t,\xi)=0\text{ for all
  }t\in\mathbbR\}
\end{align*}
is semicontinuous, i.\ e.\ $\lim\flat(v_n)\leq\flat(v)$ if $v_n\rightarrow v$.

\item
 If $M$ is a manifold of nonpositive curvature then the set 
$$\operatorname{Rank}_v:=\{\xi\in H_v\;|\;\frak j(t,\xi)=0\text{ for all
  }t\in\mathbbR\}$$
\index{Rank v@$\operatorname{Rank}_v$}
\index{$rank v$@$\operatorname{Rank}_v$}
is a subspace of $H_v$ for all $v\in UM$. Therefore the
\emph{rank}\index{rank}\index{$rank(v)$@$\rank
  j(v)$}\index{rank(v)@$\rank(v)$}, defined by
\begin{align*}
\rank\;:\;UM&\longrightarrow \{1,2,\ldots ,n\}\\
v&\longrightarrow\dim\{\xi\in H_v\;|\;\frak j(t,\xi)=0\text{ for all
  }t\in\mathbbR\}
\end{align*}
is semicontinuous, i.\ e.\ $\lim\rank(v_n)\leq\rank(v)$ if $v_n\rightarrow v$.
\end{enumerate}
\end{lemm}

\begin{lemmp}{UVRe}{}
\begin{enumerate}
\item Suppose $\xi\in \operatorname{Flat}_v$ and $\lambda\in \mathbbR$ is given. 
Since $\frak f(t,\xi)=0$ for all $t\in\mathbbR$, the geodesics
  $\gamma_v$ and $\gamma_{h(\xi)}$ span a totally geodesic flat. Now given
  $t,\lambda\in\mathbbR$ the points $p(v)$,
  $\exp(t\,p'(\lambda\xi))=\gamma_{p'(\lambda\xi)}(t)=\gamma_v(t)$ and
  $\exp(t\,h(\lambda\xi))=\gamma_{h(\xi)}(\lambda t)$ lie within this flat
  and hence form a flat geodesic triangle. Thus $\frak
  f(t,\lambda\xi)=0$ and hence $\lambda\xi\in\operatorname{Flat}_v$.\\

  Now fix  $\xi,\zeta\in\operatorname{Flat}_v$. Then the geodesics
  $\gamma_v$ and $\gamma_{h(\xi)}$ span a totally geodesic flat. This means
  that the geodesics $\gamma_1:s\rightarrow \exp(tJ_\xi(s))$ and
  $\gamma_2:s\rightarrow \exp(tJ_\zeta(s))$ are parallel to
  $\gamma_v$. Consider the distance between these two geodesics. This is a
  convex function which is bounded by
  $d(\gamma_1,\gamma_v)+d(\gamma_v,\gamma_2)$ and hence constant. So
  $\gamma_1$ and $\gamma_2$ are parallel, too, and span a flat plane. We
  find a third geodesic
  $\gamma_3:s\rightarrow\exp(\sqrt2tJ_{\xi+\zeta}(s))$ which is parallel to
  $\gamma_1$ and $\gamma_2$ and hence to $\gamma_v$, too. We conclude that
  $\frak f(t,\sqrt2(\xi+\zeta))=0$ and from this we easily deduce that
  $\xi+\zeta\in \operatorname{Flat}_v$.\\ 

\item
  For $\frak j$ notice that $\frak j(t,\lambda\xi)=\lambda^2\frak
  j(t,\xi)$. Therefore $\xi\in\operatorname{Flat}_v$ implies $\lambda\xi\in
  \operatorname{Flat}_v$ for all $\lambda\in\mathbbR$.\\
  
  Now notice that $J'_\xi\equiv 0$ for $\xi\in\operatorname{Rank}_v$. Now
  suppose 
$\xi,\zeta\in\operatorname{Rank}_v$ and consider $\frak j:t\rightarrow
\frak j(t,\xi+\zeta)$. A simple calculation shows that $\frak j(t)=2\langle
J_\xi(t),J_\zeta(t)\rangle-2\langle J_\xi(0),J_\zeta(0)\rangle$ and $\frak
j'(t)=2\langle J_\xi'(t),J_\zeta(t)\rangle+2\langle
J_\zeta'(t),J_\xi(t)\rangle\equiv 0$ and hence
$\xi+\zeta\in\operatorname{Rank}_v$. 

\item By now we know that $\frak f$ and $\frak j$ are continuous and
  $\operatorname{Flat}_v$ and $\operatorname{Rank}_v$ are subspaces of
  $H_v$ for all $v\in UM$. We will prove that $\rank$ is semicontinuous
  (the proof for $\flat$ works analogously).\\ 
  Fix $v\in UM$ and write $r$ for $\rank(v)$. Choose an orthonormal basis
  $\{\eta_1,\ldots ,\eta_n\}$ of $H_v$ such that $\{\eta_1,\ldots
  ,\eta_r\}$ is a basis of $\operatorname{Rank}_v$. We have therefore
  \begin{align*}
  \langle \eta_i,\eta_j\rangle=\delta_{i,j}\qquad&\qquad\text{for all
    }1\leq i\leq j\leq n,\\
  \frak j(t,\eta_i)=0\qquad&\qquad\text{for all }t\in\mathbbR\text{ and
    }1\leq i\leq r\\
  \frak j(t_i,\eta_i)=a_i>0\qquad&\qquad\text{for some
    }t_i\in\mathbbR,\text{ and }r<i\leq n.
  \end{align*}
  Now we can find neighbourhoods $\W_i$ of the $\eta_i$ in $\underset{u\in
    UM}H_u$ where softer versions of the above inequalities hold, namely
  for $\xi_i\in\W_i$ and $\xi_j\in \W_j$
  \begin{align*}
  |\langle
  \xi_i,\xi_j\rangle-\delta_{i,j}|<\frac{1}{n+1}\qquad&\qquad\text{for all
    }1\leq i\leq j\leq n\\
  \frak j(t_i,\xi_i)=\frac{a_i}2>0\qquad&\qquad\text{for }r<i\leq n
  \end{align*}
  $\W:=\bigcap p'(\W_i)$ is an open neighbourhood of $v$ and for any
  $u\in\W$ we may pick arbitrary elements $\xi_i\in \W_i\cap {p'}\inv(u)$,
  which will form a basis of $H_u$. Since $\frak
  j(t_i,\xi_i)=\frac{a_i}2>0$ for $r<i\leq n$ the dimension of
  $\operatorname{Rank}_u$ can be at most $r$.
\end{enumerate}

\end{lemmp}

\mysection{Sets of Constant Rank}{Constant Rank}\label{subanalytic.tex}
Consider the sets of vectors of constant rank or of bounded rank. Namely
for any manifold $X$ of nonpositive curvature define for $1\leq
k\leq n$

\begin{align*}
  \R_k&:=\{u\in UX\;|\;\rank(u)= k\}\\
  \R_{>k}&:=\{u\in UX\;|\;\rank(u)>k\}\\
  \R_{<k}&:=\{u\in UX\;|\;\rank(u)<k\}\\
\end{align*}
\index{R 1@$\R_1$}
\index{R k@$\R_k$}
\index{R >k@$\R_{>k}$}
\index{R <k@$\R_{<k}$}%
\index{$r 1$@$\R_1$}
\index{$r k$@$\R_k$}
\index{$r>k$@$\R_{>k}$}
\index{$r<k$@$\R_{<k}$}%
By Lemma~\ref{UVRe}, we know that the rank is semicontinuous and hence $\R_{<k}$ and $\R_1$ are
open subsets of $UX$ and $\R_{>k}$, $\R_n$ are closed. We will be
interested in the complementary sets $\R_1$ and $\R_>:=\R_{>1}$\index{R_>@$\R_>$} mainly. \\

In Subsection~\ref{CRAM} we will see that all these sets are subanalytic if $X$ is a real
analytic manifold with compact quotient. We will then define an equivalent
for smooth manifolds. But first we need to understand semi- and subanalytic
sets.

\mysubsection{Subanalytic Sets}{}
%\subsection{Subanalytic Sets}
%\neuesubsection{Subanalytic Sets}
The following definitions and propositions are taken
from~\cite{Hironaka1973} without proof. Suppose $X$ is any real analytic manifold.

\begin{defi}{\protect{\cite[Def. 2.1]{Hironaka1973}}}
\label{semianalytic sets}\index{semianalytic subset}\index{subset!semianalytic}
A subset $A\subset X$ is called \emph{seminanalytic} if for every point
$a\in A$ there is a
finite number of real analytic functions $g_{ij}:\U_a\rightarrow\mathbbR$
on an open neighbourhood $\,\U_a$ of $a$ in $X$ such 
that 
$$A\cap\U_a=\bigcup\limits_i\bigcap\limits_jg_{ij}\inv(I_{ij})$$
where $I_{ij}$ is any interval in $\mathbbR$.
\end{defi}

If we consider the union, intersection or difference of two semianalytic
sets, the resulting set is semianalytic itself
(cf.~\cite[Rem. 2.2]{Hironaka1973}). However, the image of a semianalytic
set under an analytic map needs not be semianalytic. We have therefore to
consider the bigger class of subanalytic sets.

\begin{defi}{\protect{\cite[Def. 3.1]{Hironaka1973}}}
A subset $A$ of $X$ is called \emph{subanalytic},
\index{subanalytic subset}\index{subset!subanalytic}
if for every point $a\in A$ there are finitely many
\emph{proper}\index{proper map}\index{map!proper}%
\footnote{A map is called proper if the preimage of any compact set is
  compact}
real analytic maps $f_j:Y_j\rightarrow \U_a$ and $h_j:Z_j\rightarrow \U_a$,
where $\,\U_a$ is an open neighbourhood of $a$ in $X$, such
that
$$A\cap\U_a=\bigcup\limits_j(\Im(f_j)\backslash\Im(h_j)).$$
\end{defi}

The union, intersection or difference of two subanalytic sets is again
subanalytic~\cite[Prop.3.2]{Hironaka1973}. Furthermore any proper real analytic
map maps subanalytic sets on subanalytic sets. The notion of
subanalyticity is in fact a generalisation of semianalyticity, since every
semianalytic set is subanalytic, too~\cite[Prop. 3.4]{Hironaka1973}. We
will need the following result about the structure of subanalytic sets:

\begin{theo}{\protect{\cite[Main Theorem (4.8)] {Hironaka1973}}}{structure}{}
Let $A\subset X$ be a subanalytic subset of a real analytic manifold
$X$. Then $A$ admits a \emph{Whitney stratification},
\index{stratification!Whitney}\index{Whitney stratification}
 i.\ e.\ we can decompose $A=\bigcup A_\alpha$ 
such that
\begin{enumerate}
\item $\{A_\alpha\}_\alpha$ is a locally finite family of pairwise disjoint 
  subsets of $A$.\label{sub1}
\item Every $A_\alpha$ is a real analytic submanifold of $X$.
\item $\partial A_\alpha\cap A_\beta\neq\emptyset$ implies
  $A_\beta\subset\partial A_\alpha$\label{whitneyfehlt}.
\item Every $A_\alpha$ is a subanalytic subset of $X$
\item In the case of~\ref{whitneyfehlt} the pair $(A_\alpha,A_\beta)$
  satisfies the \emph{Whitney condition}
\index{Whitney condition}\index{condition!Whitney}%
\footnote{We will not use the Whitney condition and therefore we do not
  define it.}
in every point of $A_\beta$.
\end{enumerate}
\end{theo}
Obviously we can replace the Whitney stratification by a stratification in
the following sense if A is closed.

\begin{defi}{Stratifications}\label{strata}{}
Let $X$ denote a smooth manifold. A
\emph{stratification}\index{stratification} (by submanifolds) of a subset
$A\subset X$ is a collection $\{A_\alpha\}$ of subsets of $X$ such that 
\begin{enumerate}
\item $\{A_\alpha\}_\alpha$ is a locally finite family of pairwise disjoint 
  subsets of $A$.
\item Every $A_\alpha$ is a smooth closed submanifold of $X$ (possibly with
  boundary).
\end{enumerate}
A subset $A\subset X$ is called \emph{stratified}\index{stratified
  subset}\index{subset!stratified} if there exists a stratification by
submanifolds of $X$. In this case we define the
\emph{dimension}\index{dimension!of a stratified subset} and the \emph{unit
  tangent bundle}\index{unit tangent bundle!of a stratified subset} of $A$
by
\begin{align*}
\dim A &:=\sup\dim(A_\alpha)\index{dim A@$\dim A$}\index{$dim A$@$\dim A$}\\
UA&:=\bigcup\overline{ U\overset\circ A_\alpha}\index{U A@$UA$}\index{$u a$@$UA$}\\
U_aA&:=UA\cap U_aX,
\end{align*}
i.\ e.\ the dimension of $A$ is the maximal dimension of a submanifold of $X$
contained in $A$. A vector $v\in UM$ is considered tangent to $A$, if it is
tangent to a submanifold of $M$ contained in $A$, or if it is the limit of
such vectors. Another way to define the tangent bundle
\index{tangent bundle!of a stratified set} would be 
$$TA:=\left\{\dot\gamma(0)\;\left|\;
\begin{array}{l}
\gamma:]-\epsilon,\epsilon[\rightarrow A\text{
  smooth and  }\\
\gamma([0,\epsilon[)\subset A\text{ or
  }\gamma(]-\epsilon,0]\subset A.
\end{array}
\right.\right\}.$$
\end{defi}
All of these definitions are independent of the choice of stratification.

\mysubsection{Compact Real Analytic Manifolds}{}\label{CRAM}
Now consider a compact real analytic manifold $M$. By the help of the
Riemannian structure on $UM$ we can define the Stiefel bundles of the
horizontal bundle\index{horizontal bundle}\index{bundle!horizontal}
$p':H\rightarrow UM$: For all $k\in\mathbbN$ and $u\in UM$ define

$$St^k_u(H):=\{(\xi_i)_{i=1\ldots k}\in
{H_u}^k\;|\;\langle\xi_j,\xi_l\rangle=\delta_{jl}\}$$

the manifold of all \emph{orthonormal
k-frames}
\index{frame!orthonormal}\index{orthonormal frame}
in $H_u$. The {$k^{\text{th}}$ Stiefel bundle}
\index{Stiefel bundle}\index{bundle!Stiefel}
is $St^k(H):=\bigcup St^k_u(H)$ 
\index{St k(H)@$St^k(H)$}\index{$st k(H)$@$St^k(H)$}
with the obvious projection which we will denote by $p'$,
  too. If $M$ is compact, then so is $St^k(H)$ and $p':St^k(H)\rightarrow
  UM$ is proper.\\
Any analytic map $a:H\rightarrow\mathbbR$ induces an analytic map
$$a^k:St^k(H)\rightarrow\mathbbR^k$$
via
$a^k((\xi_i)_{i=1..k}):=(a(\xi_i))_{i=1..k}$.\\
Now consider the map $\frak j$ introduced in Section~\ref{Defj}. The preimage
$(\frak j^k)\inv(\{0\})$ is a semianalytic set which projects down to a
subanalytic set in $UM$ under the proper projection $p':St^k(H)\rightarrow
UM$. What does it mean if a vector $u\in UM$ is in this set? In this case
there is an orthonormal $k$-frame in $T_uM$ which consists of vectors that
can be extended to parallel Jacobi fields along $\gamma_u$. Thus this
vector is of rank at least $k$. We have proved the following
\begin{prop}{}{}{}
If $M$ is a compact real analytic manifold then the sets $\R_k$ of constant rank
and the sets $\R_{<k}$, $\R_{>k}$, $\R_>$ of bounded rank are subanalytic
sets in $UM$.\\
If $\pi:X\rightarrow M$ is the universal covering then the preimage of
these sets under $d\pi$ are subanalytic sets of $UX$ and exactly the sets
with the same conditions on the rank in $X$.\\
As a result all these sets are unions of locally finite families of
disjoint real analytic submanifolds of $\,UM$, respectively $\,UX$.
\end{prop}

Now consider the base point projections of these sets,
e.\ g.\ $R_>:=p(\R_>)$. These are subanalytic subsets of $M$ (respectively
$X$) and $\R_>\subset UR_>$. This motivates the definition of s-support and
s-dimension.

\mysubsection{s-Dimension}{}
If $X$ is a real analytic manifold then the set of higher rank vectors and
the set of its base points are subanalytic and hence stratified. Suppose
that in the smooth case, these sets are stratified, too, and of low
dimension. This is the case for which our main theorems hold.\\
In general for any subset of the unit tangent bundle of $X$ we define the
\emph{structured dimension}
\index{structured dimension}\index{dimension!structured}. 

\begin{defi}{s-Dimension}{}{}\label{def s-dimension}
For a subset $\R\subset UM$ of the unit tangent bundle of a smooth manifold
$M$ define the \emph{s-dimension}
\index{s-dimension}\index{dimension!s-dimension@s-dimension}
\index{s-dim@$\sdim$}\index{$s-dim$@$\sdim$}\index{dimension!s-dim@$\sdim$}
 of $\R$ ($\sdim(\R)$) to be the
smallest dimension of a subset $R\subset M$ of $M$ that is
stratified $R=\bigcup R_i$ by submanifolds of $M$ and satisfies $\R\subset
UR$. Any such $R$ will be called \emph{an
  s-support}
\index{s-support}\index{support@s-support}
of $\R$.\\ 
Even though there
might be different s-supports for $\R$ and different stratifications for
the same s-support, the structured dimension of $\R$ is well defined.
\end{defi}

The structured dimension is defined for any subset $\R\subset UM$ and is bounded by the dimension of $M$:
$$0\leq \sdim(\R)\leq\dim(M).$$
Obviously nontrivial flow invariant subsets of $UM$ contain at least one
geodesic and hence have structured dimension at least one. For a closed totally
geodesic submanifold $N\subset M$ the structured dimension of $UN$ is just
the dimension of $N$: $\sdim(UN)=\dim(N)$ and $N$ is an s-support of $UN$. Any full subset of $UM$
(i.\ e.\ it covers all of $M$ under the base point projection) has structured
dimension $\dim(M)$.

\begin{rema}{}{sdimflowinv}{}
If $\R$ is a flow invariant subset of $UX$ then any stratified subset of
$M$ containing $p(\R)$ is an s-support of $\R$.
\end{rema}

So the simplest example of a flow invariant subset of s-dimension one is
the set of tangent vectors to a a finite collection of simply closed
geodesic. We will come back to this example when discussing the results in
Section~\ref{rankr5.tex} and Section~\ref{topologie.tex}.

\mysection{Spheres and Horospheres}{}\label{sphaere.tex}\label{flats.tex}

In this section we will take a closer look at spheres and horospheres. To
be more precise, we are interested in the corresponding structures in $UX$,
namely the vectors orthogonal to spheres and horospheres in $X$.\\
Subsection~\ref{sphaere} contains some technical results.
We will see that horospheres can be approximated by spheres of growing
radius. We can therefore think of a horosphere as a sphere with radius
infinitely large and with centre at infinity. Furthermore we have some
technical lemma that quantifies the fact that vectors in the sphere
$\sphere_v(r)$ are close if there base points are close.\\
In Subsection~\ref{flats} we consider spheres and horospheres in a real analytic
Hadamard manifold with compact quotient. Suppose on a horosphere we can find an
open subset of vectors of higher rank. Then we can integrate the parallel
Jacobi fields to find that the whole horosphere is foliated by flats. By a
result of Werner Ballmann this can only happen if the manifold is of higher
rank:%
\paragraph*{Corollary~\ref{coro2}}$ $\\
\begin{it}
  Let $X$ be an analytic rank one Hadamard manifold with compact quotient.
  Then for any horosphere and any sphere in $UX$ the subset of rank one
  vectors is dense.
\end{it}

\mysubsection{Spheres and Limits of Sphere Segments}{Limits of Sphere
  Segments} \label{sphaere}
\begin{defi}{Spheres in $UX$}
  On a Hadamard manifold $X$ define
\begin{enumerate}
\item For any given point $o\in X$ a vector $v\in UX$ is called a
  \emph{radial vector}\index{radial vector}\index{vector!radial} with origin\index{origin} in $o$ if and only
  if
  $o\in\gamma_v(\mathbb{R}_-)$
  (I.\ e.\ if $v$ is a tangent vector to a geodesic ray originating in $o$).
\item The \emph{sphere in $UX$}\index{sphere!in $UX$} of radius $r$ and centred at
  $o$ is the set of all radial vectors with origin $o$ whose base points are
  in distance $r$ of $o$.
\item Given a vector $v\in UX$ and a positive number $r\in\mathbbR$ we
  write ${\cal S}_v(r)$\index{S v(r)@${\cal S}_v(r)$}
\index{$s v(r)$@${\cal S}_v(r)$} for the sphere in $UX$ of radius $r$
centered at 
  $\gamma_v(-r)$. Write $S_v(r)$\index{S v(r)@$S_v(r)$}\index{$s v(r)$@$S_v(r)$} for the base point set of ${\cal S}_v(r)$.
  This is the \emph{sphere in $X$}\index{sphere!in $X$} with radius $r$ and centre $\gamma_v(-r)$.
\end{enumerate}
\end{defi}

\begin{rema}{}{}{}
\begin{itemize}
\item  Notice that the sphere in $UX$ is the set of outward pointing normal
  vectors of the sphere in $X$. The definition of ${\cal S}_v(r)$ assures
  that $v\in{\cal S}_v(r)$ for all $r$.
\item Compare the behaviour of spheres in $UX$ under the geodesic flow to
  that of horospheres under the geodesic flow:
$$\phi_t(\sphere_v(r))=\sphere_{\phi_t(v)}(r+t) \qquad\qquad
  \phi_t(\horo_v)=\horo_{\phi_t(v)}$$
\end{itemize}
\end{rema}

We can understand the horosphere $\horo_v$ as the limit of $\sphere_v(r)$
as $r$ goes to infinity. This fact is well known for spheres and
horospheres in $X$. It follows directly for the spheres and horospheres in
$UX$ because they are characterised as normal vectors to the objects in
$X$.\\
The convergence of the spheres to the horosphere is uniformly on compact
subsets. To understand this we take a closer look at segments.\\
\\
For a horosphere $\horo:=\horo_v$ the base point projection $p:\horo\subset
UX\rightarrow X$ induces a metric by $d_\horo(u_1,u_2):=d(p(u_1),p(u_2))$
\index{d H(.,.)@$d_\horo(.,.)$}\index{$d H(.,.)$@$d_\horo(.,.)$}
for $u_i\in\horo$. We will write $d_\horo(.,.)$ or $d_v(.,.)$
\index{d v(.,.)@$d_v(.,.)$}\index{$d v(.,.)$@$d_v(.,.)$} for this metric. The
geodesic flow $\phi$ on $UX$ satisfies the following inequalities for
vectors $u_i\in\horo$:

$$
d_{\phi_t(\horo)}(\phi_t(u_1),\phi_t(u_2))
\begin{cases}
  \geq d_\horo(u_1,u_2)\text{ for } t>0\\
  \leq d_\horo(u_1,u_2)\text{ for } t<0
\end{cases}
$$
where equality holds if and only if $u_1$ and $u_2$ are parallel
vectors.\\
Analogously we can introduce a metric on any sphere $\sphere$ in
$UX$. Notice that $\phi_t\sphere$ is a sphere for any sphere in $UX$ and
any $t\geq 0$. Obviously
$$d_{\phi_t(\sphere)}(\phi_t(u_1),\phi_t(u_2))>d_\sphere(u_1,u_2)$$
holds for all $t>0$ and $u_i\in\sphere$.\\
Now for any $\Delta>0$ we can define (open) segments of (horo)-spheres with radius $\Delta$:
\begin{align*}
\index{H v Delta@$\horo_v^\Delta$}\index{$h v Delta$@$\horo_v^\Delta$}
\horo_v^\Delta&:=\{w\in\horo_v\;|\;d_{\horo_v}(v,w)<\Delta\}\\
\index{S v Delta@$\sphere_v^\Delta(r)$}
\index{$s v Delta$@$\sphere_v^\Delta(r)$}
\sphere_v^\Delta(r)&:=\{w\in\sphere_v(r)\;|\;d_{\sphere_v(r)}(v,w)<\Delta\}
\end{align*}
The closed segment of radius $\Delta$ will be denoted by
$\overline{\horo_v^\Delta}$ and $\overline{\sphere_v^\Delta(r)}$,
respectively.\\
\begin{lemm}{}{bigspheres}{}
Let $X$ be a Hadamard manifold. Then any sequence of
sphere segments $\sphere_v^\Delta(r)$ converges with respect to the
Hausdorff metric to the corresponding horosphere segment
$\horo_v^\Delta$.
$$\underset{\Delta>0}\forall\quad\underset{v\in
  UX}\forall\quad\underset{R}\exists\quad\underset{r>R}\forall\quad
\Hd(\horo_v^\Delta,\sphere_v^\Delta(r))<\delta$$
Furthermore if $X$ has compact quotient, $\delta$ does not depend on $v$:
$$\underset{\Delta>0}\forall\quad\underset{R}\exists\quad\underset{v\in
  UX}\forall\quad\underset{r>R}\forall\quad
\Hd(\horo_v^\Delta,\sphere_v^\Delta(r))<\delta$$
\end{lemm}

We will need another technical lemma on spheres, which quantifies the fact
that on a big sphere $\sphere_v(a)$ vectors with close base  points are
close with respect to the metric in $UX$. The statement of the lemma is
illustrated in Figure~\ref{nahevektoren}.

\begin{lemm}{}{close radial vectors}{}             %Lemm close radial vectors
  Suppose $X$ is a Hadamard manifold with compact quotient $M=X/\Gamma$ and
  the \Gam-compact subset $\tilde K\subset X$ 
  and
  $a,\delta>0$ are given. Then there is a $\delta'$ such that for every
  $q\in X$ it holds:
  \bigskip\\
  Call $\sigma$ the geodesic ray starting in $\sigma(0)=q$ with
  $\sigma(s)=p\in \tilde K$ and let $w\in U_qX$ be any vector such that the
  minimal distance between $\gamma_w$ and $p$ is less than $\delta'$. If
  $d(p,q)\geq a$, then for any point $p':=\gamma_w(t)$ with
  $d(p,p')<\delta'$
  $$d^1(\dot\gamma_w(t),\dot\sigma(s))<\delta.$$

\begin{figure}[h]\caption{Close Radial Vectors}\label{nahevektoren}
\begin{center}
  \psfrag{p}{$p\in\tilde K$} \psfrag{p'}{$p'$} \psfrag{q}{$q$} \psfrag{w}{$w$}
  \psfrag{<delta'}{$\leq\delta'$} \psfrag{gammaw}{$\gamma_w$} \psfrag{sigma}{$\sigma$}
  \psfrag{dotsigma}{$\dot\sigma(s)$} \psfrag{dotgamma}{$\dot\gamma_w(t)$}
\setlength{\meinbuffer}{\textwidth * 3/4}  
\includegraphics[width=\meinbuffer]{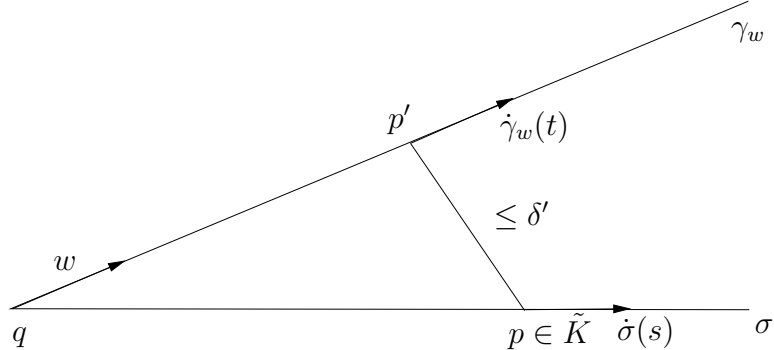}
\end{center}\end{figure}

\end{lemm}

\begin{lemmp}{close radial vectors}{}{}
  For a start suppose $\delta'<a$. We will fix $\delta'$ at the end of the
  first step.\\
  Step 1:\\
On the set
$$Def:=\underset{q\in X}\bigcup \left(\{q\}\times\underset{p\in\tilde
      K\backslash\{q\}}\bigcup\left(\{p\}\times\left(\overline\U_{\delta'}(p)\backslash\{q\}\right)\right)\right)\subset X^3$$
we define a continuous function $f$ as follows:\\
For $q\in X$, $p\in\tilde K\backslash\{q\}$ and
$p'\in\U_{\delta'}(p)\backslash\{q\}$ call $\sigma$ the geodesic connecting
$q=\sigma(0)$ to $p=\sigma(s)$ and $\gamma$ the geodesic connecting
$q=\gamma(0)$ to $p'=\gamma(t)$. Then
$$f(q,p,p'):=d^1(\dot\sigma(s),\dot\gamma(t))$$
is a \Gam-compatible map, since 
$$f(\rho(q),\rho(p),\rho(p'))=f(q,p,p')\quad\text{ for all }\quad\rho\in \Gamma.$$ 
For every $p\in\tilde K$ consider the (compact) sphere
$${\cal K}_p:=\{q\in X\;|\;d(p,q)=a\}$$
and the continuous function $f$ on the \Gam-compact set
$$Def_2:=\underset{p\in\tilde
  K}\bigcup\{p\}\times\overline{\U_{\delta'}(p)}\subset X^2$$
defined by 
$$f(p,p'):=\underset{q\in{\cal K}_p}\max \,f(q,p,p').$$
This function is continuous on $Def_2$ and satisfies $f(p,p)=0$ for all
$p\in \tilde K$. Hence we can suppose that $\delta'$ was chosen so small
that $f(p,p')<\delta/2$ for all $p,p'$ with $d(p,p')<\delta'$.\\
This means that $d^1(\dot\sigma(s),\dot\gamma(t))<\delta/2$ for all $q\in X$, $p\in \tilde K$ and $p'\in X$ with $d(p,q)=a$
and $d(p,p')<\delta'$.\\
We may assume that $\delta'$ satisfies furthermore $\arcsin(\delta'/a)<\delta/2$.
  Step 2:\\
  Now let $q\in X$ be a point with $d(p,q)>a$ for a given $p\in \tilde K$ and let
  $w\in U_qX$ be a vector
  with $d(\gamma_w,p)<\delta'$. Suppose there is a  $p'=\gamma_w(t)\in{\cal U}_{\delta'}(p)$.\\
  Let $\sigma$ be the geodesic ray through $p$ originating in $q$.  There
  is a unique point $Q$ on $\sigma$ between $q$ and $p$ with $d(Q,p)=a$. Now let
  $\gamma$ be the geodesic ray originating in $Q$ with $p'=\gamma(t')$.
  Write $v\in U_{p'}X$ for the vector you get by moving the vector
  $\dot\sigma(s)$ from $p$ to $p'$
  via parallel transport along the connecting geodesic. This is illustrated
  in Figure~\ref{zweitesbild}.

\begin{figure}[H]\caption{}
\label{zweitesbild}
  \psfrag{p}{$p$} \psfrag{p'}{$p'$} \psfrag{q}{$q$} \psfrag{Q}{$Q$}
  \psfrag{v}{$v$} \psfrag{=a}{$=a$}
  \psfrag{<delta'}{$\leq\delta'$} \psfrag{gamma}{$\gamma$}
  \psfrag{sigma}{$\sigma$}\psfrag{gammaw}{$\gamma_w$}\psfrag{w}{$w$}
\psfrag{gamma(t')}{$\dot\gamma(t')$} \psfrag{sigma(s)}{$\dot\sigma(s)$}\psfrag{gammaw(t)}{$\dot\gamma_w(t)$}\psfrag{winkel}{$\measuredangle_Q(p,p')$}
  \setlength{\meinbuffer}{\textwidth * 4/5}
  \includegraphics[width=\meinbuffer]{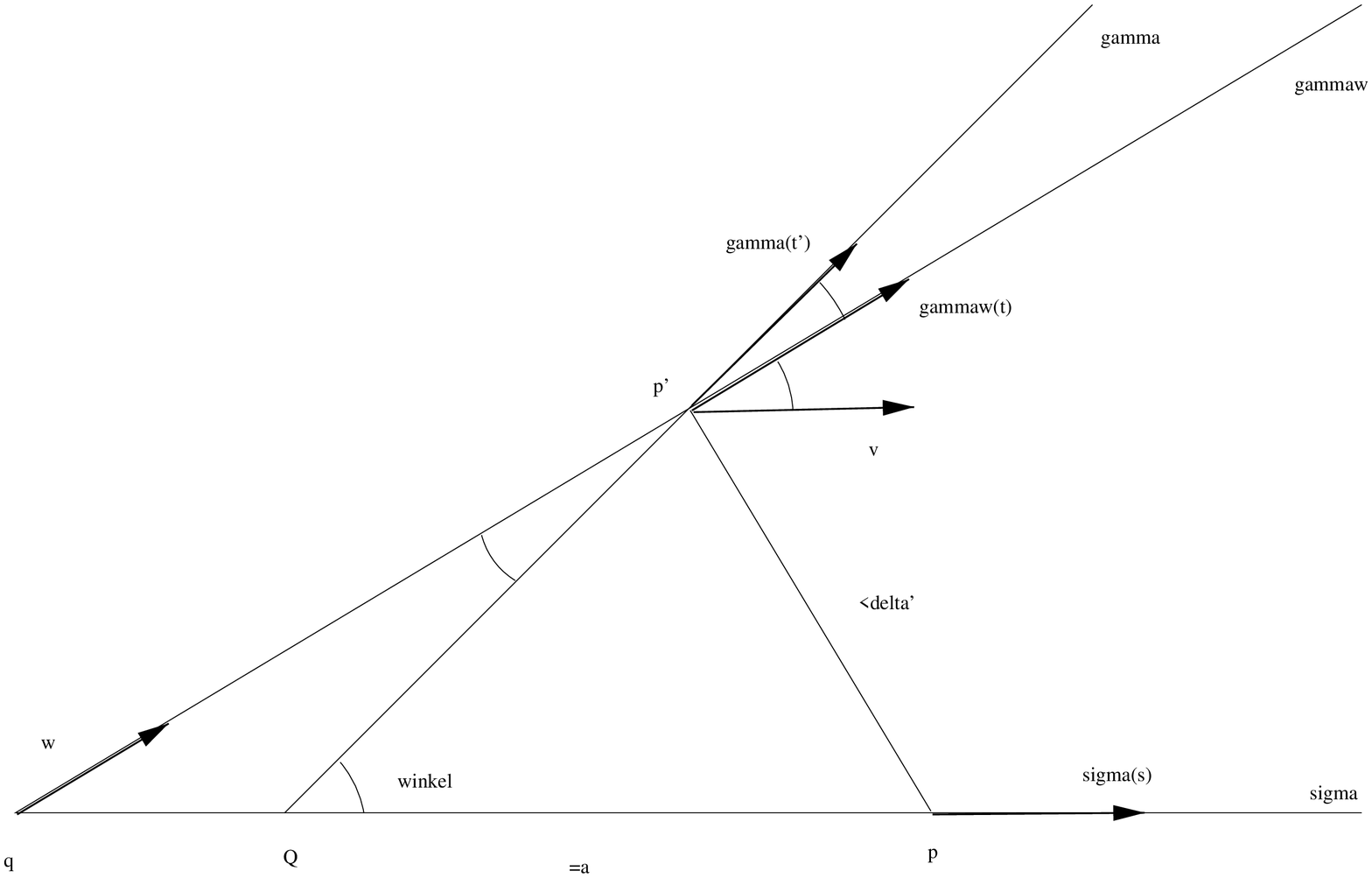}
\end{figure}
By the first step we know that
  $d(p',p)+\measuredangle(\dot\gamma(t'),v)=d^1(\dot\gamma(t'),\dot\sigma
  (s))<\delta/2$.
Furthermore it is easy to see that
  $\measuredangle(\dot\gamma_w(t),\dot\gamma(t'))\leq\measuredangle_Q(p,p')$ since, by nonpositive curvature, the sum of angles in the triangle $\Delta(q,Q,p')$ must be less or equal to $\pi$.\\
Therefore
\begin{align*}
  d^1(\dot\gamma_w(t),\dot\sigma(s))&=d(p',p)+\measuredangle(\dot\gamma_w
  (t),v)\\
  &\leq
  d(p',p)+\measuredangle(\dot\gamma_w(t),\dot\gamma(t'))+\measuredangle( \dot\gamma(t'),v)\\
  &<\measuredangle_Q(p,p')+\delta/2\\
  &\leq
  \arcsin(\delta'/a)+\delta/2\leq\delta/2+\delta/2=\delta.
\end{align*}

Where the last line uses the fact that the manifold is
nonpositively curved.
\end{lemmp}

Since $X$ itself is \Gam-compact if it has a compact quotient $M=X/\Gamma$,
it is an obvious conclusion that vectors on big spheres are close if their
base points are close. This result is formulated in the next corollary.

\begin{coro}{}{fusspunktnahe}{}
Suppose $X$ is a Hadamard manifold with compact quotient and $r_0>0$ any
given radius. Then for any $\delta>0$ there is a $\delta'>0$ such that
vectors $v,w$ on the same sphere, say $\sphere_v(r)$, of radius $r\geq r_0$
are $\delta$-close if their base points are $\delta'$-close. I.\ e.\ for
all $v\in UX$, $r\geq r_0$ 
$$\underset{w\in\sphere_v(r)}\forall\qquad\quad
d(p(v),p(w))<\delta'\qquad\Longrightarrow\qquad d(v,w)<\delta.$$
\end{coro}

\mysubsection{Flats on a Horosphere}{Flats on Horospheres}\label{flats}
In this section let $X$ denote a \emph{real analytic} Hadamard manifold.
Fix a vector $v\in UX$ and consider the set of parallel vectors $\bars
v$\index{v@$\bars v$}\index{$v$@$\bars v$}. This is a 
\index{flat(v)@$\flat(v)$}\index{$flat(v)$@$\flat(v)$}
$\flat(v)$-dimensional submanifold of $UX$ and $UX$ is
partitioned into the $\bars v$.
For any vector $u\in\bars v$ the Busemann function\index{Busemann function}\index{function!Busemann} $b_v$ along the geodesic
$\gamma_u$ grows proportionally to arc length. Hence there is a unique
vector where $\dot\gamma_u$ hits the horosphere $\horo_v$ transversally.
Obviously two geodesics intersect this horosphere in different points and
therefore the
intersection of $\bars v$ with $\horo_v$ is a ($\flat(v)-1$)-dimensional
submanifold of $UX$. If $v\in UX$ is any vector we will write $\frak H_v:=\bars
v\cap\horo_v$
\index{H v@$\frak H_v$}\index{$h v$@$\frak H_v$}
for this submanifold of all vectors parallel to $v$ on the
same horosphere.
$\frak H_v$ is a subset of the leaf $W^u(v)$ of the unstable bundle
which is homeomorphic to the manifold $X$ itself via the gradient of the
Busemann function $b_v$: 
$$
\index{grad b v@$\grad b_v$}\index{$grad b v$@$\grad b_v$}
\grad b_v:X\overset{\sim}\rightarrow W^u(v).$$

The main result of this section is Corollary~\ref{coro2}: For a rank one
manifold the vectors of rank one are dense on any sphere and horosphere.
\begin{lemm}{}{rangfolgtflach}{}
  On an analytic Hadamard manifold $X$ suppose there is an open segment
  $\horo_v^\Delta$ of a horosphere $\horo_v$ which consists entirely of
  higher rank vectors. Then all vectors in $\horo_v^\Delta$ have
  flatness at least two.\\
\end{lemm}
\begin{lemmp}{rangfolgtflach}{}{}
  The projection $p:UX\rightarrow X$ identifies $\horo_v$ with the
  horosphere in $X$. Write $H_v^\Delta:=p(\horo_v^\Delta)$ for the
  horosphere segment in $X$ consisting of all base points of vectors in
  $\horo_v^\Delta$. Denote the gradient field of the Busemann function
  $b_v$ by $V:H_v\rightarrow
  \horo_v$.\\
  
  Since all elements of $\horo_v^\Delta$ are of higher rank, given any
  point in $H_v^\Delta$ we can find locally a continuous unit vector field $Z$
  such that $Z_q$ is orthogonal to $V_q$ and can be extended to a parallel
  Jacobi vector field
  along $\gamma_{V_q}$ whenever $Z_q$ is defined.\\
  Fix any $u\in\horo_v^\Delta$ and $Z$ defined in a neighbourhood of
  $p(u)$.  Take an integral curve $\sigma:[0,\delta]\rightarrow
  \horo^\Delta$ of $Z$ with $\sigma(0)=p(u)$ and consider the geodesic
  variation
\begin{align*}
  \Gamma:[0,\delta]\times\mathbb R&\longrightarrow X\\
  (s,t)&\longmapsto\gamma_{V_{\sigma(s)}}(t).
\end{align*}

All the geodesics $\Gamma_s$ originate from the same point in $X(\infty)$.
Thus the corresponding variational vector field
$J_s:t\mapsto\frac\partial{\partial s}\Gamma(s,t)$ along $\Gamma_s=\gamma_{V_{\sigma(s)}}$ is an
unstable Jacobi field and for $t=0$ equals $\dot\sigma(s)=Z_{\sigma(s)}$.
But $J_s$ is parallel along $\Gamma_s$ hence unstable and both
fields coincide in $t=0$:
$$J_s(0)=\frac\partial{\partial s}\Gamma_s(0)=\frac\partial{\partial s}\gamma_{V_{\sigma(s)}}(0)=\frac\partial{\partial s}p(V_{\sigma(s)})=\frac\partial{\partial s}\sigma(s)=Z_{\sigma(s)}$$

Since unstable vector fields along a geodesic are characterised by one
value, the variational vector field must coincide with the parallel Jacobi
field

$$\frac{\partial}{\partial s}\Gamma(s,t)=J_s(t)$$
which is of constant length $\|J_s(t)\|\equiv\|Z_{\sigma(s)}\|=1$.
We can therefore calculate the length of the curves $\Gamma_t$ which
connect $\gamma_u(t)$ to $\gamma_{V_{\sigma(\delta)}}(t)$ to

$$L(\Gamma_t)=\int\limits_0^\delta\|\frac{\partial}{\partial
  s}\Gamma(s,t)\|ds=\int\limits_0^\delta\|J_s(t)\|ds=\delta$$

Now we can approximate the distance of the two geodesics $\gamma_u$ and
$\gamma_{V_{\sigma(\delta)}}$ by $\delta$:
$$d(\gamma_u(t),\gamma_{V_{\sigma(\delta)}}(t))\leq L(\Gamma_t)=\delta$$
or all $t\in \mathbb R$.

Thus these two geodesics span a flat strip and $\flat(u)$ is at least two by analyticity.
\end{lemmp}

\begin{prop}{}{prop1}{}
  Let $X$ be an analytic Hadamard manifold with compact quotient. Either
  for each horosphere the subset of rank one vectors is dense or there is a
  horosphere that consists entirely of vectors of higher rank.\\
\end{prop}
\begin{propp}{prop1}{}{}
  Suppose there is a horosphere on which the subset of rank one vectors is
  not dense, and suppose that every horosphere contains at least one vector
  of rank one.\\
  Define the function $l:UX\rightarrow\mathbbR$ as follows: For $v\in UX$
  write $l(v):=\inf\{d_v(v,w)\;|\;w\in\horo_v,\rank(w)=1\}$ for the minimal base point distance of $v$ to a rank one vector in
$\horo_v$. We know $l\not\equiv 0$, since at least one horosphere
contains
an open subset of higher rank vectors.\\
The foliation of $UX$ into horospheres is continuous and 
the set of rank one vectors is open. Hence the function $l$ is
semicontinuous:
$$\underset{v\in
  UX}\forall\quad\underset{\epsilon>0}\exists\quad\underset{u\in {\cal
    W}_\epsilon(v)}\forall\quad l(u)\leq l(v)
$$
Because the manifold has a compact quotient, $l$ attains its maximal
value. But we will show that $l$ is unbounded. This is a contradiction.\\
Call $v_1\in UX$ the vector in $UX$ where $l$ attains its maximal value
$\frak l:=l(v_1)$. By definition of $l$ all $u\in\horo_{v_1}^{\frak l}$ are
of higher rank.  Choose a vector $v_0$ of minimal flatness $\frak
f:=\flat(v_0)$ in this segment of $\horo_{v_1}$ and a neighbourhood
$\horo^\epsilon_{v_0}\subsetneq\horo^{\frak l}_{v_1}$ of $v_0$ of constant flatness. This segment
is foliated into submanifolds $\{\frak
H_v\cap\horo^\epsilon_{v_0}\}_{v\in\horo^\epsilon_{v_0}}$ of dimension
$\frak f-1$. Within any of these submanifolds are only parallel vectors of
constant flatness. Notice that these submanifolds
are at least one-dimensional by Lemma~\ref{rangfolgtflach}.\\

Fix $\frak B_\epsilon\subset\horo^\epsilon_{v_0}$ an ($n-\frak
f$)-dimensional simply connected submanifold transversal to the foliation
$\{\frak H_v\cap \horo^\epsilon_{v_0}\}_{v\in\horo^\epsilon_{v_0}}$
  intersecting $\frak H_{v_0}$ in $v_0$ and satisfying $\partial\frak B_\epsilon\subset\partial\horo^\epsilon_{v_0}$.

Define for all $v\in\frak B_\epsilon$ the disk of parallel vectors of distance less than $3\frak l$:

$$\frak D_v:=\frak H_{v}\cap\horo^{3\frak l}_{v}\subset\horo_{v_0}$$
All vectors in
$\frak D_v$ are of
flatness at least $\frak f$, like $v$, since they lie in $\bars v$.\\
The set
$$\frak A:=\underset{v\in\frak B_\epsilon}{\bigcup}\frak D_v$$
which
consists of higher rank vectors only, is a submanifold of $W^u(v_0)$, the
unstable leaf of $v_0$ and is diffeomorphic to $\frak B_\epsilon\times\frak
D_{v_0}$, the product of two disks of dimensions $(N-\frak f)$ and $(\frak
f-1)$ respectively. The boundary $\partial\frak A$ consists of the two sets
$$\partial_1\frak A:=\underset{v\in\partial\frak B_\epsilon}\bigcup \frak
D_v$$
and
$$\partial_2\frak A:=\underset{v\in\frak B_\epsilon}\bigcup \partial\frak
D_v.$$
For all vectors $u\in\partial\frak D_v\subset\partial_2\frak A$ and all $t>0$
$$d_{\phi_t(v_0)}(\phi_t(u),\phi_t(v_0))\geq d_{v_0}(u,v_0)\geq
d_{v_0}(u,v)-d_{v_0}(v,v_0)\geq 3\frak l-\epsilon>2\frak l.$$
The
vectors $u\in\partial_1\frak A$ are not parallel to $v_0$, but negatively
asymptotic to $v_0$. Therefore for
all $t>0$
$$d_v(u,v_0)<d_{\phi_t(v)}(\phi_t(u),\phi_t(v_0))$$
and we can define
$$\delta:=\underset{u\in\partial_1\frak A}\min \left[d_{\phi_1(v_0)}(\phi_1(u),
\phi_1(v_0))-d_{v_0}(u,v_0)\right]>0$$
since $\partial_1\frak A$ is compact; without loss of generality assume
$\delta<2\frak l$. Now
by convexity of the distance function we see that
$$d_{\phi_{2\frak l/\delta}(v_0)}(\phi_{2\frak l/\delta}(u),\phi_{2\frak
  l/\delta}(v_0))>\frac{2\frak l}\delta\delta>2\frak l.$$
$\phi_{2\frak
  l/\delta}$ is a diffeomorphism of $\horo_{v_0}$ onto $\horo_{\phi_{2\frak
    l/\delta}(v_0)}$ respecting rank and flatness of a vector. Hence $\frak
A$ is mapped onto a subset $\phi_{2\frak l/\delta}\frak A$ that is foliated
by $(\frak f-1)$-dimensional submanifolds consisting of parallel vectors and therefore contains only vectors of
flatness at least $\frak f$. The boundary of $\phi_{2\frak l/\delta}\frak
A$ is $\phi_{2\frak l/\delta}(\partial_1\frak A)\cup\phi_{ 2\frak
  l/\delta}(\partial_2\frak A)$. For all vectors in this boundary the
distance to $\phi_{2\frak l/\delta}(v_0)$ is greater than $2\frak l$. Thus
this set contains the set $\horo^{2\frak l}_{\phi_{2\frak l/\delta}(v_0)}$
and hence $l(\phi_{2\frak l/\delta}(v_0))\geq 2\frak l$ which is a
contradiction to the maximality of $\frak l$.
\end{propp}

Combining the results of Proposition~\ref{prop1} and Lemma~\ref{rangfolgtflach} we
get

\begin{coro}{}{coro1}{}
  Let $X$ be an analytic Hadamard manifold with compact quotient. Either
  for all horospheres the subset of rank one vectors is dense or there is a
  horosphere that consists entirely of vectors of flatness at least two.\\
\end{coro}

\begin{coro}{}{coro2}{}
  Let $X$ be an analytic rank one Hadamard manifold with compact quotient.
  Then for any horosphere and any sphere in $UX$ the subset of rank one
  vectors is dense.
\end{coro}

\begin{corop}{coro2}{}{}
  We start by considering the statement of the corollary for horospheres.\\
  Ballmann showed in \cite{Ballmann1982} that for a rank one manifold there
  exists a point $\xi\in X(\infty)$ which can be joined to every other
  point in
$X(\infty)$ by a geodesic in $X$.%
\footnote{Ballmann showed that the set of points with this property is
  dense in $X(\infty)$, but we will only need the existence of one point}\\
Suppose there is a horosphere $\horo_v$ consisting entirely of vectors of
flatness at least two. Write $\zeta:=\gamma_v(-\infty)\in\bar X$ for the
center of $\horo_v$. By Ballmann's result there is a geodesic $\sigma$
joining $\zeta$ to $\xi$ if $\xi\neq\zeta$, if $\xi=\zeta$ fix any geodesic
starting in $\xi$.  In both cases this geodesic meets $\horo_v$ and hence is of flatness
at least two. Thus it lies inside a flat. This is a contradiction to the
property of $\xi$: the points on the boundary of the flat can not be joined
to $\xi$ by a geodesic in $X$, since their Tits distance 
\index{Tits distance}\index{distance!Tits} is finite.\\
Now consider the statement for spheres.\\
Suppose there is a segment $\sphere_v^\Delta$ of a sphere that consists
only of higher rank vectors. The geodesic flow increases the radius of the
sphere and the
diameter of the
segment. Hence we may assume that the radius of the sphere is arbitrarily
large. By compactness there is a radius $\delta$ such that every horosphere
segment of diameter $\Delta$ contains a rank one vector which is surrounded by an open ball of
radius $\delta$. Now for spheres of radius bigger than $T$ the segment
$\sphere_v^\Delta$ is $\delta$-close to a horosphere segment
(cf. Lemma~\ref{bigspheres}) and therefore must contain a rank one
vector, which is a contradiction.
\end{corop}

\mysection{Hyperbolicity of Rank One Vectors}{Hyperbolicity}\label{hyperbolic.tex}
A parallel Jacobi field along a geodesic in a Hadamard manifold might be
considered as an infinitesimal flat. In the neighbourhood of a rank one vector
$v$ we will therefore expect to find not linear but some kind of hyperbolic
behaviour. Sergei Buyalo\index{Buyalo, Sergei} and Viktor
Schroeder\index{Schroeder, Viktor} showed in~\cite{whoknows} that this is in
fact true: The distance of our geodesic to other close geodesics
has some doubling property.
The aim of this section is to quantify this. The result is summed up in
Corollary~\ref{best hyper}. To get there we will need several rather
technical lemmas.\\ 
We introduce the notion of a \emph{hyperbolic vector}, which is closely
related to Ballmann's notion of a \emph{hyperbolic
  geodesic}\index{hyperbolic geodesic}\index{geodesic!hyperbolic} in%
~\cite{Ballmann1995}.\\
Remember that for a compact manifold $M$ of nonpositive curvature
$\pi:X\rightarrow M$ denotes the universal covering. I.\ e.\ $X$ is a
Hadamard manifold with compact quotient $M=X/\Gamma$.

\begin{defi}{Tracing Neighbourhoods}{}{}
  Let $X$ be a Hadamard manifold.\\
  Consider $K\subset UX$ and $L,\epsilon>0$ and call $\theta$ a
  \emph{tracing distance}\index{distance!tracing}\index{tracing distance} for $K$ with respect to
  the constants $L,\epsilon$ if
  and only if the following is true:\\
  For any $v\in K$ and $w\in \W_{\theta}(v)$ the end points of the
  geodesic segments $\dot\gamma_v([-L,L])$ and $\dot\gamma_w([-L,L])$ are
  within $d^1$-distance less than $\epsilon$ of each other, i.\ e.\  
$$
d(v,w)<\theta
\;\Longrightarrow \;
\max\left\{
   d^1\left(\dot\gamma_v(-L),\dot\gamma_w(-L)\right)
   ,
   d^1\left(\dot\gamma_v(L),\dot\gamma_w(L)\right)
\right\}
<\epsilon.
$$
Under this condition we call $\W_{\theta}(v)$ the \emph{tracing
    neighbourhood}\index{neighbourhood!tracing}\index{tracing neighbourhood} of $v$ with respect to the
  constants $L,\epsilon$ or simply the \emph{$(L,\epsilon)$-tracing
    neighbourhood} of $v$.
\end{defi}

\begin{rema}{}{}{}
  Since $X$ is a Hadamard manifold the function $t\mapsto
  d(\gamma_v(t),\gamma_w(t))$ is convex.  It follows immediately that for
  any $v\in K$ and $w\in \W_{\theta}(v)$ the geodesic segment
  $\dot\gamma_w([-L,L])$ is contained in the $\epsilon$-neighbourhood of
  the geodesic $\dot\gamma_v$ provided $\theta$ is the tracing distance
  for $K$ with
  respect to the constants $L,\epsilon$.\\
  Notice that, by $\,p^* d\leq d^1$, this implies that the geodesic
  segment $\gamma_w([-L,L])$ is contained in the $\epsilon$-neighbourhood
  of the geodesic $\gamma_v$.
\end{rema}

\begin{lemm}{}{tracing}{}                            %Lem tracing
  \begin{enumerate}
\item
  Let $X$ denote a Hadamard manifold and $K\subset UX$ a compact
  subset. Then for any choice of the constants $T$ and $\epsilon$ there is a
  tracing distance $0<\theta<\epsilon$ for $K$ with respect to the
  constants $T, \epsilon$.
\item\label{tracing Gamma}
  Let $X$ denote a Hadamard manifold with compact quotient
  $M=X/\Gamma$. Then for 
  any $\Gamma$-compact subset $K\subset UX$ and 
  any choice of the constants
  $T$ and $\epsilon$ there is a tracing distance $0<\theta<\epsilon$ for $K$ with
  respect to the constants $T, \epsilon$.
\item
  Let $X$ denote a Hadamard manifold with compact quotient. Then $UX$ is
  $\Gamma$-compact and hence by~\ref{tracing Gamma} for
  any choice of the constants
  $T$ and $\epsilon$ there is a global tracing distance $0<\theta<\epsilon$ with
  respect to the constants $T, \epsilon$.
\end{enumerate}
\end{lemm}

\begin{lemmp}{tracing}{}{}
  Suppose there are constants $T,\epsilon$ such that there is no tracing
  distance.  Then for every $\epsilon_i:=\epsilon/i$ there are vectors
  $v_i\in K$ and $w_i\in \W_{\epsilon_i}(v_i)$ such that
\begin{equation}\label{sternli}%\tag{*}
d^1(\dot\gamma_{w_i}(-T),\dot\gamma_{v_i}(-T))>\epsilon\;\text{ or
 }\;d^1(\dot\gamma_{w_i}(T),\dot\gamma_{v_i}(T))>\epsilon.
\end{equation}
In case $K$ is compact there are subsequences $(v_i)$ and $(w_i)$
converging in $K$. Obviously they converge to the same vector, say $v\in
K$.\\
In case $K$ is $\Gamma$-compact we do not necessarily find
subsequences. But surely we find two converging sequences (we will call
them  $(v_i)$ and $(w_i)$, again) satisfying~(\ref{sternli}). Again, these
sequences converge to the same vector, say $v\in K$.\\
In both cases we can, by taking subsequences, assume that the
  inequality for $T$ holds for all $i$. (If we do not find a subsequence
  like this we can use an analogous argument for $-T$.)  Thus
  $$0<\epsilon<d^1(\dot\gamma_{w_i}(T),\dot\gamma_{v_i}(T))\rightarrow
  d^1(\dot\gamma_v(T),\dot\gamma_v(T))=0$$
  which is a contradiction.
\end{lemmp}

Before we can introduce the notion of hyperbolicity, it is necessary to
recall the definition of two different distances between subsets in a
metric space. The \emph{(minimal) distance}
\index{distance!minimal}\index{minimal distance} 
is defined by
$$d(A,B):=\underset{a\in A, b\in B}{\inf}d(a,b)\index{d(A,B)@$d(A,B)$}$$
while the
\emph{Hausdorff distance}\label{defHd}\index{distance!Hausdorff}\index{Hausdorff distance} is defined as the
smallest number $\epsilon$ such that each set is completely contained in
the $\epsilon$-neighbourhood of the other.
$$\Hd(A,B):=\max(\underset{a\in A}{\sup}\,d(a,B),\underset{b\in
  B}{\sup}\,d(A,b))\index{Hd@$\Hd(A,B)$}\index{$hd$@$\Hd(A,B)$}$$

Notice that the Hausdorff distance is a metric taking values in
$\mathbbR_+\cup\{\infty\}$. The minimal distance is finite but not a metric.

\begin{defi}{Hyperbolicity}{}{}  %Def Hyperbolicity
  A unit vector $v\in UX$ on a Hadamard manifold is said to satisfy the \emph{hyperbolicity
    condition}\index{condition!hyperbolicity}\index{hyperbolicity condition} for
  $\mu\in]0,1[$, $T>0$ and $\delta>0$ if the following holds:\\
  For any two vectors $w_1,w_2$ in the $(T,\delta)$-tracing neighbourhood
  of $v$ the minimal distance
  $d(\gamma_{w_1}([-T,T]),\gamma_{w_2}([-T,T]))$ is less than $\mu$ times
  the Hausdorff distance
  $\Hd\left(\gamma_{w_1}([-T,T]),\gamma_{w_2}([-T,T])\right)$.\bigskip\\
  The vector $v\in UX$ is called \emph{hyperbolic}\index{vector!hyperbolic}\index{hyperbolic vector}, if there is a $\delta$ such that for any $\mu\in ]0,1[$ we can
  choose $T:=T(\mu,\delta)$ to make $v$ satisfy the hyperbolicity condition
  for $\mu$, $T$ and $\delta$.
\end{defi}

\begin{prop}{}{rank 1 hyper}{}%Prop rank 1 hyper
  On a Hadamard manifold $X$ %with compact quotient 
  every rank one vector is hyperbolic.\smallskip\\
  Furthermore:\\
  Suppose $K\subset UX$ consists entirely of vectors of rank one and either
  $K$ is
  compact or $X$ has compact quotient $M=X/\Gamma$ and $K$ is
  $\Gamma$-compact.\\
  Choose $\delta$ such that the $3\delta$-neighbourhood
  of $K$ consists of 
  rank one vectors, too.  Then for any $\mu\in]0,1[$ we can find 
  $T:=T(\mu,\delta,K)>0$ such that every vector $v\in K$ satisfies the
  hyperbolicity condition for $\mu,T,\delta$.
\end{prop}

%%%%%%%%%%%%%%%%%%%%%%%%%%%%%%%%%%%%%%%%%%%%%%%%%%%%%%%%%%%%%%%%%%%%%%%%%%%%

\begin{propp}{rank 1 hyper}{}{}
  Suppose this is false. Then there is a $\mu\in]0,1[$ for which we will
  find a sequence of vectors $z_n\in K$ and
  $v_n,w_n\in\W_{\theta_n}(z_n)$, where $\theta_n$ is the tracing
  distance of $z_n$ with respect to $n,\delta$. Thus we have
\begin{xalignat*}{2}
  d(\gamma_{v_n}(-n),\gamma_{z_n}(-n))&<\delta&d(\gamma_{w_n}(-n),\gamma_{z_n}(-n))&<\delta\\
  d(\gamma_{v_n}(n),\gamma_{z_n}(n))&<\delta&d(\gamma_{w_n}(n),\gamma_{z_n}(n))&<\delta
\end{xalignat*}
$$d(\gamma_{w_n}[-n,n],\gamma_{v_n}[-n,n])>\mu\,
\Hd(\gamma_{w_n}[-n,n],\gamma_{v_n}[-n,n])=:\mu\, \Hd_n.$$

So for any two points on the two geodesic segments the distance is bigger
than $\mu$ times the Hausdorff distance $\Hd_n$ which is the maximal distance
between corresponding end points
$\gamma_{v_n}(\pm n)$ and $\gamma_{w_n}(\pm n)$ and therefore $\Hd_n<2\delta$.\\

First, by ($\Gamma$-)compactness, we can assume that $z_n$ converges to some
$z\in K$. 
Then, for $n$ large enough, we see that
$v_n,w_n\in\W_{\theta_n}(z_n)\subset \W_\delta(z_n)\subset
\W_{2\delta}(z)$ and we can suppose that $v_n\rightarrow \bar v$ and
$w_n\rightarrow \bar w$ with $\bar v,
\bar w\in \overline{\W_{2\delta}(z)}$.\\
Now fix $k\in \mathbb N$ and take an arbitrary $n\gg k$. Take a look at the
convex function
\begin{align*}
  f_n:[-n,n]&\longrightarrow[\mu,1]\\
  r&\longmapsto \frac{1}{\Hd_n}d(\gamma_{w_n}(r),\gamma_{v_n}[-n,n])
\end{align*}
\begin{figure}[h]\caption{\label{abbhyper}}
  \psfrag{mn}{$-n$} \psfrag{mk}{$-k$} \psfrag{fn}{$f_n$} \psfrag{n}{$n$}
  \psfrag{k}{$k$} \psfrag{r}{$r$} \psfrag{re}{$r+\epsilon$}
  \psfrag{fnmn}{$f_n(-n)\leq 1$} \psfrag{fnn}{$f_n(n)\leq 1$}
  \psfrag{fnr}{$f_n(r)\geq \mu$} \psfrag{fnre}{$f_n(r+\epsilon)\geq \mu$}
  \setlength{\meinbuffer}{\baselineskip} \multiply \meinbuffer by 12
  \includegraphics[width=\textwidth,height=\meinbuffer]{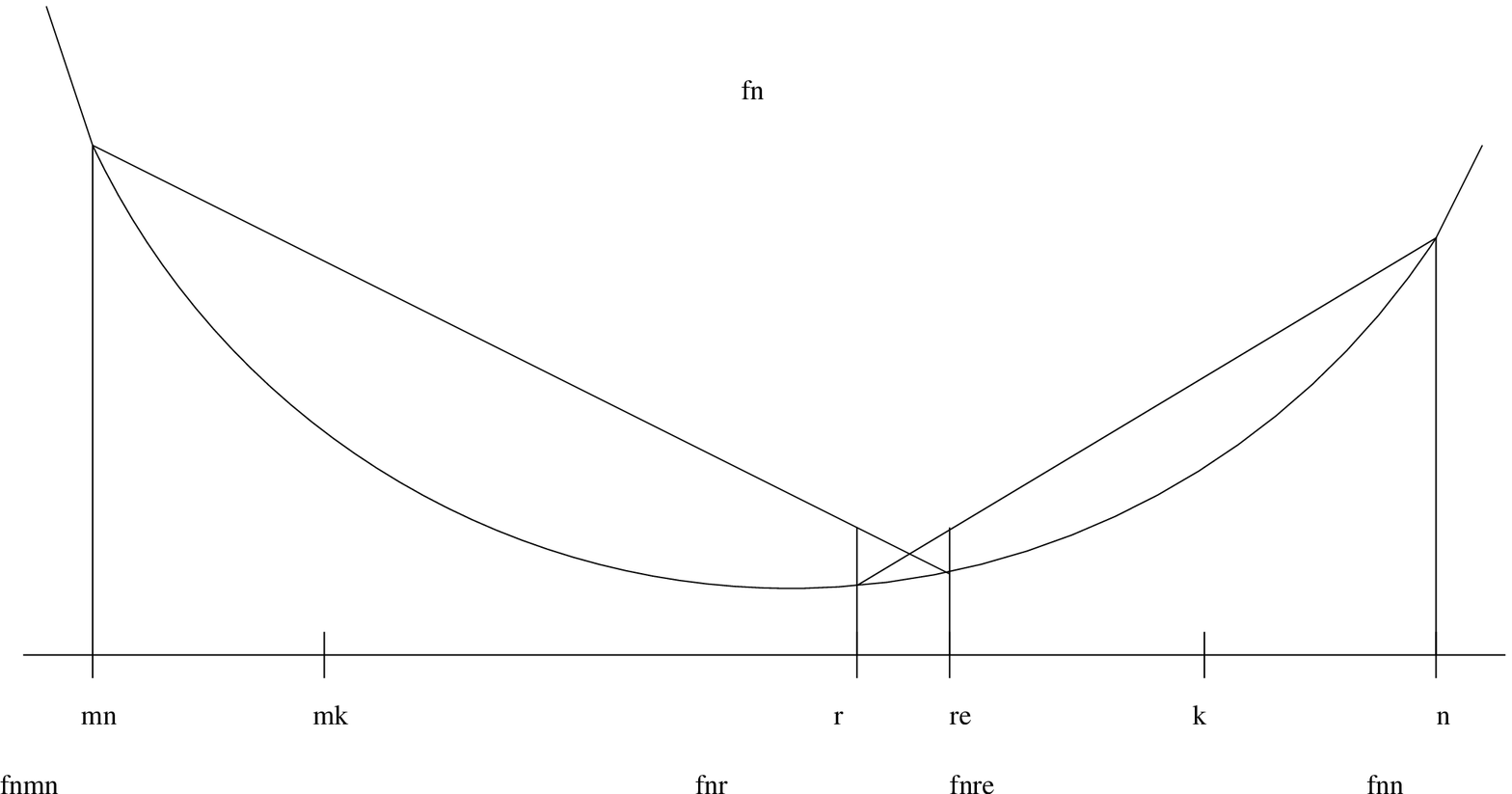}
\end{figure}

For $r\in[-k,k]$, $\epsilon$ small, by convexity we have (as illustrated in Figure~\ref{abbhyper}):
\begin{align}
  f_n(r+\epsilon)&\leq f_n(r)+\frac{\epsilon}{n-r}(f_n(n)-f_n(r))\label{oben}\\
  f_n(r)&\leq
  f_n(r+\epsilon)+\frac{\epsilon}{n+r}(f_n(-n)-f_n(r+\epsilon)).\label{unten}
\end{align}
Rewriting (\ref{oben}) we get
$$f_n'(r)\overset{0\leftarrow\epsilon}\longleftarrow\frac{f_n(r+\epsilon)-f_n(r)}{\epsilon}\leq
\frac{f_n(n)-f_n(r)}{n-r}\leq\frac{1-\mu}{n-k}$$
and a similar computation
using (\ref{unten}) yields a lower bound for the derivative of $f_n
\Big|_{[-k,k]}$. Thus
$|f_n'(r)|<\frac{1-\mu}{n-k}\underset{n\rightarrow\infty}\longrightarrow
0$ for all $r\in[-k,k]$.\\

\begin{figure}[h]\caption{Definition of $\alpha_r^n$}{\label{alpharn}}
  \psfrag{gvn}{$\gamma_{v_n}$} \psfrag{gwn}{$\gamma_{w_n}$}
  \psfrag{alphan}{$\alpha_r^n$} \psfrag{gwnr}{$\gamma_{w_n}(r)$}
  \setlength{\meinbuffer}{\baselineskip} \multiply \meinbuffer by 12
\begin{center}
  \includegraphics[width=\textwidth/2]
{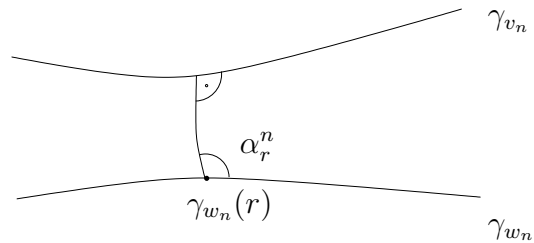}
\end{center}
\end{figure}

Now we can consider
$L_n(r):=\Hd_nf_n(r)=d(\gamma_{w_n}(r),\gamma_{v_n}([-n,n]))$.  This is the
length function of a geodesic variation. Applying the first variation
formula we find that
\begin{align}
\left|\cos(\alpha_r^n)\right|
&=\left|\cos(\alpha_r^n)+cos(\pi/2)\right| \notag\\
&=\left|\dot L_n(r)\right| = \Hd_n\left|f_n'(r)\right| 
\leq2\delta\,\frac{1-\mu}{n-k}\rightarrow 0  
\label{dotL}
\end{align}
where $\alpha_r^n$ is the angle between
$\gamma_{w_n}$ and the geodesic linking $\gamma_{w_n}(r)$ to
$\gamma_{v_n}$, as illustrated in Figure~\ref{alpharn}.

Thus for $|r|<k$ this angle goes to $\pi/2$ for $n\rightarrow \infty$.\\

Now consider the subsequence such that $v_n\rightarrow \bar v$ and
$w_n\rightarrow \bar w$.  We have to consider two cases. In each case we
will prove the existence of a
parallel Jacobi field along $\gamma_{\bar w}$ of length at least $2k$.\\
\paragraph{In case} $\bar v\neq\bar w$ in the limit we get a flat strip of length $2k$ 
between $\gamma_{\bar v}$ and $\gamma_{\bar w}$. This corresponds to a Jacobi 
field along $\gamma_{\bar w}$ that is parallel on $[-k,k]$.\\
\paragraph{In case} $\bar v=\bar w$ we have to use a more elaborated
proof. Refer to Figure~\ref{wild} for the following definitions.
\begin{figure}[h]\caption{}\label{wild}
  \psfrag{an}{$a_n=\sigma_n(0)$} \psfrag{bn}{$b_n=\rho_n(0)$}
  \psfrag{sn}{$\sigma_n$} \psfrag{rn}{$\rho_n$} \psfrag{snt}{$\sigma_n(t)$}
  \psfrag{rnt}{$\rho_n(\frac{l_b^n}{l_a^n}t)$} \psfrag{gvn}{$\gamma_{v_n}$}
  \psfrag{gwn}{$\gamma_{w_n}$}
  \psfrag{gwnmk}{$\gamma_{w_n}(-k)=\sigma_n(l_a^n)$}
  \psfrag{gwnk}{$\gamma_{w_n}(k)=\rho_n(l_b^n)$} \psfrag{ln}{$\lambda_t^n$}
  \psfrag{alphan-k}{$\alpha^n_{-k}$}
  \psfrag{alphank}{$\alpha^n_k$}
  \includegraphics[width=\textwidth]{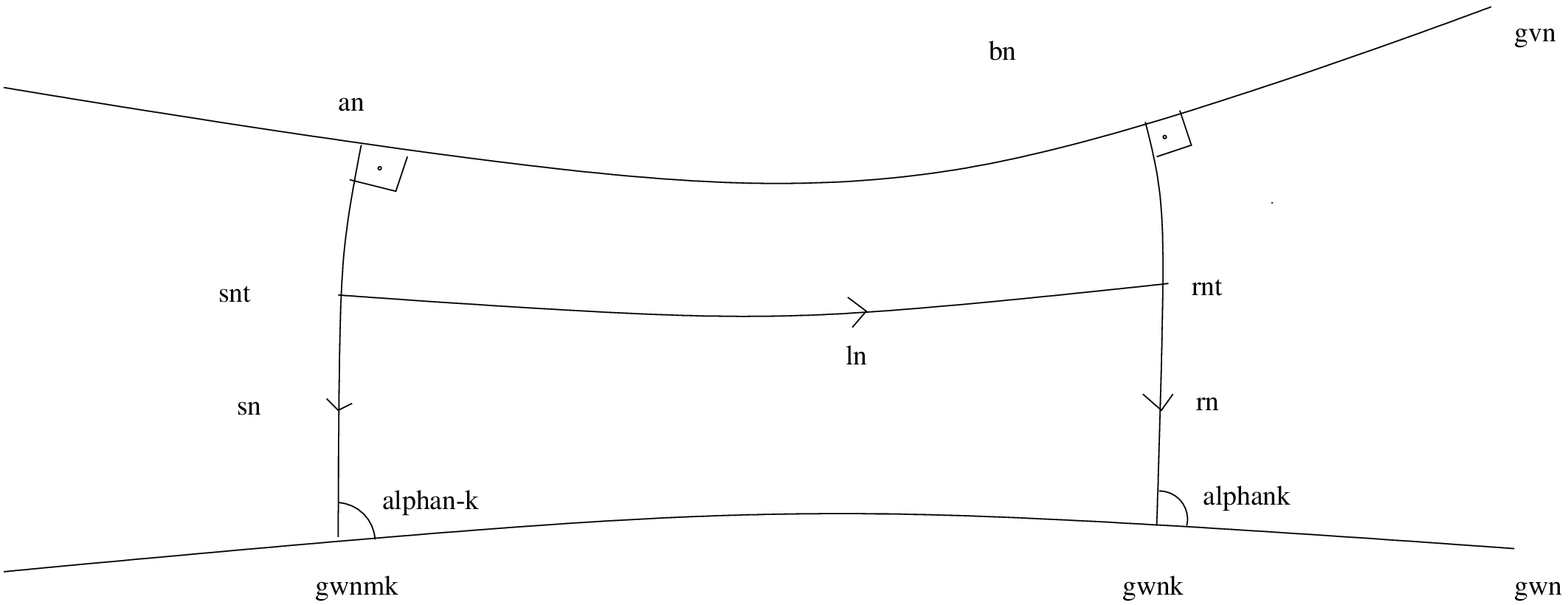}
\end{figure}

In this case let $l_a^n:=d(\gamma_{w_n}(-k),\gamma_{v_n})=\Hd_nf_n(-k)$ and
$l_b^n:=d(\gamma_{w_n}(k),\gamma_{v_n})=\Hd_nf_n(k)$ and write $a_n$ and
$b_n$ respectively for the points on $\gamma_{v_n}$ closest to
$\gamma_{w_n}(-k)$ and $\gamma_{w_n}(k)$ respectively. Call
$\sigma_n:[0,l_a^n]\rightarrow X$ and $\rho_n:[0,l_b^n]\rightarrow X$ the
geodesics linking $a_n=\sigma_n(0)$ to $\gamma_{w_n}(-k)=\sigma_n(l_a^n)$
and $b_n=\rho_n(0)$ to $\gamma_{w_n}(k)=\rho_n(l_b^n)$ respectively. For
$t\in[0,l_a^n]$ let $\lambda_t^n$ denote the geodesic from
$\sigma_n(t)=\lambda_t^n(0)$ to
$\rho(\frac{l_b^n}{l_a^n}t)=\lambda_t^n(1)$. We get a geodesic variation by
defining

\begin{align*}
  h:[0,1]\times[0,l_a^n]&\longrightarrow X\\
  (s,t)&\longmapsto \lambda_t^n(s).
\end{align*}

Consider
$L^n(t):=L(\lambda_t^n)=d(\sigma_n(t),\rho_n(\frac{l_b^n}{l_a^n}t))$.
Applying the first variation formula we see that $\dot L^n(0)=0$, due to
the right angles in $a_n$ and $b_n$, and
$$\dot L^n(1)=\cos(\alpha_{-k}^n)+\cos(\pi/2-\alpha_{k}^n)=\dot
L_n(-k)-\dot L_n(k)\overset{(\ref{dotL})}\leq 4\delta \frac{1-\mu}{n-k}.$$
By the Mean Value Theorem
there must be a $t^n\in[0,l_a^n]$ with 
\begin{align}\label{ddot L}
\ddot L^n(t^n)&=\frac{\dot
  L^n(1)-\dot L^n(0)}{1-0}\leq 4\delta\frac{1-\mu}{n-k}.
\end{align}
Now take a look at the
Jacobi field of the variation $h$ along the geodesic $\lambda_{t^n}^n$.
This is
$$J^n(s):=\frac{\partial}{\partial
  t}h(s,t)\Big|_{t=t^n}=\frac{\partial}{\partial
  t}\lambda_t^n(s){\Big|}_{t=t^n}.$$
For $s=0$ we have
$\lambda_t^n(0)=\sigma_n(t)$ and therefore $J^n(0)=\dot\sigma_n(t^n)$ and
$||J^n(0)||=1$. Notice that $J^n(0)$ is almost
orthogonal because $\alpha_{-k}^n$ equals almost $\pi/2$.\\
We apply the second variation formula to get
\begin{align*}
4\delta\frac{1-\mu}{n-k}\overset{(\ref{ddot L})}\geq \ddot
  L^n(t^n)&=\frac{1}{||\dot\lambda_{t^n}^n||}\int_0^1\left<\frac{\nabla
      {J^n}^\perp}{dt},\frac{\nabla {J^n}^\perp}{dt}\right>-
      \drunter{\qquad\qquad\quad\leq 0 \;(\text{ since } K\leq 0)}
              {\left<R({J^n}^\perp,\dot\lambda_{t^n}^n) 
               \dot\lambda_{t^n}^n,{J^n}^\perp
               \right> } 
      dt\\
  &\geq\frac{1}{\|\dot\lambda_{t^n}^n\|}\int^1_0{\left\|\frac{\nabla
        {J^n}^\perp}{dt}\right\|}^2dt.
\end{align*}
This yields\footnote{ use
\begin{align*}
  ||\dot\lambda^n_{t^n}||&=d(\sigma_n(t^n),\rho_n(\frac{l_b^n}{l_a^n}t^n))\\
  &\leq
  d(\sigma_n(t^n),\gamma_{w_n}(-k))+d(\gamma_{w_n}(-k),\gamma_{w_n}(k))+d(\gamma_{
    w_n}(k),\rho_n(\frac{l_b^n}{l_a^n}t^n))\\
  &\leq d(a_n,\gamma_{w_n}(-k))+2k+d(\gamma_{w_n}(k),b_n)\leq 2k+4\delta.
\end{align*}
}

$$0\leq\int_0^1{\left\|\frac{\nabla {J^n}^\perp}{dt}\right\|}^2dt\leq
4\delta\|\dot\lambda_{t^n}^n\|\frac{1-\mu}{n-k}\leq
8\delta(k+2\delta)\frac{1-\mu}{n-k}\overset{n\rightarrow\infty}\longrightarrow 0.$$

Now, taking a subsequence for which $w_n$ and $v_n$ converge to $\bar w$,
obviously $\lambda_{t^n}^n$ converges to $\gamma_{\bar w}$ and $J^n$
converges to an orthogonal Jacobi field of length 1 along $\gamma_{\bar w}$
that is
parallel on $[-k,k]$.\\

Thus the construction can be completed in any case for any $k\in\mathbb N$
and we obtain a $\bar w(k)\in\overline{\W_{2\delta}(z)}$ with an orthogonal, parallel
Jacobi field $J_k$ along $\gamma_{\bar w(k)}\Big|_{[-k,k]}$. Choosing a
converging subsequence of the $\bar w(k)\rightarrow w'\in \W_{3\delta}(z)$,
we can find a subsubsequence such that the Jacobi fields $J_k$ converge to
a Jacobi field along $\gamma_{w'}$ that must be orthogonal and globally
parallel.
Since $w'\in\W_{3\delta}(z)$ is of rank one, this is a contradiction.
\end{propp}%rank 1 hyper

\begin{rema}{}{}{}
This proof is a variation of the preprint version of~\cite[Lemma
4.1.]{whoknows}. The published proof is shorter. Since we are interested in
procompact sets and not in single vectors we stick to this version of the proof.
\end{rema}

In the rest of this section $X$ will denote a Hadamard manifold with compact
quotient $M=X/\Gamma$. All the results remain true even for Hadamard
manifolds without compact quotient if the condition `$\Gamma$-compact' is
replaced by `compact'.\\

Geodesic rays originating in the same point diverge faster than in the flat
case, if one of the rays is of rank one. Proposition~\ref{radial hyper}
quantifies this behaviour for subsets of rank one.

\begin{prop}{}{radial hyper}{}                %Prop radial hyper
  Let $K\subset UX$ be a $\Gamma$-compact subset consisting only of
  rank one vectors.  Given $N\in\mathbb N$ there are constants $\delta',a>0$
  such that the following holds for any
  $\epsilon<\delta'$:\\
  For any point $o\in X$ and any radial vectors $v\in K$ and $w\in UX$ with
  origin $o$ the inequality $d(\gamma_v(a),\gamma_w)\leq \epsilon$ implies
  $d(p(v),\gamma_w)<\epsilon/N$.
\end{prop}

For the proof of Proposition~\ref{radial hyper} we need to control the
tracing neighbourhoods.
Lemma~\ref{radial trace} provides a simple criterion for radial vectors
by which to decide whether they are within tracing distance.
This lemma is an easy
consequence of Lemma~\ref{close radial vectors}%

\begin{lemm}{}{radial trace}{}                 %Lemm radial trace
  Let $K\subset UX$ be $\Gamma$-compact and $T$, $a$, $\delta$ given. If $\theta$
  is a tracing distance for $K$ with respect to $a$, $\delta$, then there
  is a $\delta'>0$ with the following property: If $o\in X$ satisfies
  $d(o,p(K))>T$ and $v\in K$ and $w\in UX$ are radial with origin $o$ then
  the two inequalities $d(\gamma_v(a),\gamma_w(a))<\delta$ and
  $d(p(v),p(w))<\delta'$ imply that $w$ lies in the tracing neighbourhood
  of $v$ (i.\ e.\ $w\in \W_{\theta}(v)$).
\end{lemm}

\begin{propp}{radial hyper}{}{}                %Prop radial hyper
  For technical reasons assume $N>2$. Fix any $\mu\in]0,1/N[$. Choose
  $\delta$ and $T=T(\mu,\delta)$ for $K$ as in Proposition~\ref{rank 1
    hyper}, i.\ e.\ every vector $v\in K$ satisfies the hyperbolicity
  condition for $\mu$, $T(\mu,\delta)$, $\delta$. Choose
  $a>\frac{N-1}{1-N\mu}T>2T$ and let $\theta<(\frac{1}{1-N\mu}-1)T$ be a
  tracing distance for $K$ with respect to the constants $a$, $\delta$. Let
  $\delta'$ be the constant from Lemma \ref{radial trace} for $K$, $T$,
  $a$, $\delta$. Now fix $\epsilon<\delta'$.\\
  Given an origin $o\in X$ and radial vectors $v\in K$ and $w\in UX$
  suppose $d(\gamma_v(a),\gamma_w) <\epsilon$.

\paragraph{In case $d(p(v),o)\leq T+\epsilon$} let $o=\gamma_v(-r)$ where
$r\leq T+\epsilon\leq T+\theta<\frac{1}{1-N\mu}T$. Applying the convexity
of $g(t):=d(\gamma_v(t),\gamma_w)$ 

\begin{align*}
d(p(v),\gamma_w)=g(0)&\leq g(-r)+\frac{r}{r+a}(g(a)-g(-r))\\
&=\frac{r}{r+a}g(a)\leq\frac{r}{r+a}\epsilon<\frac{\epsilon}{N}
\end{align*}

where the last inequality is due to $r\leq\frac{1}{1-N\mu}T<\frac{a}{N-1}$
which is equivalent to $\frac{r}{r+a}<\frac{1}{N}$. 

\paragraph{In case $d(p(v),o)>T+\epsilon$} define $s$ by
$d(\gamma_v(a),\gamma_w(s))=d(\gamma_v(a),\gamma_w)<\epsilon$. Then

\begin{align*}
  d(\gamma_w(s),o)&\geq d(\gamma_v(a),o)-d(\gamma_v(a),\gamma_w(s))\\
  &> a+T+\epsilon-\epsilon=a+T.
\end{align*}

Thus, setting $w':=\dot\gamma_w(s-a)$ we get a situation as illustrated in
Figure~\ref{dreieck}.

\begin{figure}[h]\caption{}\label{dreieck}
  \psfrag{-T}{$-T$} \psfrag{a}{$a$} \psfrag{0}{0} \psfrag{v}{$v$}
  \psfrag{w}{$w'$} \psfrag{o}{$o$} \includegraphics{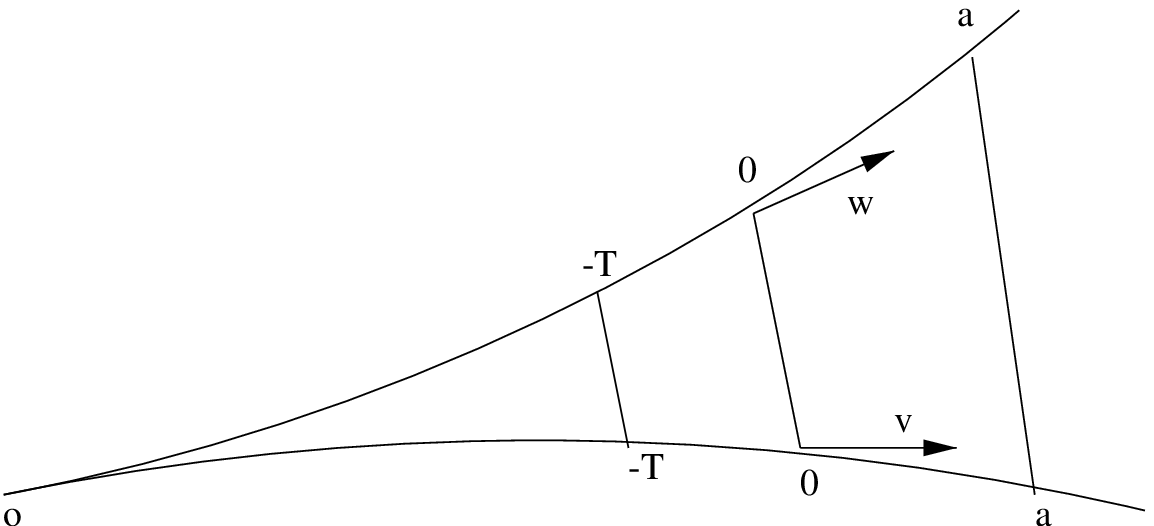}
\end{figure}

Obviously
$d(p(v),p(w'))<d(\gamma_v(a),\gamma_{w'}(a))=d(\gamma_v(a),\gamma_w)<\epsilon<\delta'$
and therefore, by Lemma \ref{radial trace}, $w'$ lies in the tracing
neighbourhood of $v$. We can now apply the hyperbolicity of $v$ to conclude
that
$$d:=d(\gamma_v([-T,T]),\gamma_{w'}([-T,T]))\leq\mu\,
\Hd(\gamma_v([-T,T]),\gamma_{w'}([-T,T]))=:\mu\,\Hd$$
Obviously the Hausdorff
distance is
$\Hd=d(\gamma_v(T),\gamma_{w'}(T))<d(\gamma_v(a),\gamma_{w'}(a))<\epsilon$
and the minimal distance can be estimated by $d\geq
d(\gamma_v(-T),\gamma_{w})=d(\gamma_v(-T),\gamma_{w'})$.\\
So $d(\gamma_v(-T),\gamma_{w})\leq \mu\,\epsilon$ and applying the convexity
of $g(t):=d(\gamma_v(t),\gamma_w)$ we find

\begin{align*}
  d(p(v),\gamma_w)=g(0)&\leq g(-T)+\frac{T}{T+a}(g(a)-g(-T))\\
  &=\frac{a}{T+a}d(\gamma_v(-T),\gamma_w)+\frac{T}{T+a}d(\gamma_v(a),\gamma_w)\\
  &\leq
  \frac{a}{T+a}\mu\,\epsilon+\frac{T}{T+a}\epsilon=\frac{a\mu+T}{T+a}\epsilon<\frac{\epsilon}{N}
\end{align*}

Where the last inequality is due to the fact that $a>\frac{(N-1)T}{1-N\mu}$
is equivalent to $\frac{a\mu+T}{a+T}<\frac{1}{N}$.
\end{propp}

We will see that in fact we can choose $\delta'$ arbitrarily. I.\ e.\ given
a $\delta$ we can choose $a$ big enough so that $\delta$ and $a$ satisfy
the conditions for the constants in Proposition~\ref{radial hyper}. To
prove this we need a topological lemma which is based on the properties of
the distance function.\\

\begin{lemm}{}{nice distance}{}                 %Lemm nice distance
  Fix two different points $p,o\in X$ on a Hadamard manifold (not
  necessarily with compact quotient). Consider the
  map
\begin{align*}
  f:U_oX&\longrightarrow \mathbb R\\
  w&\longmapsto d(p,\gamma_w)
\end{align*}
Write $\Delta:=d(o,p)$. Call $\bar v\in U_oX$ the tangent vector to the
geodesic ray starting in $o$ and passing through $p$. For $\eta>0$ we find
\begin{enumerate}
\item $f\leq\Delta$, hence $f\inv(\eta)=\emptyset$ for $\eta>\Delta$
\item if $\eta=\Delta$ then $f\inv(\eta)=\{w\in U_oX\big|w\perp\bar v\}$.
\item if $\eta<\Delta$ then $f\inv(\Delta)=\{w\in U_oX\big|w\perp\bar v\}$
  divides $U_oX$ into two connected components (disks). Consider the one
  containing $\bar v$ and call it $C$. This is the halfsphere facing $p$.
  The set $f\big|_C\inv(\eta)$ is a topological sphere and divides $C$ into
  two connected components, which are characterised by $\Delta>f>\eta$ (an
  annulus) and $0<f<\eta$ (a disk) respectively.
\end{enumerate}
Obviously in all cases $f(\bar v)=0$ and for $w\perp\bar v$ the maximum
$f(w)=\Delta$ is attained.
\end{lemm}

\begin{lemmp}{nice distance}{}{}
  The first two properties are clear. Now consider the third case where we
  are interested in $f\big|_C$. We will show that the gradient vanishes
  only in $\bar v$.\\
  For $w\in C$ let $p_w$ be the point on $\gamma_w$ closest to $p$. Call
  $N$ the set of all these points. $N$ is a $(\dim X-2)$-dimensional
  submanifold of $X$, diffeomorphic to $C$. Obviously for every $w\in C$ we
  know $f(w)=d(p_w,p)$.\\
  Now fix any point $p_w\in N\backslash p$. Consider the geodesic $\sigma$
  in $X$ with $\sigma(0)=p_w$ and $\sigma(1)=p$. For small $t$ let $w_t\in
  U_oX$ denote the vector tangent to the geodesic ray from $o$ to
  $\sigma(t)$. Then $w_t$ is a differentiable curve in $C$. Take a look at
  $f(t):=f(w_t)$. Note that $f(t)=d(\gamma_{w_t},p)\leq
  d(\sigma(t),p)=\Delta(1-t)$. So
  $$f'(0)=\underset{t\rightarrow
    0}\lim\frac{f(t)-f(0)}{t}\leq\underset{t\rightarrow
    0}\lim\frac{\Delta(1-t)-\Delta}{t}=-\Delta<0$$
  and therefore the gradient does not vanish in $p_w$.\\
  Now we have a gradient on the disk $C$ that vanishes only in one point
  $\bar v$, which is the minimum, and therefore the geodesic flow on this
  disk has a simple structure resulting in the property quoted.
\end{lemmp}

It is easy to see that this implies the following lemma.

\begin{lemm}{}{gleich reicht}{}     %Lemm: gleich reicht
  Suppose there are $v\in UX$, $\eta_1,\eta_2>0$ and $a>0$
  such that for any radial vector $w\in UX$ with origin
  $o\in\gamma_v(\mathbb R_-)$ we have
  $$d(p(v),\gamma_w)=\eta_1\quad\Longrightarrow\quad
  d(\gamma_v(a),\gamma_w)>\eta_2.$$
  \emph{Claim:} Then the implication
  $$d(p(v),\gamma_w)\geq\eta_1\quad\Longrightarrow\quad
  d(\gamma_v(a),\gamma_w)>\eta_2.$$
  holds, too. Or equivalently:
  $$d(p(v),\gamma_w)<\eta_1\quad\Longleftarrow\quad d(\gamma_v(a),\gamma_w)
  \leq\eta_2$$
\end{lemm}

Now we have the tools to prove an improved version of Proposition~\ref{radial hyper}. In Proposition~\ref{radial hyper} for a given
$N\in\mathbb N$ we get the existence of constants $a$ and $\delta$. We will
see now, that in fact we have one more degree of freedom, since for any
given $\delta$ we can find an $A$ guaranteeing the desired property.\\

\begin{prop}{}{better hyper}{}        %Prop: better hyper
  Let $K\subset UX$ be a $\Gamma$-compact subset consisting only of rank
  one vectors.
  Given a shrinking factor $N\in\mathbb N$ and a shield radius $\delta>0$ there is a distance
  $A>0$ such that the following holds for any $\epsilon<\delta$:\\
  For any point $o\in X$ and any radial vectors $v\in K$ and $w\in UX$ with
  origin $o$ the inequality $d(\gamma_v(A),\gamma_w)\leq \epsilon$ implies
  that $d(p(v),\gamma_w)<\epsilon/N$ holds.
\end{prop}

\begin{propp}{better hyper}{}{}
  By Proposition~\ref{radial hyper} we can find $\delta'$ and $a$ with
  the desired property. If $\delta$ is less or equal than $\delta'$ we are
  done. So suppose $\delta>\delta'$. Define
  $A:=\frac{\frac{\delta}{\delta'}N-1}{N-1}a$. Given $\epsilon<\delta$,
  define $\epsilon':=\frac{\epsilon}{\delta}\delta'<\delta'$. Proposition~\ref{radial hyper} applies to $\epsilon'$ and so for all radial vectors
  $v\in K$, $w\in UX$ with common origin $o\in X$
  $$d(p(v),\gamma_w)=\frac{\epsilon'}{N}=:\eta_1\quad\Longrightarrow\quad
  d(\gamma_v(a),\gamma_w)>\epsilon'.$$
  By convexity we get
\begin{align*}
  d(\gamma_v(A),\gamma_w)&\geq
  d(p(v),\gamma_w)+\frac{A}{a}\left(d(\gamma_v(a),\gamma_w)-d(p(v),\gamma_w)\right)\\
  &=\frac{a-A}{a}d(p(v),\gamma_w)+\frac{A}{a}d(\gamma_v(a),\gamma_w)\\
  &>\left(1-\frac{\frac{\delta}{\delta'}N-1}{N-1}\right)\frac{\epsilon'}{N}+
  \left(\frac{\frac{\delta}{\delta'}N-1}{N-1}\right)\epsilon' \\
  &=\frac{\delta}{\delta'}\epsilon'=\epsilon=:\eta_2.
\end{align*}
Now Lemma~\ref{gleich reicht} yields the desired implication:
$$d(\gamma_v(A),\gamma_w)\leq\epsilon\quad\Longrightarrow\quad
d(p(v),\gamma_w)<\frac{\epsilon'}{N}<\frac{\epsilon}{N}$$
\end{propp}

\begin{figure}[h]\caption{}\label{Bild best hyper}
  \psfrag{gammaw}{$\gamma_w$}
  \psfrag{gammav}{$\gamma_v$}
  \psfrag{gammavA}{$\gamma_v(A)$}
  \psfrag{A}{$=A$}
  \psfrag{r0}{$\geq r_0$}
  \psfrag{o}{$o$}
  \psfrag{w}{$w$}
  \psfrag{v}{$v$}
  \psfrag{K}{$K$}
  \psfrag{abstand}{$d(\gamma_v(A),\gamma_w)$}
  \psfrag{implikation}
  {$d(\gamma_v(A),\gamma_w)\leq\epsilon<\Delta
    \quad\Longrightarrow \quad
    d^1(v,\dot\gamma_w)<\frac\epsilon{N}$}

  \includegraphics[width=\textwidth]{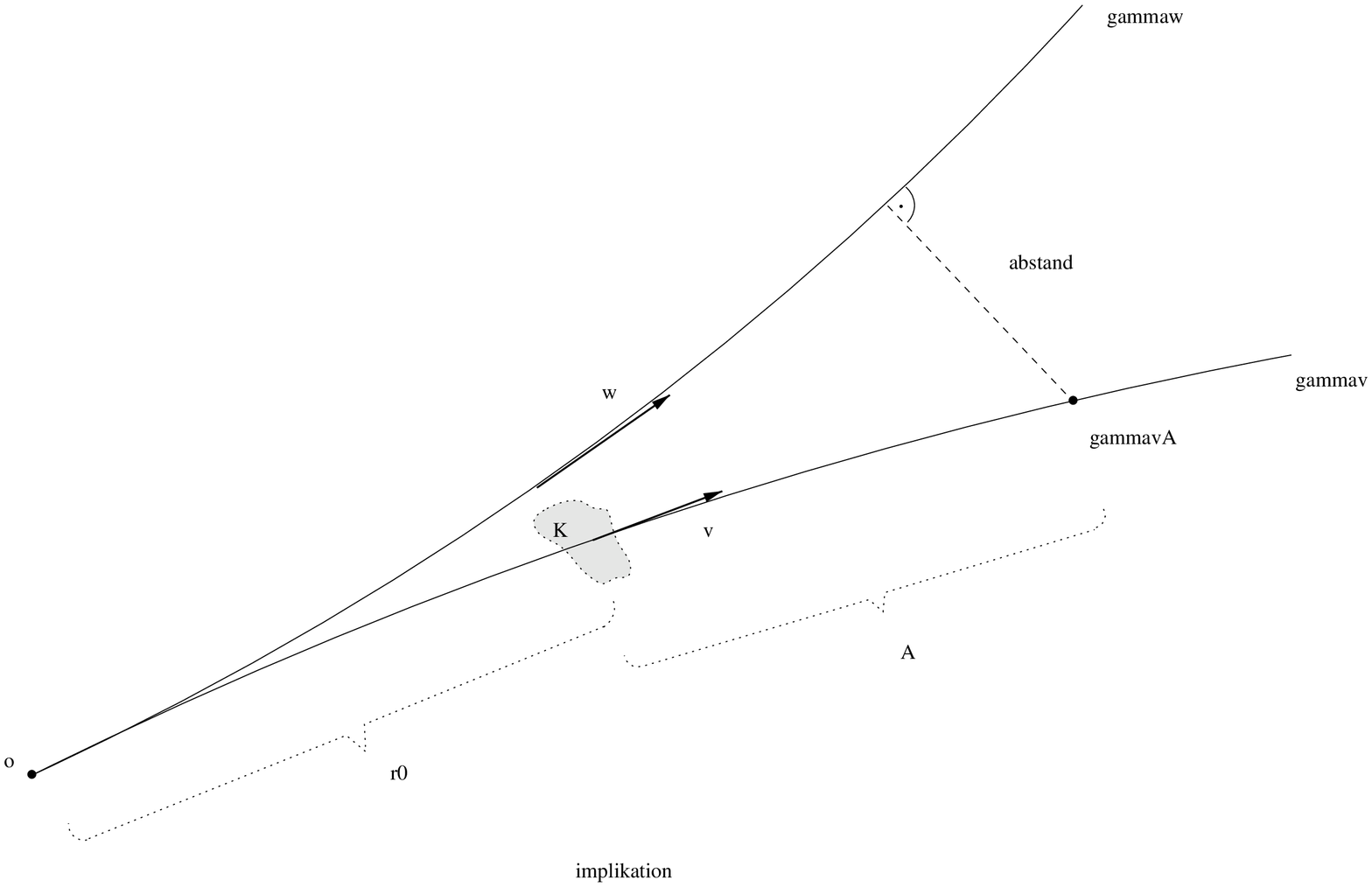}
\end{figure}

Combining Proposition~\ref{better hyper} with Lemma~\ref{close radial
  vectors} we get Corollary~\ref{best hyper} which is illustrated in
Figure~\ref{Bild best hyper}.

\begin{coro}{}{best hyper}{}
  Let $K\subset UX$ be a $\Gamma$-compact subset consisting only of rank
  one vectors.
  Given a shrinking factor $N\in\mathbb N$, a shield radius $\Delta>0$ and
  a radius $r_0>0$ there is a distance
  $A>0$ such that the following holds for any $\epsilon<\Delta$:\\
  For any point $o\in X$ and any radial vectors $v\in K$ and $w\in UX$ with
  origin $o$ and $d(p(v),o)\geq r_0$ the inequality $d(\gamma_v(A),\gamma_w)\leq \epsilon$ implies
  that $d^1(v,\dot\gamma_w)<\epsilon/N$ holds.
\end{coro}

\mysection{Avoiding a Rank One Vector}{}\label{rank1.tex}
This section explains the construction of Sergei Buyalo and Viktor
Schroeder in~\cite{whoknows}. They proved that for any given vector $v_0\in
UM$ of rank one on a compact manifold $M$ of nonpositive curvature and of
dimension $\geq 3$ it is possible to construct a closed, flow invariant,
full subset $Z\subset UM$ which does not contain $v_0$. Furthermore, for
every $\epsilon>0$ there is such a set $Z_\epsilon$ which is
$\epsilon$-dense in $UM$.\\
The obvious way to construct this set, is to fix some small $\eta>0$ and
define
$$Z:=\left\{ v\in UM\;|\;d(\dot\gamma_v,v_0)\geq\eta\right\}$$
but for all we know, this set could be empty or neither full nor
$\epsilon$-dense. So how do we prove that this set is full, provided $\eta$
is small enough?\\
Pick any point $o\in M$ and show that, for $\eta$ small enough, we can find
many geodesics passing through $o$ and avoiding an $\eta$-neighbourhood
of $v_0$. These geodesics can be explicitly constructed.

\mysubsection{Idea}{}
Consider a compact connected manifold $M$ of nonpositive curvature, a rank
 one  vector $v_0\in UM$ and a point $o\in M$. Let $\pi:X\rightarrow M$ be the
universal covering and choose one element in $\pi\inv(o)$ which we will
identify with $o$.\\
We want to find a geodesic in $M$ that avoids the direction $v_0$ in the unit tangent bundle of $M$. In the
universal covering this corresponds to avoiding the set
$\Omega:=d\pi\inv(v_0)\subset UX$.
\index{Omega@$\Omega$}\index{$omega$@$\Omega$}
We can identify this set with the set of its
base points $\Sigma:=p(d\pi\inv(v_0))=\pi\inv(p(v_0))$
\index{Sigma@$\Sigma$}\index{$sigma$@$\Sigma$} 
via the base point projection. Elements of $\Sigma$ will be denoted with
$\omega$ 
\index{omega@$\omega$}\index{$omega$@$\omega$}
and the corresponding vector in $\Omega$ with
$v_\omega$.
\index{v(omega)@$v_\omega$}\index{$v(omega)$@$v_\omega$}  
Thus $p(v_\omega)=\omega$. For
technical reasons in Proposition~\ref{geodesic} we would like $\Omega$ to
be symmetric. Hence suppose that $\pm v_\omega\in\Omega$ holds for all
$\omega\in\Sigma$.\\ 
Obviously, a geodesic that avoids some
neighbourhood of an $\omega$, avoids an even
bigger neighbourhood of the vector $v_\omega$.
($d(\gamma,\omega)>\eta\Rightarrow d(\dot\gamma,v_\omega)>\eta$). Hence we will
work in $X$ rather than in $UX$ and try to avoid
$\Sigma$. But indeed we only need to consider elements of
$\Sigma$ which correspond to vectors $v_\omega$ which are almost
radial with respect to our starting point $o$. So for small
$\tilde\epsilon>0$ ($\tilde\epsilon$ will be specified later on) we define the subsets

$$\Sigma(o,\tilde\epsilon):=\left\{\omega\in\Sigma\;|\;d\left(v_\omega,\phi_{\mathbbR_+}(U_oX)\right)<\tilde\epsilon\right\}
\index{Sigma(o,epsilontilde)@$\Sigma(o,\tilde\epsilon)$}
\index{$sigma(o,epsilontilde)$@$\Sigma(o,\tilde\epsilon)$}$$
To every $\omega\in\Sigma(o,\tilde\epsilon)$ we now assign a region in $X$
that should not be hit (the protective disk) and a region that is save to
hit (the good annulus) and a sphere inside this annulus that is far away
from the boundary of the annulus (the good sphere).\\ 
 \\
A ray originating in $o$ that comes close to an obstacle\index{obstacle}
$\omega\in\Sigma(o,\tilde\epsilon)$
will be replaced by a close ray passing the good sphere.

\begin{figure}[h]\caption{How to Avoid Obstacles}
\label{Scheiben}
  \psfrag{o1}{$\omega_1$} \psfrag{g0}{$\gamma_0$} \psfrag{Ro}{$R_\omega$}
  \psfrag{o2}{$\omega_2$} \psfrag{g1}{$\gamma_1$} \psfrag{So}{$S_\omega$}
  \psfrag{o3}{$\omega_3$} \psfrag{g2}{$\gamma_2$} \psfrag{Do}{$D_\omega$}
  \psfrag{o}{$\omega$} \psfrag{g3}{$\gamma_3$}
  \includegraphics{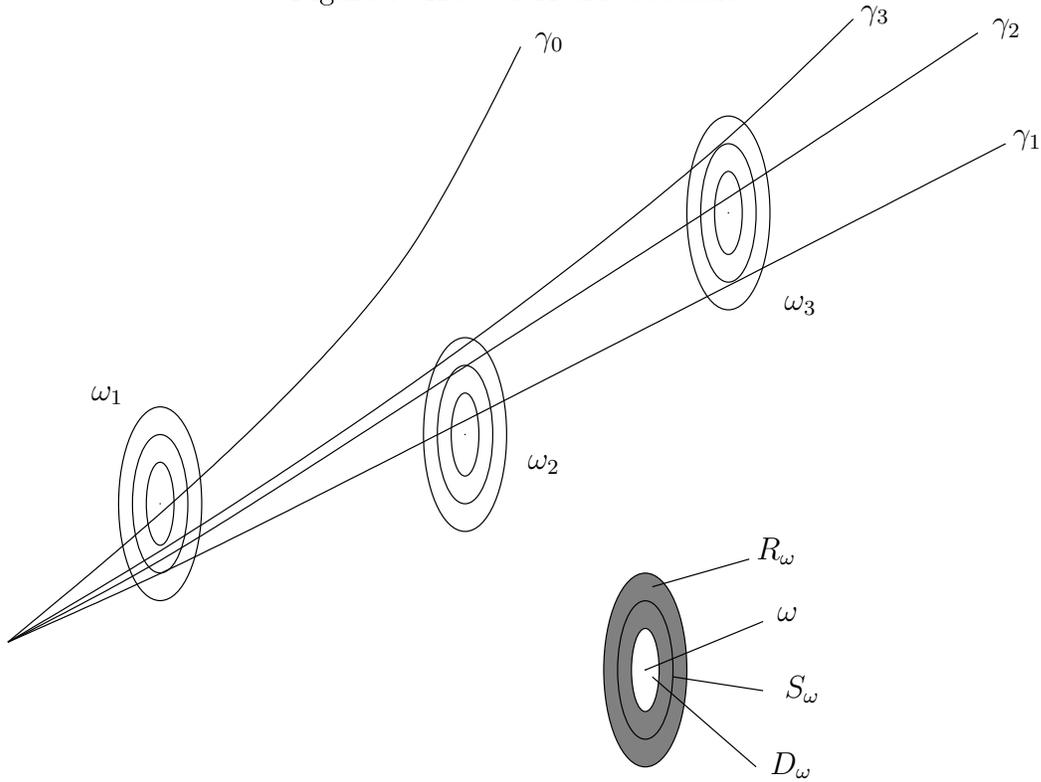}
\end{figure}

Proposition \ref{far away} will show that for small displacements it holds: 
If the new ray hits another protective screen $D_{\omega'}$ then any
ray hitting the new annulus $R_{\omega'}$ will have passed through the old
annulus $R_\omega$ before. This is illustrated in Figure~\ref{Scheiben}.

\begin{defi}{}{}{}          %Defi protective screen
  For $v\in U_oX$ and $\omega\in\Sigma(o,\tilde\epsilon)$ define
  $p(\omega,v)$, to be the point of minimal distance to $\omega$ on the ray
  $\gamma_v([0,\infty[)$. For this distance we write
  $d_\omega(v):=d(p(\omega,v),\omega)=d(\gamma_v,\omega)$.
\index{d omega@$d_\omega(.)$}\index{$d omega$@$d_\omega$}
  \bigskip\\
  Consider a small $\delta'$ (it will be specified later on). Choose a
  $\delta<\delta'/6$ and define for every
  $\omega\in\Sigma(o,\tilde\epsilon)$ the \emph{protective
    screen}\index{protective screen}\index{screen!protective}\index{disk!protective} (disk)
  $$D_\omega:=\{p(\omega,v)|v\in U_oX \text{ und } d_\omega(v)<2\delta\},
\index{D omega@$D_\omega$}\index{$d omega$@$D_\omega$}$$
  the \emph{good annulus}\index{annulus!good}\index{good annulus}
  $$R_\omega:=\{p(\omega,v)|v\in U_oX \text{ und }
  d_\omega(v)\in]\delta,3\delta[\}
\index{R omega@$R_\omega$}\index{$r omega$@$R_\omega$}
$$
  and contained in the good annulus, the
  \emph{good sphere}\index{good sphere}\index{sphere!good}
  $$S_\omega:=\{p(\omega,v)|v\in U_oX \text{ und } d_\omega(v)=2\delta\}.
\index{S omega@$S_\omega$}\index{$s omega$@$S_\omega$}
$$
\end{defi}

\mysubsection{Choice of Constants}{}

Choose a neighbourhood ${\cal U}_0$ of $v_0$ such that it consists of rank
 one  vectors only and such that its inverse image under $d\pi$ consists of
disjoint sets ${\cal U}_\omega$.%
\index{$u omega$@${\cal U}_\omega$}\index{U omega@${\cal U}_\omega$}
Choosing the neighbourhood small enough,
we can suppose that it is contained in a compact set with the same
properties. The preimage of this set under $d\pi\inv$ is a \Gam-compact set
in $UX$. 
For this \Gam-compact set and $N=6$ Proposition~\ref{radial
  hyper} yields constants $\delta',a$ where we may assume that
$\delta'<a/2$. We deduce that the property in Proposition~\ref{radial
  hyper} holds for the subset $d\pi\inv{\cal U}_0$ with the same constants
$\delta',a,N=6$.\\
\\
In fact we want $\delta'$ to satisfy two further conditions.
\begin{enumerate}%1
\item Since ${\cal U}_0$ is open, we can find an $\eta>0$ such that
  $\W_{\eta}(v_0)\subset{\cal U}_0$ and thus
  $\W_\eta(v_\omega)\subset{\cal U}_\omega$ for all $\omega$.\\
  By Lemma~\ref{close radial vectors} we may suppose that $\delta'$ was
  chosen so small that, if we put $\delta:=\delta'/6$, for any
  $\omega\in\Sigma(o,\tilde\epsilon)$ satisfying $d(\omega,o)\geq a$ all
  the radial rays originating in $o$ and hitting the protective disk or the
  good annulus of $\omega$ are within $\eta$-distance of
  $v_\omega$ and therefore are in ${\cal U}_\omega$ and of rank one.\\
  So from now on assume, that for $d(o,\omega)>a$, all radial vectors with
  origin $o$ hitting $D_\omega$ or $R_\omega$ are rank one
  vectors which have the property described in Proposition~\ref{radial hyper}.\\
\item Now consider $\Omega_a:=\Sigma(o,\tilde\epsilon)\cap{\cal U}_{a}(o)$.
  This is the finite set of points in $\Sigma(o,\tilde\epsilon)$ that are
  within distance less than $a$ of our starting point $o$. For every
  $\omega\in\Omega_a$ call $\theta_\omega\in U_oX$ the starting direction
  of the geodesic ray from $o$ to $\omega$; and for $\eta>0$ define
  $$V_\eta(\omega):=\{v\in U_oX\,|\, {d}%frueher d^1
(v_\omega,\dot\gamma_v(\mathbb
  R_+))<\eta\}.$$
  This is an open subset of $U_oX$ and $\diam
  V_\eta(\omega)\overset{\eta\rightarrow 0}{\longrightarrow}0$. Choose
  $\rho<\underset{\omega_1,\omega_2\in\Omega_a}{\min}{d}%frueher d^1
(\theta_{\omega_1},\theta_{\omega_2})$
  and find an $\eta>0$ so small that $\diam V_\eta(\omega)<\rho/5$ for all
  $\omega\in\Omega_a$. Now the sets
\begin{align*}
  U_1&:=\{v\in U_oX\big|{d}%frueher d^1
(v,\theta_\omega)<\frac{1}{5}\rho\text{ for all
    }\omega\in\Omega_a\} \text{ and } \\ 
  U_2&:=\{v\in U_oX\big|{d}%frueher d^1
(v,\theta_\omega)<\frac{2}{5}\rho\text{ for all
    }\omega\in\Omega_a\}
\end{align*}
have the following properties:
\begin{enumerate}%2
\item $\{\theta_\omega\}_{\omega\in\Omega_a}\subset U_1\subset U_2$
\item $v\in U_oX\backslash U_1\quad\Rightarrow\quad
  {d}%frueher d^1
(v_\omega,\dot\gamma_v(\mathbb R_+))\geq\eta$ for all
  $\omega\in\Omega_a$
\item $U_oX\backslash U_2$ is closed and a sphere with finitely many holes
  of diameter $\frac{4}{5}\rho$.
\end{enumerate}%2
Now choose $\delta'<a\sin\frac{\rho}{10}$.
\end{enumerate}%1

\mysubsection{Construction}{}
Start with a vector $x_0\in U_oX\backslash U_2$. Consider the geodesic ray
$\gamma_{x_0}([a,\infty[)$. Call $\omega_1\in\Sigma(a,\tilde\epsilon)$ the
first obstacle\footnote{here first obstacle means that $d(\omega_1,o)$ is
  minimal} with $d_{\omega_1}(x_0)<2\delta$. Move $x_0$ to a vector $x_1\in
U_oX$ with $p(\omega_1,x_1)\in S_{\omega_1}$, i.\ e.\ the geodesic ray
$\gamma_{x_1}$ passes through the good sphere of $\omega_1$.\\
Now define recursively a sequence of vectors $x_i\in U_oX$ and elements
$\omega_i\in\Sigma(o,\tilde\epsilon)$ as follows: The first time
$\gamma_{x_i}([d(p(\omega_i,x_i),o,\infty[)$ hits the protective disk of an
element of $\Sigma(o,\tilde\epsilon)$, call this element $\omega_{i+1}$ and
move $x_i$ to a vector $x_{i+1}\in U_oX$ where $\gamma_{x_{i+1}}$ passes
the good sphere of $\omega_{i+1}$. We get a sequence $\{x_i\}$ of vectors
in $U_oX$ with the following properties.
\begin{enumerate}
\item $d(\omega_i,\gamma_{x_i})=d(\omega_i,p(\omega_i,x_i))=2\delta$
\item $d(\omega_{i+1},\gamma_{x_i})=d(\omega_{i+1},p(\omega_{i+1},x_i))<2\delta$
\item $d(\omega_{i+1},\omega_i)>a+\delta'$
\item $d(\omega_i,o)\geq i\frac{a}{2}$.
\end{enumerate}

\mysubsection{The Construction Has the Desired Properties}{Properties}

We need another technical definition:
\begin{defi}{}{}{}
A subset $\Sigma\subset X$ is called \emph{radially
  $(a,\delta)$-separated}
\index{radially separated}
\index{separated!radially} 
with respect to $o\in X$, if for any vector $v\in U_oX$ and any two points
$\omega_1,\omega_2\in \Sigma$  
$$\max\left\{d(\gamma_v,\omega_1),d(\gamma_v,\omega_2)\right\}\leq
\delta \quad \Longrightarrow\quad d(\omega_1,\omega_2)>a.$$
\end{defi}

In~\cite[Lemma 4.3]{whoknows} it is proved for our setting, that if $\gamma_{v_0}$ is
nonperiodic then for any $a>0$ we can find $\delta$ and $\tilde\epsilon$ small
enough such that $\Sigma(o,\tilde\epsilon)$ is radially
$(a+\delta,\delta)$-separated for all $o\in X$.\\
If $\gamma_{v_0}$ is periodic of period smaller than $a$, the construction
is more difficult. But still we can find a radially separated set
$\Sigma^N(o,\tilde \epsilon)$
\index{$sigma N(o,epsilon)$@$\Sigma^N(o,\tilde\epsilon)$}
\index{Sigma N(o,epsilon tilde)@$\Sigma^N(o,\tilde \epsilon)$}
if we consider the set $\Omega^N$ instead of $\Omega$ (consider only every $N^\text{th}$ element of $\Omega$ along any
geodesic $\gamma_{v_\omega}$. For more details see~\cite{whoknows}). For
the next few pages we will suppose that
$\Sigma(o,\tilde\epsilon)$ is radially separated, but we will return to
the problem of $\Sigma^N(o,\tilde \epsilon)$ at the end ot this section.\\

Notice that the geodesic ray $\gamma_{x_j}$ passes through all
the good annuli of the $\omega_i$ for $i<j$. This is proved in a slightly
more general form in the following proposition.

\begin{prop}{}{far away}{}                    %Prop far away
  Suppose $\Sigma(o,\tilde\epsilon)$ is radially
  $(a+\delta',\delta')$-separated with respect to $o$. Choose
  $\delta<\delta'/6$. Let $\omega\in\Sigma(o,\tilde\epsilon)$ and let
  $v\in U_oX$ denote the initial direction of a geodesic ray, hitting the
  good sphere $S_\omega$ in $p(\omega,v)$. Call
  $\omega'\in\Sigma(o,\tilde\epsilon)$ the first obstacle with
  $d(\omega',o)>d(\omega,o)$ for which the ray hits the protective disk
  $D_{\omega'}$  (i.\ e.\ $d_{\omega'}(v)<2\delta$). Then:\\
  
  Any ray that meets the annulus $R_{\omega'}$ has passed the annulus
  $R_\omega$ before.\\
  Phrased differently this means for any $x\in U_oX$ the following
  implication holds:
  $$\delta<d_{\omega'}(x)<3\delta\qquad \Longrightarrow \qquad\delta<d_\omega(x)<3\delta$$
\end{prop}

\begin{propp}{far away}{}{}
  We know $d(p(\omega,v),\omega)=2\delta<\delta'$ and
  $d(p(\omega',v),\omega')<2\delta<\delta'$. Applying the
  $(a+\delta',\delta')$-separatedness $d(\omega,\omega')>a+\delta'$
  and thus
  $$d(p(\omega',v),p(\omega,v))\geq
  d(\omega,\omega')-d(\omega,p(\omega,v))-d(\omega',p(\omega',v))\geq
  a+\delta'-2\delta-2\delta\geq a.$$
  
  For a vector $x\in U_oX$ with $d(p(\omega',x),\omega')<3\delta$ it
  follows $d(p(\omega',v),\gamma_x)\leq
  d(p(\omega',x),p(\omega',v))<3\delta+2\delta<\delta'$. Now
  Proposition~\ref{radial hyper} can be applied and we find a point $q$
  on $\gamma_x$ 
  satisfying
  $$
  d(p(\omega,v),q)<\frac{1}{6}d(p(\omega',v),p(\omega',x))<\frac{1}{6}5
  \delta=\delta
  $$
  and therefore
  $$d(p(\omega,x),\omega)\leq d(q,\omega)\leq
  d(q,p(\omega,v))+d(p(\omega,v),\omega)<3\delta$$
  
  Suppose now $d(p(\omega,x),\omega)<\delta$. By a similar argument there
  is a $q''$ on $\gamma_v$ satisfying $d(p(\omega,x),q'')<\delta$ and
  therefore $d(q'',\omega)<2\delta=d(\gamma_v,\omega)$ in contradiction to
  the minimal property of $p(\omega,v)$.
\end{propp}

The next lemma will show that all the $x_i$ respect the minimal distance of
$\delta$ to these $\omega$ if
$d(\omega_i,o)>d(\omega,o)$.\\

\begin{lemm}{}{global avoid}{}                    %Lemm global avoid
  For every  $\omega\in\Sigma(o,\tilde\epsilon)$ with
  $d(\gamma_{x_i},\omega)>2\delta$ we have
  $d(\gamma_{x_j},\omega)>\delta$ for all $j>i$.
\end{lemm}

\begin{lemmp}{global avoid}{}{}
  Suppose $d(p(\omega,x_j),\omega)\leq\delta$. As we know, $\gamma_{x_j}$
  hits the good annulus of $\omega_{i+1}$. Applying the radial
  separatedness we get
\begin{align*}
  d(p(\omega,x_j),p(\omega_{i+1},x_j))&\geq
  d(\omega,\omega_{i+1})-d(\omega,p(\omega,x_j))-d(\omega_{i+1},p(\omega_{i+1},x_j))\\
  &\geq a+\delta'-\delta-3\delta>a.
\end{align*}
Thus there is a $q$ on $\gamma_{x_i}$ such that
$$d(q,p(\omega,x_j))\leq
\frac{1}{6}d(p(\omega_{i+1},x_j),p(\omega_{i+1},x_{i}))\leq\frac{5}{6
  }\delta<\delta,$$
hence $d(\gamma_{x_i},\omega)\leq d(q,\omega)<\delta
+\delta=2\delta$, a contradiction.
\end{lemmp}

Now we have found a sequence of geodesic rays avoiding more and more
obstacles. In fact this sequence converges to a geodesic ray avoiding all
the obstacles:

\begin{lemm}{}{convergence}{}                     %Lemm convergence
  For every vector $x_0\in U_oX\backslash U_2$ the sequence
  $\{x_i\}_{i\in\mathbb N}$ converges to a vector $z\in U_oX\backslash
  U_1$.
\end{lemm}

\begin{lemmp}{convergence}{}{}
  Suppose the sequence is not finite. We know for $j\geq i\geq 1$ that
  $d(\omega_i,\gamma_{x_j})<3\delta$ since $\gamma_{x_j}$ passes the good
  annulus of $\omega_i$.  Since $d(\omega_i,o)>i\frac{a}{2}$ and the angle
  $\measuredangle(o,p(\omega_i,x_j),\omega_i)$ is $\frac{\pi}{2}$, we know
  
  $$\sin(\measuredangle(\theta_{\omega_i},x_j))\leq
  \frac{3\delta}{i\frac{a}{2}}=\frac{6\delta}{ia}\overset{i\rightarrow\infty}{\longrightarrow}0$$
  
  and thus $\{x_j\}$ is a Cauchy sequence converging to a vector $z\in
  U_oX$. Notice that
  $$\measuredangle(\theta_{\omega_1},z)=\underset{j\rightarrow\infty}{\lim}
  \measuredangle(\theta_{\omega_1},x_j)\leq\arcsin\left(\frac{6\delta}{a}\right)<\frac{\rho}{10}$$
  and since $\gamma_{x_0}$ hits the protective disk we know

  $$\measuredangle(\theta_{\omega_1},x_0)\leq\arcsin\left(\frac{2\delta}{\frac{a}{2}}\right)
  <\frac{\rho}{10}.$$
  Hence $z$ is within $\frac{\rho}{5}$-distance of
  $x_0$ and therefore in $U_oX\backslash U_1$.
\end{lemmp}

So, by mapping $x_0$ to $z$, we can define a map
$$\Phi:U_oX\backslash U_2 \longrightarrow U_oX\backslash U_1$$
where every
vector $z=\Phi(x_0)$ in the image of $\Phi$ satisfies:
\begin{enumerate}
\item The speed vectors $\dot\gamma_z$ of the geodesic ray $\gamma_z$ respect a distance of $\tilde\epsilon$
  (with respect to ${d}%frueher d^1
(.,.)$) to all $v_\omega$ with
  $\omega\in\Sigma\backslash\Sigma(o,\tilde\epsilon)$.
\item The geodesic ray $\gamma_z$ respects a distance of $\delta$ to all
  $\omega\in\Sigma(o,\tilde\epsilon)\backslash\Omega_a$.
\item The geodesic ray $\gamma_z$ respects a distance of $\eta$ to all
  $\omega\in\Omega_a$.
\end{enumerate}
Furthermore $\Phi$ satisfies
$$d(x_0,\Phi(x_0))<\frac{\rho}{5}\text{ for all }x_0\in
U_oX\backslash U_2.$$
Notice that we can choose $\rho$ arbitrarily small
and still find
appropriate $\eta$, $\delta<\delta'/6$ and $\tilde\epsilon$.\\

\begin{lemm}{}{manifold}{}
  Fix $\rho>0$. For any submanifold $Y\subset U_oX\backslash U_2$ we can
  define a continuous map
  $$\Phi_Y:Y\rightarrow U_oX\backslash U_1$$
  and an $\eta>0$ such that the
  image $\Phi_Y(Y)$ consists of initial directions of geodesic rays that
  miss the neighbourhood
  $\underset{\omega\in\Omega}{\bigcup}\W_\eta(v_\omega)$ of $\Omega$ and
  such that
  $${d}%frueher d^1
(\Phi_Y(y),y)<\frac{\rho}{5}\text{ for all }y\in Y.$$
\end{lemm}

\begin{lemmp}{manifold}{}{}
  Instead of moving every starting vector individually and looking at the
  limit we now consider its neighbourhood in $Y$. So consider an obstacle
  $\omega$.  Notice that $\gamma_Y\cap D_\omega$ is an at most $\dim
  Y$-dimensional submanifold of $D_\omega$ which may therefore be deformed
  continuously into $S_\omega$ keeping $\gamma_Y\cap S_\omega$ fixed. This
  deformation in $\gamma_Y\cap D_\omega$ can be projected down to $U_oX$
  and we get a series of deformations of $Y$ with a continuous limit
  deformation $\Phi_Y$.
\end{lemmp}

Now we are in the position to prove our main proposition:

\begin{prop}{}{geodesic}{}
  If $M$ is a compact manifold of nonpositive curvature and $v_0\in UM$ is
  a vector of rank one, then for any $o\in M$ we can find a geodesic
  $\gamma:\mathbb R\rightarrow M$ and an open neighbourhood $\cal U$ of
  $v_0$ in $UM$ such that $\gamma(0)=o$ and $\gamma$ does not meet $\cal
  U$.
\end{prop}

\begin{propp}{geodesic}{}{}
  Of course we work in the universal covering $X$ of $M$ as before. Use all
  the notations as introduced in this section.\\
  Fix one vector $S\in U_oX\backslash U_2$, the south-pole. Take a
  $\dim(X)-2$-dimensional subspace $W$ of $T_SU_oX$ and its one dimensional
orthogonal complement $W^\perp$. Define two submanifolds% 
\footnote{in fact these are two great spheres in $U_oX$.}
of $U_oX$ of dimensions $\dim(X)-2$ and $1$ respectively by $Y:=\exp_SW$
and $Y^\perp:=\exp_SW^\perp$. Because $U_oX\backslash U_2$ looks like a
sphere with small holes of diameter $\frac{2}{5}\rho$, we can choose a
deformation $\phi:U_oX \rightarrow U_oX$ that deforms $Y$ and $-Y^\perp$
into two submanifolds $Y_1$ and $-Y_1^\perp$ of $U_oX\backslash U_2$ and
every point of $U_oX$ is moved by at most $\frac{2}{5}\rho$.\\
Now apply Lemma~\ref{manifold} to $Y_1$ and $-Y_1^\perp$ and deform them
into simplices in $U_oX$ homotopic to $Y$ and $-Y^\perp$, respectively.
Notice that $\Phi_{Y_1}\circ\phi$ and $\Phi_{-Y^\perp_1}\circ\phi$ move
every vector by less than $\frac{3}{5}\rho$. Since $Y$ and $Y^\perp$
intersect in two antipodal points and the intersection number of homotopic
simplices is a constant, $\Phi_{Y_1}(\phi(Y))$ and
$-\Phi_{-Y_1^\perp}(\phi(-Y^\perp))$ still
intersect in two points, provided $\frac{6}{5}\rho<\pi$.\\
Take one of these intersection points $z$. Since $z\in\Phi_{Y_1}(\phi(Y))$
we know that $\gamma_z(\mathbb R_+)$ avoids a neighbourhood of $v_0$ and
since $-z\in\Phi_{-Y_1^\perp}(\phi(-Y))$ so does $\gamma_z(\mathbb R_-)$.
\end{propp}

\begin{rema}{}{}{}
\begin{itemize}
\item In fact all the constants were chosen globally. So all the geodesic we
constructed avoid the same neighbourhood of $v_0$. So the set 
$$ Z:=\left\{v\in UM\;|\;d(\dot\gamma_v,\Omega)\geq\eta\right\}$$
is \emph{full}.
\item Take a closer look at the proof of
  Proposition~\ref{geodesic}: Given any vector $v\in U_oX$  there is a
  vector $S\in U_oX\backslash\U_2$ that is at most $(\frac25\rho)$-far from
  $v$. The construction yields a vector $v'$ such that $\dot\gamma_{v'}$ avoids an $\eta$-neighbourhood of
  $\Omega$ and $v'$ is $(\frac65\rho)$-close
  to $S$. Thus $Z$

is \emph{$(\frac85\rho)$-dense} in $UM$.

\item As remarked  before, if $\gamma_{v_0}$ is $L$-periodic we have to
  modify the proof. Proposition~\ref{geodesic} holds for $\Sigma^N(o,\tilde
  \epsilon)$ instead of $\Sigma(o,\tilde \epsilon)$. Now suppose $\dot\gamma_v$
  avoids an $\eta$-neighbourhood of $\Omega^N$. But if
  $\gamma_v$ came very close, say $\nu$, to $v_\omega\in\Omega$ then it
  would be close to $\gamma_{v_\omega}$ for a long time, say $T$. But if
  this time is much bigger than $N$ times the period of $\gamma_{v_0}$ then
  it would pass at least one of the obstacles in $\Omega^N$very
  closely. Therefore $\nu$ can not be arbitrarily small 
  and Proposition~\ref{geodesic} holds for periodic $v_0$, too.
\end{itemize}
\end{rema}

\mysection{Avoiding Subsets of Small s-Dimension}{}\label{ausweichen.tex}

Consider a Riemannian manifold $X$ and a closed submanifold $R\subset
X$ of smaller dimension $\dim(R)<\dim(X)$. In this section we discuss deformations of $UX\backslash UR$ such that
every vector is moved away from $UR$.\\
$UR$ is a submanifold of the Riemannian manifold $UX$. So if there is a
$2\epsilon$-neighbourhood $\W_{2\epsilon}(UR)$ of $UR$ on which the
gradient $\grad d_{UR}$ of the distance function $d_{UR}:=d(UR,.)$ is defined, we
can use the flow $\phi$ generated by $\grad d_{UR}$ to deform $UX\backslash
UR$ into $UX\backslash\W_{\epsilon}(UR)$. Notice that the gradient and
hence $\phi$ might not be defined outside $\W_{2\epsilon}(UR)$ but still
the deformation
$$\Phi(v):=\phi_{\max\{0,2\epsilon-\frac12d_{UR}(v)\}}(v)$$
is well defined and moves every vector by at most $2\epsilon$ into
$UX\backslash\W_{\epsilon}(UR)$.\\

We need, however, a more elaborated deformation, since it will be
essential for us to have control on the displacement of the base points,
too. Furthermore we want our deformation to respect the foliation of $UX$
into spheres $\sphere_v(r)$ of fixed radius $r$.\\
In Subsection~\ref{UTB} we will discuss a deformation which gives control on
$d_{UR}$ and on
the base point distance from $R$ at the same time. Based on this deformation we will find another
deformation in Subsection~\ref{rSF} that respects the sphere foliation of
$UX$.\\
Finally in Corollary~\ref{subanalytischesausweichen} we generalize the
deformation to the case where $R$ is a stratified subset of $X$:

\paragraph{Corollary~\ref{subanalytischesausweichen} }{$ $}\\
\begin{it}
Let $X$ denote a Hadamard manifold with compact quotient $M=X/\Gamma$.
Suppose a \Gam-compact stratified subset $R\subset X$ and $r_0>0$ are given.\\

Then there is a constant $k>0$ such that for small $\lambda>0$ and given $r>r_0$, $v\in UX$
and manifold $Y$ with
$$\dim Y<\dim X -\dim R$$
we can deform any
continuous map $c:Y\rightarrow\sphere_r(v)$ into a continuous map
$c_\lambda :Y\rightarrow \sphere_r(v)$ such that $c_\lambda$ is
$k\lambda$-close to $c$ and $c_\lambda$ avoids a $\lambda$-neighbourhood of
$UR$, i.\ e.\ for all $x\in Y$ 
$$d(c(x),c_\lambda(x))\leq k\lambda\qquad\text{ and }\qquad d_{UR}(c_\lambda(x))\geq\lambda.$$
\end{it}

\mysubsection{Avoiding the Unit Tangent Bundle of a Submanifold}
{Avoiding a Submanifold}\label{UTB}

Suppose $X$ is a Hadamard manifold with compact quotient $M=X/\Gamma$ and
$R\subset X$ is a $\Gamma$-compact submanifold of $X$.
We will try to avoid its unit tangent bundle $UR$ in $UX$. Use the shorthand notation $d_R(.):=d(R,.)$%
\index{$d r(.)$@$d_R(.)$}\index{d R@$d_R(.)$}
for the
distance to $R$ in $X$ and $d_{UR}:=d(UR,.)$%
\index{$d UR$@$d_{UR}$}\index{d UR@$d_{UR}$}
for the distance in $UX$ to $UR$.\\

If $R$ is a submanifold with boundary, notice that we can find a
slightly larger, $\Gamma$-compact submanifold $R'$ in $X$ with $\dim
R'=\dim R$, $R\subset R'$ and $d(R,\partial R')>\rho$ for some $\rho>0$. If
$R$ is a closed manifold without boundary, set $R'=R$. The proof still works.\\
Suppose furthermore that $\rho>0$ is so small that 
\begin{itemize}
\item $\rho$ is smaller than the injectivity radius of $M$
\item $\rho$ is smaller than the injectivity radius of $UM$
\item $\rho$ is so small that $\grad_pd_{R'}$ is defined for all $p\in X$
  with $0<d_{R'}(p)<\rho$.\\
Hence $\gamma_p:=\gamma_{\grad_pd_{R'}}$ realizes the distance from ${R'}$ to
$\gamma_p(t)$ for all $t\in[-d_{R'}(p),\rho-d_{R'}(p)]$.
\end{itemize}

Fix $\epsilon<\rho$.\\
We define the distortion $\Pi_\epsilon$
\index{$pi epsilon$@$\Pi_\epsilon$}\index{Pi epsilon@$\Pi_\epsilon$}
of $UX\backslash p\inv {R'}:=\{v\in UX\;\big|\;p(v)\notin {R'}\}$.
To shorten notation we define the distortion length (depending on $\epsilon$)
$$\begin{array}{ccccc}\index{alpha@$\alpha$}\index{$alpha$@$\alpha$}
\alpha:& UX\backslash p\inv {R'}&\longrightarrow&[0,2\epsilon]\\
       & v                   &\longmapsto    &
\begin{array}\{{c}.0\\
             2\epsilon-\frac12{d_{U{R'}}(v)}
\end{array}&
\begin{array}{c} \text{ if }d_{U{R'}}(v)\geq 4\epsilon\\\text{ if }d_{U{R'}}(v)\leq 4\epsilon.\end{array}
\end{array}$$
The distortion is defined by
$$\begin{array}{cccc}
\Pi_\epsilon:&UX\backslash p\inv {R'}&\longrightarrow&UX\backslash p\inv {R'}\\
    &v&\longmapsto&_{\gamma_{p(v)}}\big\|_0^{\alpha(v)}v.
\end{array}$$
This is well defined, since $\alpha(v)=0$ for those $v$ where the gradient
of $d_{R'}$ is not defined in $p(v)$.

\begin{prop}{}{distortion}{}
The distortion $\Pi_\epsilon$ has the following properties if we choose $\epsilon<\frac\rho4$:
\begin{enumerate}
\item $d(v,\Pi_\epsilon(v))=\alpha(v)$ for all $v\in UX\backslash p\inv {R'}$.\\
Furthermore $d(v,_{\gamma_{p(v)}}\big\|_0^tv)=t$ for all $t\in[0,\alpha(t)]$.
\item \label{aussage ueber dR}
$d_{R'}(p(\Pi_\epsilon(v)))=d_{R'}(p(v))+\alpha(v)>0$ for all $v\in UX\backslash
  p\inv {R'}$.\\
Furthermore $\gamma_{p(v)}$ is a shortest from ${R'}$ to $\gamma_{p(v)}(t)$
for all $t\in[0,\alpha(t)]$. 

\item For all $v\in UX\backslash p\inv {R'}$ we have
$$\begin{array}{ccccc}
d_{U{R'}}(v)\geq 4\epsilon&\Longleftrightarrow&d_{U{R'}}(\Pi_\epsilon(v))\geq
4\epsilon&\Longleftrightarrow&\Pi_\epsilon(v)=v\\
d_{U{R'}}(v)<4\epsilon&\Longleftrightarrow&d_{U{R'}}(\Pi_\epsilon(v))<4\epsilon.
\end{array}$$
\item $\Pi_\epsilon$ is continuous, injective and a homeomorphism onto the image.
\item The image of $\Pi_\epsilon$ is contained in $UX\backslash\W_\epsilon U{R'}$,\\
i.\ e.\ $d_{U{R'}}(\Pi_\epsilon(v))\geq\epsilon$ for all $v\in UX\backslash p\inv {R'}$.
\item $\Pi_\epsilon$ is compatible with the \Gam-action on $X$ and $UX$.
\end{enumerate}

\end{prop}

\begin{propp}{distortion}{}
Recall that for any $v\in UX\backslash p\inv {R'}$ the curve 
$$\sigma_v:t\longmapsto _{\gamma_{p(v)}}\big\|_0^tv$$
is a geodesic in $UX$. Since the projection $p\circ\sigma_v=\gamma_{p(v)}$
is a shortest from ${R'}$ to $\gamma_{p(v)}(t)$ for
$t\in]-d_{R'}(p(v)),\rho-d_{R'}(p(v))]$ it does not intersect ${R'}$ and hence
$\sigma_v(t)\notin U{R'}$ for these $t$.\\
We conclude that $\phi_v:t\mapsto d_{U{R'}}(\sigma_v(t))$ is
differentiable for these $t$ and $\phi_v'(t)\in[-1,1]$. 

\begin{enumerate}
\item %%%%%%%%%%%%%%%%%%%%%%%%%%Beweis 1%%%%%%%%%%%%%%%%%%%%%%%%%%%
This holds since $\sigma_v$ is a geodesic and
  $\alpha(v)$ is smaller than the injectivity radius of $UX$. 
\item %%%%%%%%%%%%%%%%%%%%%%%%%%Beweis 2%%%%%%%%%%%%%%%%%%%%%%%%%%%
$\gamma_{p(v)}$ realises the distance to ${R'}$ for $t-d_{R'}(p(v))<\rho$.
Since
  $p(\Pi_\epsilon(v))=p(_{\gamma_{p(v)}}\big\|_0^{\alpha(v)}v)=\gamma_{p(v)}(\alpha(v))$ we get $d_{R'}(p(\Pi_\epsilon(v)))=d_{R'}(p(v))+\alpha(v)$.

\item %%%%%%%%%%%%%%%%%%%%%%%%%%Beweis 3%%%%%%%%%%%%%%%%%%%%%%%%%%%
The implications
$$ \Pi_\epsilon(v)=v\;\;\Longleftarrow\;\;d_{U{R'}}(v)\geq 4\epsilon\;\;\Longrightarrow\;\;
d_{U{R'}}(\Pi_\epsilon(v))\geq 4\epsilon$$
are obvious. We need to show
$$d_{U{R'}}(v)<4\epsilon\;\;\Longrightarrow\;\; d_{U{R'}}(\Pi_\epsilon(v))<4\epsilon.$$
So take $v\in UX\backslash p\inv {R'}$ with $d_{U{R'}}(v)<4\epsilon$. Then
\begin{align*}
d_{U{R'}}(\Pi_\epsilon(v))&\leq d_{U{R'}}(v)+d(v,\Pi_\epsilon(v))\\
              &= d_{U{R'}}(v)+\alpha(v)\\
              &=2\epsilon+\frac12\drunter{<4\epsilon}{d_{U{R'}}(v)}<4\epsilon.
\end{align*}

\item %%%%%%%%%%%%%%%%%%%%%%%%%%Beweis 4%%%%%%%%%%%%%%%%%%%%%%%%%%%%
$\Pi_\epsilon$ is continuous since $\alpha$ is.\\
Suppose we find $v_1,v_2\in UX\backslash p\inv {R'}$ with
$z:=\Pi_\epsilon(v_1)=\Pi_\epsilon(v_2)$.\\
\begin{itemize}
\item Suppose $d_{U{R'}}(z)\geq 4\epsilon$ and therefore
  $d_{U{R'}}(v_1),d_{U{R'}}(v_2)\geq 4\epsilon$. But
  $v_1=\Pi_\epsilon(v_1)=z=\Pi_\epsilon(v_1)=v_2$ for these vectors. 

\item
Now suppose $d_{U{R'}}(z)< 4\epsilon$ and hence $d_{U{R'}}(v_i)<
  4\epsilon$ for both $i$ and $\alpha(v_i)=2\epsilon-\frac12d_{U{R'}}(v_i)$.\\
We use the notation $p_i:=p(v_i)$, $\gamma_i:=\gamma_{p_i}$, $\sigma_i:=\sigma_{v_i}$,
$\phi_i:=\phi_{v_i}$.\\
Without loss of generality we suppose 
\begin{equation}\label{deltas}
d_{R'}(p_1)\leq d_{R'}(p_2).
\end{equation}
By Part~\ref{aussage ueber dR} of the
proposition
\begin{equation}\label{alphas}
d_{R'}(p_1)+\alpha(v_1)=d_{R'}(z)=d_{R'}(p_2)+\alpha(v_2)
\end{equation}
and by condition (\ref{deltas})
$\alpha(v_1)-\alpha(v_2)=d_{R'}(p_2)-d_{R'}(p_1)\geq 0$ which translates into
\begin{equation}\label{1kleiner2}
d_{U{R'}}(v_1)\leq d_{U{R'}}(v_2).
\end{equation}
By definition of $\Pi_\epsilon$ the points $p_i$ and $q:=p(z)$ lie on the unique
shortest geodesic from ${R'}$ to $q$. And the vectors $v_i$ and $z$ can be
obtained via parallel transport along that geodesic $\sigma_1$. So
\begin{align}
\sigma_1(0)&=v_1\label{sigma1}\\ 
\sigma_1(d_{R'}(p_2)-d_{R'}(p_1))&=v_2\label{sigma2}\\
\sigma_1(\alpha(v_1))&=z\label{sigma3}
\end{align}
Now recall that $\phi_1(t)$ is differentiable for $t\in[0,\alpha(v_1)]$
with derivative in $[-1,1]$. Therefore
$$d_{U{R'}}(v_2)=\phi_1(d_{R'}(p_2)-d_{R'}(p_1))\leq\drunter{=d_{U{R'}}(v_1)}{\phi_1(0)}+(d_{R'}(p_2)-d_{R'}(p_1))$$
and thence
$$d_{U{R'}}(v_2)-d_{U{R'}}(v_1)\leq d_{R'}(p_2)-d_{R'}(p_1)\overset{(\ref{alphas})}=
\alpha(v_1)-\alpha(v_2)=\frac12(d_{U{R'}}(v_2)-d_{U{R'}}(v_1)).$$
Combining this with (\ref{1kleiner2}) we see that
$d_{U{R'}}(v_1)=d_{U{R'}}(v_2)$, hence $\alpha(v_1)=\alpha(v_2)$ and by
(\ref{alphas}) $d_{R'}(p_2)=d_{R'}(p_1)$. By (\ref{sigma1}) and (\ref{sigma2}) we
get
$$v_1\overset{(\ref{sigma1})}=\sigma_1(0)=\sigma_1(d_{R'}(p_2)-d_{R'}(p_1))\overset{(\ref{sigma2})}=v_2.$$
\end{itemize}
%%%%%%%%%%%%%%%%%%%%%%%%%%%%%homeo%%%%%%%%%%%%%%%%%%%%%%%%%%%%%%%%%
It remains to show that $\Pi_\epsilon\inv$ is continuous. Write $Y:=U{R'}\backslash
p\inv {R'}$. \\
Take any $v\in Y$ and fix $\eta<\frac12d_{R'}(p(v))$. We claim that
$\Pi_\epsilon(\W_\eta(v))$ is open in $\Pi_\epsilon(Y)$.\\
Write $\U$ for $\W_\eta(v)$ and assume that $\Pi_\epsilon(\U)$ is not open.
 Then we find a
sequence $\Pi_\epsilon(v_i)\in \Pi_\epsilon(Y)\backslash\Pi_\epsilon(\U)=\Pi_\epsilon(Y\backslash\U)$ converging to a vector
$\Pi_\epsilon(w)\in\Pi_\epsilon(\U)$. But $\{v_i\}_i$ is a bounded sequence in $Y$ since 
$$d(v_i,w)\leq \drunter{\leq
  2\epsilon}{d(v_i,\Pi_\epsilon(v_i))}+\drunter{\rightarrow
  0}{d(\Pi_\epsilon(v_i),\Pi_\epsilon(w))}+\drunter{\leq 2\epsilon}{d(\Pi_\epsilon(w),w)}$$
and therefore a subsequence converges to a vector $z\in UX\backslash\U$
(closed set) whose base 
point lies not in ${R'}$ since
$$d_{R'}(p(z))\leftarrow d_{R'}(p(v_i))\geq d_{R'}(p(v))-\eta>\frac12d_{R'}(p(v))>0.$$
So $z\in Y\backslash\U$ and $\Pi_\epsilon(z)=\Pi_\epsilon(\lim v_i)=\lim\Pi_\epsilon(v_i)=\Pi_\epsilon(w)$ and
by injectivity $w=z$ which is a contradiction since $w\in\U$ and $z\notin
\U$.

\item %%%%%%%%%%%%%%%%%%%%%%%%Beweis 5%%%%%%%%%%%%%%%%%%%%%%%%%%%%%%
Take $v\in UX\backslash p\inv {R'}$. To see that
  $d_{U{R'}}(\Pi_\epsilon(v))>\epsilon$ we need to distinguish three cases:
  \begin{itemize}
   \item If $d_{U{R'}}(v)\geq 4\epsilon$ then $\Pi_\epsilon(v)=v$ and
     $d_{U{R'}}(\Pi_\epsilon(v))\geq 4\epsilon$.
   \item If $d_{U{R'}}(v)\leq 2\epsilon$ then
   \begin{align*}
   d_{U{R'}}(\Pi_\epsilon(v))&\geq d_{R'}(p(\Pi_\epsilon(v)))\\
                 &=d_{R'}(p(v))+\alpha(v)\\
                 &> \alpha(v)\\
                 &=2\epsilon-\frac12d_{U{R'}}(v)\geq\epsilon.
   \end{align*}
   \item If $d_{U{R'}}(v)> 2\epsilon$ then 
   \begin{align*}
   d_{U{R'}}(\Pi_\epsilon(v))&\geq d_{U{R'}}(v)-d(v,\Pi_\epsilon(v))\\
                 &=d_{U{R'}}(v)-\alpha(v)\\
                 &=\frac32d_{U{R'}}(v)-2\epsilon> \epsilon.
   \end{align*}
  \end{itemize} 
\item This is obvious since the distance functions and hence the gradient
  field are invariant under deck transformations.
\end{enumerate}
\end{propp}

\begin{rema}{}{}{}
The construction uses a displacement of maximally $2\epsilon$ to achieve a
distance to $U{R'}$ of at least $\epsilon$. The quotient of these two values is 
$$\frac{\text{distance}}{\text{displacement}}\;=\;\frac12.$$
Optimal would be a value of $1$ while any value $>0$ might be called
satisfactory. A natural approach to maximize the quotient is to consider a
different distortion length function. For $\lambda\in]0,2[$ we could use 
$\alpha_\lambda:=\lambda\alpha$ (For $\lambda\geq 2$ the construction would
fail to be injective). Then Proposition~\ref{distortion} holds for all
$\epsilon<\frac\rho{2\lambda}$ with maximal displacement
$2\lambda\epsilon$ and minimal distance to $U{R'}$ at least
$\frac{2\lambda}{\lambda+1}\epsilon$. So in this way the quotient is 
$$\frac{\text{distance}}{\text{displacement}}\;=\;\frac1{\lambda+1}$$
and hence for small $\epsilon$ the distortion becomes very effective.
\end{rema}

\mysubsection{Respecting the $r$-Sphere Foliation of $UX$}
{Respecting the $r$-Sphere Foliation}
\label{rSF}
Next we will use the distortion $\Pi_\epsilon$ to define distortions on
$UX\backslash p\inv R'$ that map vectors normal to a sphere to vectors that
are normal to the same sphere. This distortion will be denoted by $\Psi_{r,\epsilon}$
where $r$ is the radius of the sphere.
We will use the following notation:\\
For $p\in R'\backslash\partial R'$ define the \emph{normal tangent space} by
$${\cal N}_pR':=\{u\in T_pX\;|\;u\perp v\text{ for all } v\in T_pR'\}.$$
Notice that $R\subset R'\backslash\partial R'$ and hence we can define the 
\emph{normal bundle} \emph{normal
bundle}\index{normal bundle}\index{bundle!normal} ${\cal N}R\subset
TX$\index{N R@${\cal N}R$}\index{$n r$@${\cal N}R$} of $R$ by
$$ {\cal N}R:=\underset{p\in R}\bigcup {\cal N}_pR'.$$
This is a proper euclidean bundle over $R$ which is independent of the
choice of $R'$.\\

The \emph{unit normal bundle}\index{bundle!unit normal}\index{unit normal
  bundle} is the subset of all vectors of unit length and will
be denoted by ${\cal N}^1R\subset UX$.
\index{N 1 R@${\cal N}^1R$}\index{$n 1 r$@${\cal N}^1R$}
This sphere
bundle over $R$ is \Gam-compact. \\

\begin{figure}[h] \caption{Construction of $\Psi_{r,\epsilon}$}\label{konstruktion}
\psfrag{b}{$\beta$}
\psfrag{r}{$r$}
\psfrag{l}{$l(r):=d(m,q)$}
\psfrag{a}{$\alpha(v)$}
\psfrag{pi}{$\pi-\beta$}
\psfrag{w}{$w_r(v)$}
\psfrag{Pi}{$\Pi_\epsilon(v)$}
\psfrag{v}{$v$}
\psfrag{Psi}{$\Psi_{r,\epsilon}(v)$}
\psfrag{p}{$p=p(v)$}
\psfrag{q}{$q:=p(\Pi_\epsilon(v))$}
\psfrag{q1}{$q':=p(\Psi_{r,\epsilon}(v))$}
\psfrag{m}{$m:=\gamma_v(-r)$}
\psfrag{lr=}{}%{$l(r):=d(m,q)$}
\psfrag{beta=}{$\beta=\measuredangle(v,\grad_pd_{R'})$}
\psfrag{dvpi=}{$d(v,\Pi_\epsilon(v))=d(p,q)=\alpha(v)$}
\psfrag{dpiw=}{$d(\Pi_\epsilon(v),w_r(v))=\measuredangle(\Pi_\epsilon(v),w_r(v))$}
\psfrag{hintenunten}{}
\psfrag{dwpsi=}{$d(w_r(v),\Psi_{r,\epsilon}(v))=d(q,q')=l(r)-r$}
\includegraphics[width=\textwidth]{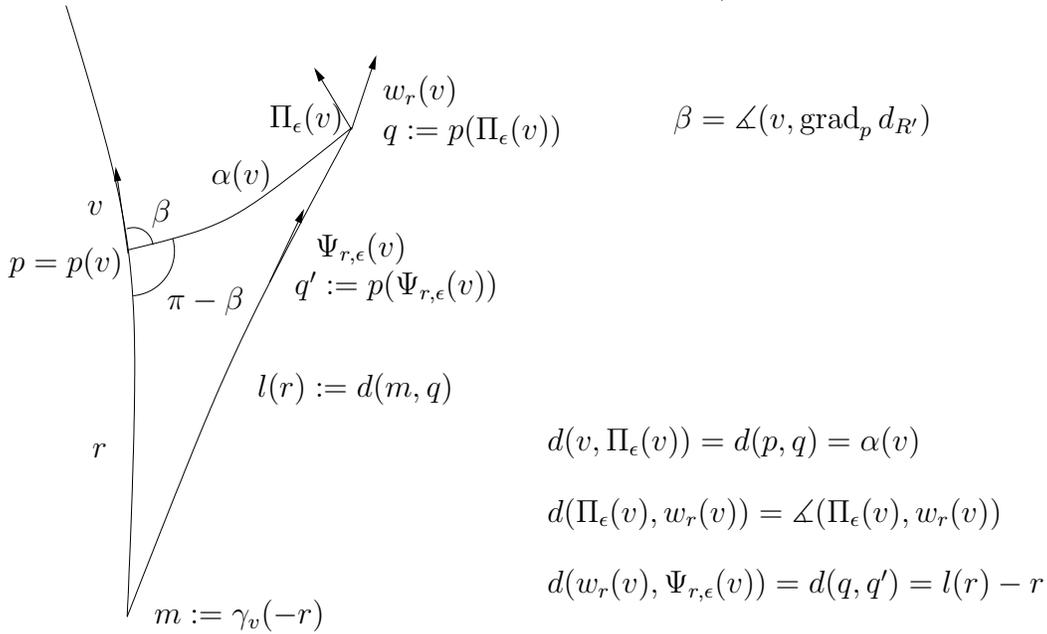}
\end{figure}

For $r>0$ and $v\in UX$ recall the definition of \emph{spheres of radius
  $r$ defined by $v$}
\begin{description}
\item[in $X$:] $S_r(v):=\{q\in X\;\big|\;d(\gamma_v(-r),q)=r\}$ and
\item[in $UX$:] ${\cal S}_r(v):=\{\dot\sigma(r)\in UX\;|\;\sigma\text{ geodesic}\wedge\sigma(0)=\gamma_v(-r)\wedge\sigma(r)\in S_r(v)\}$.
\end{description}
So $S_r(v)$ is the sphere of radius $r$ centered at $\gamma_v(-r)$ (hence
$p(v)\in S_r(v)$) and ${\cal S}_r(v)$ is the set of outward unit vectors to
this sphere (and by construction $v\in {\cal S}_r(v)$). \\
Now fix a radius $r$. We want to define a distortion 
$$\Psi_{r,\epsilon}:UX\backslash p\inv {R'}\longrightarrow UX\text{
  such that }\Psi_{r,\epsilon}(v)\in\sphere_r(v)\text{ for all }v.
\index{$psi r,epsilon$@$\Psi_{r,\epsilon}$}
\index{Psi r epsilon@$\Psi_{r,\epsilon}$}
$$
 The idea of the construction is
illustrated in figure~\ref{konstruktion}.
For any vector $v$ let $m:=\gamma_v(-r)$ denote the
centre of the sphere to which $v$ is a normal vector. Write $\sigma$ for
the geodesic starting in $m$ and passing through $q:=p(\Pi_\epsilon(v))$. Call $q'$
the point on $\sigma$ at distance $r$ from $m$, so $p:=p(v)$ and $q'$ lie
on the sphere around $m$ of radius $r$. Define the vector tangent to
$\sigma$ in $q'$ to be $\Psi_{r,\epsilon}(v)$. We might as well describe the map by
$$\Psi_{r,\epsilon}(v):=\frac{d}{dt}\Big|_{t=r}\exp_{\gamma_v(-r)}\left(t\frac{\exp\inv_{\gamma_v(-r)}(p(\Pi_\epsilon(v)))}{\|\exp\inv_{\gamma_v(-r)}(p(\Pi_\epsilon(v)))\|}\right)$$
which shows that the map is surely continuous for $r>2\epsilon\geq\alpha(v)$.\\
We will see that for small $\epsilon$ and big $r$ this distortion has the
property that the displacement is small and still the distance to $UR$ of
the image is strictly bigger than zero by Proposition~\ref{Psir}. To see this we have to consider the
behaviour of the following functions first:\\

Fix a point $p\in X$ and vectors $v,z\in U_pX$. For $a\in\mathbbR_+$
let $\Pi_{\epsilon,a}(v,z):=_{\gamma_z}\big|_0^av$ denote the vector we get by parallel transport of $v$
along the geodesic segment $\gamma_z\big\|_{[0,a]}$.\\
For $r\in\mathbbR_+$ map
$v$ to the vector $w_{a,r}(v,z)$ that is tangent in $\gamma_z(a)$ to the
geodesic $\sigma$ starting
in $\gamma_v(-r)$ and passing through $\gamma_z(a)$. So for fixed $r$ and $v$ the
vectors $w_{a,r}(v,z)$ are radial vectors originating in $\gamma_v(-r)$ with
base  point on $\gamma_z$.\\

Write $q_{a}(v,z)$ for the base  point $\gamma_z(a)$ of $\Pi_{\epsilon,a}(v,z)$ and
$w_{a,r}(v,z)$ (notice that  $q_a(v,z)$ is independent of $r$). The point on the geodesic $\sigma$ that is at distance $r$
from $\gamma_v(-r)$ will be denoted by $q'_{a,r}(v,z)$. Write
$l_{a}(r,v,z)$ for the distance $d(q_a(v,z),q'_{a,r})$.\\

We can estimate the displacement of $\Psi_{r,\epsilon}$ if we understand the following
displacements:
\begin{align}
v&\rightarrow\Pi_\epsilon(v)=\Pi_{\epsilon,\alpha(v)}(v,\grad_pd_{R'})\label{vnachpi}\tag{Step 1}\\
&\rightarrow w_{\alpha(v),r}(v,\grad_pd_{R'})=:w_r(v)\label{pinachw}\tag{Step 2}\\
&\rightarrow \Psi_{r,\epsilon}(v)\label{wnachpsi}\tag{Step 3}
\end{align}
By Proposition~\ref{distortion} we understand quite well~\ref{vnachpi}. We want $\Psi_{r,\epsilon}$ to inherit
the main properties of this map. So we need to estimate the other two
steps.~\ref{wnachpsi} is parallel transport along the geodesic
$\sigma$. So the displacement is exactly the distance of the base  points
$d(q,q')=|l_r(v)-r|$. The displacement in~\ref{pinachw} is the angle
between the radial vector and the parallely transported vector. We will see
that these two displacements can be controlled if we make sure $\alpha(v)$
is small enough.

\begin{lemm}{\ref{pinachw}}{step2}{}
Given any $r_0>0$ there exists  $\mu\geq0$ such that for every $\nu>0$
there is a length $A>0$ such that for all $a<A$ and $r>r_0$ 
$$d(w_{a,r}(v,z),\Pi_{\epsilon,a}(v,z))<(\mu+\nu)a\quad\text{ for all }v,z\in UX\text{
  with }p(v)=p(z).$$ 
\end{lemm}

\begin{lemmp}{step2}{}\\
Idea: Define $$\mu:=\frac1{r_0}+\underset{v,z}{\max}\frac{\partial}{\partial a}\Big|_{a=0}d(w_{a,r_0}(v,z),\Pi_{\epsilon,a}(v,z)).$$
The first summand times $a$ is less than the angle between the vectors in
$p(\Pi_{\epsilon,a}(v,z))$ which are radial to $\gamma_v(-r_0)$ and $\gamma_v(-r)$ for
$r>r_0$. The second obviously is an upper bound for the case $r=r_0$.
\end{lemmp} 

\begin{lemm}{\ref{wnachpsi}}{step3}{}
For every $\nu>0$ and $r_0>0$ there is a distance $A>0$ such that for all
$a<A$ the following holds:
$$\cos(\measuredangle(v,z))-\nu<\frac{l_a(r_0,v,z)-r_0}a<\cos(\measuredangle(v,z))+\nu$$
for all $v,z\in UX$ with $p(v)=p(z)$ and $\measuredangle(v,z)\leq\frac\pi2$.\\
Furthermore holds 
$$0\leq l_a(r,v,z)-r<l_a(r_0,v,z)-r_0\quad\text{ for }r>r_0$$
if $r_0>a$.
\end{lemm} 

\begin{lemmp}{step3}{}
For simplicity we suppose that the curvature is bounded by $-1\leq K\leq
0$. This rescaling does not change the result.\\
In the comparison spaces we can apply the Law of Cosine to get an estimate
of $l_a(r,v,z)$:
\begin{align}
l_a(r,v,z)^2&\geq r^2+a^2-2ra\cos(\pi-\measuredangle(v,z))\tag{Euclidean}\\
\notag \\
\cosh(l_a(r,v,z))&\leq\cosh(r)\cosh(a)-\cos(\pi-\measuredangle(v,z))\sinh(r)\sinh(a)\tag{Hyperbolic}
\end{align}
Fix $v,z$ and write $l_a(r)$ for one of the estimations 
\begin{align*}
l_a(r)&:=\sqrt{r^2+a^2+2ra\cos(\measuredangle(v,z))}\qquad\qquad\qquad\text{ or }\\
l_a(r)&:=\Arcosh\left(\cosh(r)\cosh(a)+\cos(\measuredangle(v,z))\sinh(r)\sinh(a)\right).
\end{align*} 
Notice
that $l_a(r)\geq r$ and 
$$l_a(r_0)-r_0\;\overset{a\rightarrow 0}\longrightarrow\;0.$$
So we can use de l'Hôpital to calculate
$$\underset{a\rightarrow0}\lim\frac{l_a(r_0)-r_0}{a}=\frac{\partial}{\partial
  a}\Big|_{a\rightarrow 0}(l_a(r_0)-r_0)=\cos(\measuredangle(v,z))$$
and we can find a value $A(v,z)$ depending on $v$ and
$\measuredangle(v,z)$. Since $UX$ is \Gam-compact and
$\measuredangle(v,z)\in[0,\frac\pi2]$ we can choose $A$ globally. To prove the
second part of the lemma calculate that
$$\frac{\partial}{\partial r}(l_a(r)-r)<0$$
for $r>a$.
\end{lemmp}

We need to distinguish one more case, namely when $v$ is close to $UR$ and
the gradient $\grad_{p(v)}d_{R'}$ points almost in the same direction as
$v$. In this case~\ref{wnachpsi} almost undoes the displacement
of~\ref{vnachpi}. By Corollary~\ref{gradn} we will see that this case can
be avoided. But first we need the following lemma.

\begin{lemm}{}{winkelklein}{}
Pick a sequence of vectors $\{v_n\}\subset UX$ with base  points
$p_n:=p(v_n)\in{\cal
  U}_\rho R\backslash {R'}$ (hence $w_n:=\grad_{p_n}d_{R'}$ is
defined). Suppose the base point sequence converges to a point $p\in
R$. Then there is a subsequence of $\{w_n\}$ converging to a normal vector
$w\in{\cal N}^1_pR$.
\end{lemm}
\begin{lemmp}{winkelklein}{}
The gradient defines a shortest geodesic $\gamma_{-w_n}$ from $p_n$ to some
$q_n\in R'$. Thus
$$q_n=\gamma_{-w_n}(d_{R'}(p_n))\in R'$$
and obviously $q_n\rightarrow p\in R$. So
$q_n\notin\partial R'$ holds for $d(p,p_n)$ small and hence ${\cal N}_{q_n}R'$ is well defined and 
$$z_n:=\dot
\gamma_{-w_n}(d_{R'}(p_n))\in{\cal N}_{q_n}^1R'.$$
Notice that $z_n$ is obtained from $w_n$ by parallel transport along
$\gamma_{-w_{n}}$ and hence $d(z_n,w_n)=d_R(p_n)\rightarrow 0$. Choose a
converging subsequence $w_n\rightarrow w\in UX$. Obviously $z_n$ converges
to the same vector $w$. But all the $z_n$ are normal vectors hence
$w\in{\cal N}_p^1R$. We conclude 
$$\grad_{p_n}d_{R'}=w_n\rightarrow w\in {\cal N}_p^1R$$
\end{lemmp}

\begin{prop}{Distance Gradient and Unit Tangent Bundle}{gradn}{}
Choose $\rho<\frac\pi2$ so small that the gradient is defined on ${\cal
  U}_\rho R'\backslash R'$. Then for every
given $\eta'<\frac\pi2$ we can find $\eta>0$ such that the following implication holds for all 
$v\in UX$ with base  point in
$\U_\rho R\backslash R'$:
$$\measuredangle(v,\grad_{p(v)}d_{R'})<\eta'\quad\Longrightarrow\quad
d(v,UR)>\eta.$$
I.\ e.\ for every $\eta'$-neighbourhood of the gradient vectors on
$\U_\rho R\backslash R'$ we can find an $\eta$-neighbourhood of $UR$ such
that these neighbourhoods do not intersect.
\end{prop}

\begin{propp}{gradn}{}
Given $\eta'<\frac\pi2$ choose
$$\eta<\inf\{d(v,UR)\;\big|\;v\in UX,\;p(v)\in {\cal U}_\rho R\backslash
R',\;\measuredangle(v,\grad_{p(v)}d_{R'})<\eta'\}.$$ 
Suppose the infimum on the right hand side is $0$. Then we can construct a
sequence contradicting Lemma~\ref{winkelklein}.

\end{propp}

Now we can define the desired deformation:
\begin{prop}{}{Psir}{}
Given $r_0>0$ we can find $C,\tau,\epsilon>0$ such that for all
$\lambda\in]0,1]$ the maps $\Psi_{r,\lambda\epsilon}$ have the following
properties for all $r\geq r_0$:
\begin{enumerate}
\item\label{Psistetig}
$\Psi_{r,\lambda\epsilon}:UX\backslash p\inv R'\rightarrow UX$ is continuous and respects the $r$-sphere foliation of $UX$:\\ 

$$\Psi_{r,\lambda\epsilon}(v)\in\sphere_v(r)\qquad\text{ for all }\qquad v\in UX\backslash p\inv R'.$$
\item\label{Psiabstand}
The image of $\Psi_{r,\lambda\epsilon}$ has no intersection with the $\lambda\tau$-neighbourhood of $UR$:\\

$$d_{UR}(\Psi_{r,\lambda\epsilon}(v))\geq\lambda\tau\qquad\text{ for all
  }\qquad v\in UX\backslash p\inv R'.$$
\item\label{Psiverschiebung}
The displacement by $\Psi_{r,\lambda\epsilon}$ is bounded by the global
constant $\lambda C$:\\
$$d(v,\Psi_{r,\lambda\epsilon}(v))\leq \lambda C\qquad\text{ for all
  }\qquad v\in UX\backslash p\inv R'.$$
\end{enumerate}
\end{prop}

\begin{propp}{Psir}{}
First we have to choose the constants that appear in the proposition.

Fix $\eta':=\arccos(\frac14)$ and $\nu:=\frac14$.
Take $\eta$ as provided by Proposition~\ref{gradn} and
$\mu$ and $A$ from Lemma~\ref{step2} and Lemma~\ref{step3}. 

Next choose $\tau$, $\kappa$, $\epsilon$ and $C$.
$$\tau:=4\frac{\min(2A,\eta,\frac\rho2)}{(3+\mu)(17+4\mu)}$$
\[
\kappa:=\frac{9+4\mu}2\tau\qquad
\epsilon:=\frac{17+4\mu}8\tau\qquad
C:=\frac{7+4\mu}2\epsilon.
\]
Notice that for the scaled constants $\tau':=\lambda\tau$,
$\kappa':=\lambda\kappa$, $\epsilon':=\lambda\epsilon$ and $C':=\lambda C$

$$\tau'\leq 4\frac{\min(2A,\eta,\frac\rho2)}{(3+\mu)(17+4\mu)}$$

\[
\kappa':=\frac{9+4\mu}2\tau'\qquad
\epsilon':=\frac{17+4\mu}8\tau'\qquad
C':=\frac{7+4\mu}2\epsilon'.
\]

It is then easy to see that all the following calculations hold if we
replace the constants $\tau$, $\kappa$, $\epsilon$ and $C$ by the scaled
constants $\tau'$, $\kappa'$, $\epsilon'$ and $C'$. So suppose $\lambda=1$.\\

With the chosen constants we can now prove the proposition:

\begin{enumerate}

\item This property is obvious from the definition of $\Psi_{r,\epsilon}$.

\item
First of all notice that $\kappa\leq 4\epsilon$ and $(6+2\mu)\epsilon\leq \eta$. Therefore, by
Proposition~\ref{gradn}, holds for all $v$ with $d_{UR}(v)\leq (6+2\mu)\epsilon$ that
$\cos(\measuredangle(v,\grad_pd_{R'}))\leq\cos(\eta')\leq\frac14$ and we can
apply Lemma~\ref{step3} to see that for $v$ with $d_{UR}(v)\leq
(6+2\mu)\epsilon$ 
\begin{align*}
d(w_r(v),\Psi_{r,\epsilon}(v))&=l_{\alpha(v)}(r,v,\grad_p(d_{R'}))-r\\
                   &<(\cos(\eta')+\nu)\alpha(v)\\
                   &=\frac{\alpha(v)}2.
\end{align*}
For other $v$ we get the weaker estimation 
\begin{align*}
d(w_r(v),\Psi_{r,\epsilon}(v))&=l_{\alpha(v)}(r,v,\grad_p(d_{R'}))-r\\
                   &<(1+\nu)\alpha(v)\\
                   &=\frac54\alpha(v).
\end{align*}

Furthermore, by Lemma~\ref{step2}, 
$$d(\Pi_\epsilon(v),w_r(v))\leq(\mu+\nu)\alpha(v)$$
for these $v$.\\
We will
distinguish four cases.
\begin{description}
\item[If \mathversion{bold}$\;d_{U{R'}}(v)\geq 4\epsilon$\mathversion{normal}] then $\Pi_\epsilon(v)=v$ and hence $\Psi_{r,\epsilon}(v)=v$
  and therefore 
$$d_{UR}(\Psi_{r,\epsilon}(v))\geq
d_{U{R'}}(\Psi_{r,\epsilon}(v))\geq4\epsilon.$$
\item[If \mathversion{bold}$\;\kappa\leq d_{U{R'}}(v)\leq 4\epsilon\;$ and $\;d_{UR}(v)\leq(6+2\mu)\epsilon$\mathversion{normal}] then
\begin{align*}
d_{UR}(\Psi_{r,\epsilon}(v))\qquad\\ \geq 
d_{U{R'}}(\Psi_{r,\epsilon}(v))&\geq 
d_{U{R'}}(v)
-\drunter{=\alpha(v)}{d(v,\Pi_\epsilon(v))}
-\drunter{\leq(\mu+\nu)\alpha(v)}{d(\Pi_\epsilon(v),w_r(v))}
-\drunter{\leq\alpha(v)/2}{d(w_r(v),\Psi_{r,\epsilon}(v))}\\
&\geq d_{U{R'}}(v)-(\frac74+\mu)\drunter{=2\epsilon-\frac12d_{U{R'}}(v)}{\alpha(v)}\\
&=(\frac{15}8+\frac\mu2)\drunter{\geq\kappa=\frac{9+4\mu}{2}\tau}{d_{U{R'}}(v)}
-(\frac72+2\mu)\drunter{=\frac{17+4\mu}{8}\tau}{\epsilon}\\
&\geq\tau
\end{align*}
\item[If  \mathversion{bold}$\;0<d_{U{R'}}(v)\leq\kappa\;$ and $\;d_{UR}(v)\leq
  (6+2\mu)\epsilon$\mathversion{normal}] then notice that $w_r(v)$ and $\Pi_{\epsilon,r}(v)$
  have the same base  point in $X$ and by Property~\ref{aussage ueber dR} in
  Proposition~\ref{distortion} 
$$d_{{R'}}(p(w_r(v)))=d_{R'}(p(\Pi_\epsilon(v)))>\alpha(v).$$
Now
\begin{align*}
d_{UR}(\Psi_{r,\epsilon}(v))\geq 
d_{U{R'}}(\Psi_{r,\epsilon}(v))&\geq
\drunter{\geq d_{{R'}}(p(w_r(v)))>\alpha(v)}{d_{U{R'}}(w_r(v))}
-\drunter{\leq\alpha(v)/2}{d(w_r(v),\Psi_{r,\epsilon}(v))}\\
&\geq\frac12\alpha(v)\\
&=\drunter{=\frac{17+4\mu}8\tau}{\epsilon}-\frac14\drunter{\leq\kappa=\frac{9+4\mu}2\tau}{d_{U{R'}}(v)}\\
&\geq \tau.
\end{align*}

\item[If
  \mathversion{bold}$\;d_{UR}(v)\geq(6+2\mu)\epsilon$\mathversion{normal}]
  then 
\begin{align*}
d_{U{R}}(\Psi_{r,\epsilon}(v))&\geq 
d_{U{R}}(v)
-\drunter{=\alpha(v)}{d(v,\Pi_\epsilon(v))}
-\drunter{\leq(\mu+\nu)\alpha(v)}{d(\Pi_\epsilon(v),w_r(v))}
-\drunter{\leq\frac54\alpha(v)\leq\frac52\epsilon}{d(w_r(v),\Psi_{r,\epsilon}(v))}\\
&\geq \drunter{\geq (6+2\mu)\epsilon}{d_{U{R}}(v)}-\drunter{=5+2\mu}{(\frac92+2\mu+2\nu)}\epsilon\\
&\geq\epsilon\\
&\geq\tau
\end{align*}  

\end{description}
\item For $d_{UR'}(v)>4\epsilon$ there is no displacement, so suppose $d_{UR'}(v)\leq4\epsilon<\eta$. Then
\begin{align*}
d(v,\Psi_{r,\epsilon}(v))&\leq
\drunter{\leq\alpha(v)}{d(v,\Pi_\epsilon(v))}
+\drunter{\leq(\mu+\nu)\alpha(v)}{d(\Pi_\epsilon(v),w(v))}
+\drunter{\leq\alpha(v)/2}{d(w(v),\Psi_{r,\epsilon}(v))}\\
&\leq(\frac74+\mu)\drunter{\leq 2\epsilon}{\alpha(v)}\\
&\leq\frac{7+4\mu}2\epsilon=C
\end{align*}

\end{enumerate}
\end{propp}

So by Proposition~\ref{Psir} we can deform
$\sphere_v(r)\cap\left(UX\backslash p\inv R'\right)$ inside $\sphere_v(r)$
away from $UR$. But what happens if we want to deform some structure, say
$Y$, inside $\sphere_v(r)$ away from $UR$ and the intersection $p(Y)\cap
R'$ is not empty? In this case we have to consider the dimensions of $Y$
and $R$. If the dimensions add up to less then the dimension of
$\sphere_v(r)$ (i.\ e.\ $\dim Y+\dim R<\dim X-1$), then the intersection is not transversal and can be undone
by an infinitesimal displacement of $Y$. If we assume that the intersection
of $R$ and $\sphere_v(r)$ is transversal, which we can by Sard's
Theorem if we vary the radius r slightly,
\index{Sard's Theorem}\index{Theorem!Sard's}
then the dimensional condition becomes
$\dim Y +\dim R <\dim X$.

\begin{coro}{}{minderdimensional}{}
Let $X$ denote a Hadamard manifold with compact quotient $M=X/\Gamma$.
Suppose a \Gam-compact submanifold $R\subset X$ and a radius $r_0>0$ are given.\\
Then there are constants $C,\tau>0$ such that for any $\lambda\in]0,1]$, $r>r_0$, $v\in UX$
and manifold $Y$ with 
$$\dim Y<\dim X -\dim R$$
we can deform any
smooth map $c:Y\rightarrow\sphere_r(v)$ into a smooth map $c_\lambda:Y\rightarrow \sphere_r(v)$ such that
\begin{itemize}
\item $c_\lambda$ is $(C+\frac\tau2)\lambda$-close to $c$, i.\ e.\ 
$$d(c(x),c_\lambda(x))\leq \left(C+\frac\tau2\right)\lambda\quad\text{ for all } x\in Y$$ and
\item $c_\lambda$ avoids a $\frac\tau2\lambda$-neighbourhood of $UR$, i.\ e.\ 
$$d_{UR}(c_\lambda(x))\geq\frac\tau2\lambda\quad\text{ for all }x\in Y.$$
\end{itemize}

\end{coro}

\begin{corop}{minderdimensional}{}\\
In the following we will apply Sard's Theorem, and the Preimage Theorem as
stated in~\cite{DiffTop}. Thom's Transversality Theorem is applied as
stated in~\cite{Demazure1989}.\\
Choose $C$, $\epsilon$ and $\tau$ as provided by Proposition~\ref{Psir} for
$r_0/2$. Without loss of generality we can assume that $\tau<r_0/2$ since
the constants are scalable. Fix these constants and consider any $r\geq r_0$
and
$\lambda\in]0,1]$. By Corollary~\ref{fusspunktnahe} we can find $\tau'>0$
such that on big spheres (radius $>r_0/2$) vectors are
$\frac\tau6\lambda$-close if their footpoints are $\tau'$-close.\\
Now write $o:=\gamma_v(-r_0)$ for the center of our sphere and consider the
differentiable map
$$
\begin{array}{cccc}
\rho:&R'\backslash{o}&\longrightarrow&\mathbbR\\
     &p             &\longmapsto    &d(o,p),
\end{array}
$$
which measures the radius in polar coordinates around $o$. By Sard's
Theorem the set of regular values of $\rho$ is dense in $\mathbbR$. Hence
we can choose a regular value $r'\in]r-\frac\tau6\lambda,r+\frac\tau6\lambda[$. The inverse image
$R'_{r'}:=\rho\inv(r')$ is the intersection of $R'$ with the sphere $S_v(r')$ in $X$.
Since $r'$ is regular, this intersection $R'_{r'}$ is a smooth submanifold of
dimension at most $\dim R -\dim \mathbbR=\dim R-1$ in the sphere $S_v(r')$
by the Preimage Theorem.
Next notice that the geodesic flow $\phi_{r'-r}$ identifies the two spheres
$\sphere_v(r)\rightarrow\sphere_{\phi_{r'-r}v}(r')$, moving every vector by exactly
$|r'-r|<\frac\tau6\lambda$. Furthermore we have natural identifications of spheres
with center $o$ in $X$ and in $UX$, given by the gradient of $\rho$:
$\grad \rho\big|_{S_v(r)}:S_v(r)\rightarrow\sphere_v(r)$ and
$\grad\rho\big|_{S_v(r')}:S_v(r')\rightarrow \sphere_v(r')$.\\

Now $y:=p\circ\phi_{r'-r}\circ c:Y\rightarrow S_v(r')$ is a
smooth map of $Y$ into the sphere $S_v(r')$. By Thom's Transversality
Theorem we can find a $\tau'$-close smooth map $y':Y\rightarrow
S_v(r')$ which avoids $R'_{r'}$. Write $\alpha$ for the homotopy of
$S_v(r')$ with $\alpha_0=id_{S_v(r')}$ and $\alpha_1\circ y=y'$. Consider the map
$$c_\lambda:=\phi_{r-r'}\circ\Psi_{r,\lambda\epsilon}\circ\grad\rho
\circ\alpha_1\circ p\circ \phi_{r'-r}\circ c:Y\rightarrow \sphere_v(r)$$
which looks quite monstrous at first glance but we will explain it step by
step:
\begin{itemize}
\item $\phi_{r'-r}\circ c$ is $\frac\tau6\lambda$-close to $c$, since $r'-r<\frac\tau6\lambda$.
\item The base point distance between $\phi_{r'-r}\circ c$ and
  $\grad\rho\circ\alpha_1\circ p\circ \phi_{r'-r}\circ c$  is just the
  distance between $y$ and $y'$ which is smaller than $\tau'$ by definition
  of $\alpha$. But $\tau'$-close base points imply
  $\frac\tau6\lambda$-close radial 
  vectors and hence $\grad\rho\circ\alpha_1\circ p\circ \phi_{r'-r}\circ
  c=\grad\rho\circ y'$ is $\frac\tau6\lambda$-close to  $\phi_{r'-r}\circ c$
  and $\frac\tau3\lambda$-close to $c$.
\item $y'$ avoids $R'$ and hence $\grad\rho\circ y'$ maps $Y$ into
  $UX\backslash p\inv R'$ and the map
  $\Psi_{r,\lambda\epsilon}\circ\grad\rho\circ y'$ is well defined. By
  Proposition~\ref{Psir} it is $C\lambda$-close to $\grad\rho\circ y'$ and
  hence $(C\lambda+\frac\tau3\lambda)$-close to our original map $c$. Again by
  Proposition~\ref{Psir} it avoids a $\tau\lambda$-neighbourhood of $UR$.
\item Projecting this map back to $\sphere_v(r)$ by the geodesic flow
  $\phi_{r-r'}$ gives a further displacement of $|r-r'|\leq\frac\tau6\lambda$
  which leaves us with the properties that $c_\lambda$ is
  $(C\lambda+\frac\tau3\lambda+\frac\tau6\lambda)$-close to $c$ and avoids a
  $(\tau\lambda-\frac\tau6\lambda)$-neighbourhood of $UR$.
\end{itemize}
This ends the proof since
$C\lambda+\frac\tau3\lambda+\frac\tau6\lambda=(C+\frac\tau2)\lambda$ and $\tau\lambda-\frac\tau6\lambda>\frac\tau2\lambda$.

\end{corop}

Corollary~\ref{minderdimensional} stays true (with different constants) if
$R$ is not a submanifold of $X$ but a subset stratified by closed
submanifolds of $X$.%
\footnote{E.\ g.\ if $R$ is a subanalytic subset
\index{subanalytic subset}\index{subset!subanalytic}}

\begin{coro}{}{subanalytischesausweichen}{}
Let $X$ denote a Hadamard manifold with compact quotient $M=X/\Gamma$.
Suppose a \Gam-compact stratified subset $R\subset X$ and $r_0>0$ are given.\\

Then there is a constant $k>0$ such that for small $\lambda>0$ and given $r>r_0$, $v\in UX$
and manifold $Y$ with
$$\dim Y<\dim X -\dim R$$
we can deform any
continuous map $c:Y\rightarrow\sphere_r(v)$ into a continuous map
$c_\lambda :Y\rightarrow \sphere_r(v)$ such that $c_\lambda$ is
$k\lambda$-close to $c$ and $c_\lambda$ avoids a $\lambda$-neighbourhood of
$UR$, i.\ e.\ for all $x\in Y$ 
$$d(c(x),c_\lambda(x))\leq k\lambda\qquad\text{ and }\qquad d_{UR}(c_\lambda(x))\geq\lambda.$$
\end{coro}

\begin{corop}{subanalytischesausweichen}{}
For a start suppose that $c$ is a smooth map of $Y$ into $\sphere_v(r)$. 
$R$ is a locally finite union $R=\bigcup R_i$ of submanifolds of $X$. By
\Gam-compactness the union can be assumed to be finite. Suppose the index set
is $i\in\{1\ldots j\}$. For each $i$ we can find constants
$\tau_i$ and $C_i$
such that Corollary~\ref{minderdimensional} holds. Start with $\tau_1$ and $C_1$. Choose
$\lambda_2$ small enough to assure that $\lambda_2 C_2<\tau_1/8$ and
$\lambda_2\tau_2\leq\tau_1/4$. So if we deform 
$c:Y\rightarrow \sphere_r(v)$ first away from $R_1$ and then away from
$R_2$ the resulting map, by Corollary~\ref{minderdimensional} still avoids a
$\tau_1/2-\lambda_2C_2-\lambda_2\tau_2/2\geq\tau_1/4$-neighbourhood of $UR_1$. Iterating this
concept we get constants $\tau=\frac{\tau_1}{2^j}$ and $k'=\frac{\sum C_i}\tau$ with
the desired property, where $\lambda$ is small, if $\lambda<\tau$.\\
Now return to the case where $c$ is not smooth but only continuous. In this
case for every $0<\lambda<\tau$ there is a smooth map $c':Y\rightarrow
\sphere_v(r)$ which is $\lambda$-close to $c$. We can apply what we
just learned to $c'$ to get $c_\lambda$, which is $k'\lambda$-close to $c'$
and $\lambda$-far from $UR$. Hence $c_\lambda$ satisfies
$$d(c(x),c_\lambda(x))\leq (k'+1)\lambda\qquad\text{ and }\qquad d_{UR}(c_\lambda(x))\geq\lambda$$
and hence the corollary holds for $k:=k'+1$ and $\lambda<\tau$. 

\end{corop}

\mysection{Avoiding Higher Rank Vectors}{}\label{rankr5.tex}
Suppose we are given a Hadamard manifold $X$ with compact quotient
$M=X/\Gamma$. Suppose 
furthermore that $R$ is an s-support of $\R_>$, the set of vectors of
higher rank. In this section we prove that for every vector $v\in UM$ there
is a geodesic ray $\gamma$ avoiding a neighbourhood of $\R_>$ such that
$\gamma$ starts in the base point of $v$ and $\dot\gamma(0)$ is close to
$v$. Furthermore our construction does not only work for a single vector
but for manifolds of directions. If we start with a
whole manifold $Y$ of directions in one point we
can deform this manifold to a close set of geodesic rays which avoids a
neighbourhood of $\R_>$. This works if the dimension of $Y$ is bounded by
$\dim Y<\dim X-\sdim R_>$:

\paragraph*{Proposition~\ref{summary}}$ $\\
\begin{it}
Let $X$ denote a rank one Hadamard manifold with compact quotient
$M=X/\Gamma$. Suppose $\sdim(\R_>)<\dim X-1$ then there are constants $\epsilon$ and $c$ such that for any $\eta<\epsilon$
there is an $\eta'<\eta$ for which the following holds: 
\begin{quote} 
For any compact manifold $Y$ with 
$$\dim Y<\dim X -\sdim \R_>$$ and any continuous map
$v_0:Y\rightarrow U_oX$ from $Y$ to the unit tangent sphere at a point
$o\in X$ we can find a continuous map $v_\infty:Y\rightarrow U_oX$ which is
$c\eta$-close to $v_0$ and satisfies $$d({\R_>},\phi_{\mathbbR_+}(v_{\infty}))\geq\eta'.$$
\end{quote}
\end{it}
 
So from now on suppose that $Y$ is a compact
manifold which satisfies $\dim Y<\dim X-\sdim R_>$.
Given a continuous map $v:Y\rightarrow U_{o}X$ we want to find a $C$-close
map $v_{\infty }:Y\rightarrow U_{o}X$ such that the geodesic ray
$\gamma_{v_{\infty }(x)}([0,\infty [)$ 
avoids an $\eta $-neighbourhood of all vectors of higher rank. Under a
condition on $\sdim \R_>$ this will
be possible for all $\eta $ small enough and $C/\eta $ will be a constant
independent of $\eta $. 
To find $v_\infty$ we will construct a sequence of continuous maps $v_{i}:Y\rightarrow U_{o}X$
such that $\gamma _{v_{i}(x)}([0,t_{i}[)$ avoids an $\eta $ -neighbourhood
of all vectors of higher rank, where $t_{i}\rightarrow \infty $ and the
sequence $v_{i}$ converges uniformly to $v_{\infty }$.\\
Our strategy for the construction of the converging sequence is twofold. On
the one hand we use the construction of Section~\ref{ausweichen.tex} to
avoid $UR$. This will work well for fixed times, but under the geodesic
flow we lose control. So we have to use the fact that $\R_>$ is
$\Gamma$-compact and of higher rank while we will work entirely on a
\Gam-compact set of rank one vectors. By \Gam-compactness there is a finite
distance between these two sets.

\mysubsection{Rough Idea for a Single Vector}{}
%\neuesubsection{Rough Idea}
To motivate the construction let us consider one vector and try to deform the corresponding ray to a ray that avoids a neighbourhood of $\R_>$. We
work with the constants $0<\beta\epsilon<\epsilon<\delta\epsilon$
which will be defined later on.\\

Starting with a given vector $v\in U_oX$ we construct a sequence of
geodesic rays $\gamma_{v_i}=:\gamma_i$ with $v_i\in U_oX$ recursively. At
the same time we construct a sequence of times $t_i\in \mathbbR_+$ with
$t_i>i\tilde B$ and hence $t_i\rightarrow\infty$. Inductively we want these
sequences to satisfy that $\dot\gamma_i$ is 
\begin{itemize}
\item $\beta\epsilon$-far away from $\R_>$ at all times $t\in[0,t_{i-1}]$,
\item $\epsilon$-far away from $UR$ at some control times $t_{i-1,k}\in[0,t_{i-1}]$ and even
\item $10\epsilon$-far away from $UR$ at the time $t_{i-1}$.
\end{itemize}
Now suppose this is true for some $i>0$. If $\dot\gamma_i$ is
$\beta\epsilon$-far away from $\R_>$ for all times $t>0$ we are done, if
not wait until time $t_{i-1}+\tilde B$ and let $t_i$ denote the first time
after $t_{i-1}+\tilde B$, where $\gamma_i$ is less than $11\epsilon$-far
away from $UR$. At this time move away from $UR$ to get a geodesic
$\gamma_{i+1}$ that is $10\epsilon$-far away from $UR$ at the time
$t_i$. This movement may be chosen such that $\gamma_{i+1}$ and $\gamma_i$
are within distance less than $\delta\epsilon$ from each other at time
$t_i$. Now, since $\gamma_i$ is of rank one, the displacement of
$\delta\epsilon$ at time $t_i$ has a very small effect at times
$t<t_i-\tilde B$. Namely we can assume that it is shrunk by a factor of
$1/N$. This assures that at the control times $t_{i,k}\approx t_{i-1,k}$
the geodesic $\gamma_{i+1}$ is close to $\gamma_i$ and therefore is
$\epsilon$-far away from $UR$. From this we can deduce that $\gamma_{i+1}$
is $\beta\epsilon$-far away from $\R_>$ for all $t\in[0,t_i]$.

\mysubsection{Construction of the Maps $v_{i}$}{Deformation}\label{howto}

Suppose we have $v_{i}:Y\rightarrow U_{o}X$ and we know that this map has
the desired properties at time $t_{i}$. Then we have a control on how close
the map will come to $\R_>$ until the time $t_{i}+\tilde{B}$. At this
time check, whether the map is too close to $UR$.
If it is then define $t_{i+1}:=t_{i}+\tilde{B}$. Else start
to check for all times $t>t_{i}+\tilde{B}$ whether
$v_{i}$ gets too close to $UR$ at that time. If it does, define
$t_{i+1}:=t$. If it does not then define $t_{i+1}:=t_{i}+\tilde{B}$ and
obviously the following deformations will not change $v_i$.\\ 

Now apply the distortion provided by Corollary~\ref{subanalytischesausweichen} to the
image of $v_{i}$ at time $t_{i+1}$ on the sphere $\sphere _{t_{i+1}}(o)$,
the set of outwards normal vectors to the sphere around $o$ of radius $t_{i+1}$. Using the notation of the corollary, the map $c=\phi _{t_{i+1}}\circ v_{i}$, gets distorted into the map $\bar{c}:Y\rightarrow \sphere_{t_{i+1}}(o)$. Applying
the geodesic flow again we pull this map back to $U_{o}X$ and call the
resulting map $v_{i+1}:=\phi _{-t_{i+1}}\circ \bar{c}:Y\rightarrow U_oX$ the \emph{deformation of
$v_i$ at time $t_i$}\index{deformation}. The corollary
tells us that the resulting map is $\tau /2$ -good at time $t_{i+1}$
and $d\left (\phi _{t_{i+1}}(v_{i}),\phi _{t_{i+1}}(v_{i+1})\right )\leq C+\frac{\tau }{2}$.
This will suffice to guarantee that $d(v_{i},v_{i+1})\leq \frac{1}{2^{i-1}}(C+\frac{\tau }{2})$
and hence the sequence $\{v_{i}\}_{i\in \mathbb N}$ converges uniformly.
We will need to fix some constants to assure that our construction has the
desired properties.

\mysubsection{Choice of Constants}{}
\begin{enumerate}
\item\label{deftau}\label{defC}Fix a starting radius $r_{0}>2$ and find
  constants $\tau $ and $C$ by Corollary~\ref{minderdimensional}. So on big
  spheres we can deform a submanifold of the sphere by less than $C+\tau/2$
  to get a submanifold which avoids a $\tau/2$-neighbourhood of $UR$. (Keep in mind that $\tau$ and $C$ may be scaled down).
\item\label{defepsilon} Fix $\epsilon:=\frac{\tau }{20}$, so $\epsilon$
  depends on $\tau$. 

\item\label{defdelta} Define $\delta:=\frac{C+\frac{\tau }{2}}{\epsilon
    }$. Notice that $\delta$ does not depend on the choice of $\tau$, since
  $C/\tau$ and $\epsilon/\tau$ are constants. $\delta\epsilon$ is a first
  estimate of the shield radius.

\item\label{defN}Fix $N\in\mathbb N$ with
  $$N>1+\frac\delta 2$$ keeping in mind that $N$ is independent of $\tau$.
  This will be the shrinking factor we need.
\item\label{deflambda} By the definition of $N$ we can choose $\lambda$
  with
  $$\frac{\delta}{N-1}\;<\;\lambda\;<\;2N-\delta.$$
  We need to enlarge the shield radius by $\lambda\epsilon$. 
\item\label{defDelta} Define $\Delta>(\delta+\lambda)\epsilon$.

\item\label{defA} Fix a distance $A>0$ as provided by Corollary~\ref{best
    hyper} for the \Gam-compact set $K_{\epsilon}:=UX\backslash\W_\epsilon UR$, the shield radius
  $\Delta$, the start radius $r_0$ and the shrinking factor $N$.\\
  So we need to consider distances bigger than $A$ before we can apply the
  shrinking property.
\item\label{defBtilde} Fix $\tilde B>A+8\epsilon$ and demand that
  $\tilde B>r_0+4\epsilon$. This is a technical choice. We have to work
  with a number bigger than $A$ to assure that even after small changes we
  can still apply the shrinking property.
\item\label{defbeta}% BETA WAR FRUEHER TAU!!!
Notice that for any $t$ the \Gam-compact set
  $\phi_t(K_{\epsilon})$ consists of rank one vectors only and define
  $0<\beta<1$ by 
  $$\beta\epsilon<\underset{t\in[-\tilde B,\tilde B]}\min\quad
  d(\phi_t(K_{\epsilon}),\R_>).$$ 
  During our construction we do not care what happens at times
  $t\in]t_i,t_i+\tilde B[$. But since we are only working on the \Gam-compact
  set $K_\epsilon$ of rank one
  vectors and try to avoid the \Gam-compact set of higher rank, there is a
  bound on how close we can get to $\R_>$ during a fixed time. Taking the
  min of $t\in[-\tilde B,\tilde B]$ instead of $t\in[0,\tilde B]$ is
  necessary to assure that after the first step we will only deal with rank
  one vectors.
\end{enumerate} 

\mysubsection{Construction of a Ray}{}

\begin{figure}[h]
\caption{Construction of $v_{i+1,i}$ and $t_{i+1,k}$}
\label{verschiebung}
\psfrag{o}{$o$}
\psfrag{gk1}{$\gamma_{v_{k+1}}$}
\psfrag{gi}{$\gamma_{v_i}$}
\psfrag{gi1}{$\gamma_{v_{i+1}}$}
\psfrag{vk1k1}{$v_{k+1,k+1}$}%{$\phi_{t_{k+1}}(v_{k+1})$} 
\psfrag{vk1k}{$v_{k+1,k}$}%{$\phi_{t_{k}}(v_{k+1})$} 
\psfrag{vik}{$v_{i,k}$}%{$\phi_{t_{i,k}}(v_{i})$} 
\psfrag{vi1k}{$v_{i+1,k}=\phi_{t_{i+1,k}}(v_{i+1})$}%{$\phi_{t_{i+1,k}}(v_{i+1})$}
\psfrag{vii}{$v_{i,i}$}%{$\phi_{t_i}(v_i)$}
\psfrag{w}{$w$}
\psfrag{vi1i}{$v_{i+1,i}$}%{$\phi_{t_i}(v_{i+1})$}
\psfrag{vi1i1}{$v_{i+1,i+1}$}%{$\phi_{t_{i+1}}(v_{i+1})$}

\psfrag{H}{$\horo_{v_{i,i}}^D$}
\psfrag{S}{$\sphere_{v_{i}}(t_{i})$}

\setlength{\meinbuffer}{\baselineskip} \multiply \meinbuffer by 20
  \includegraphics[width=\textwidth,height=\meinbuffer]{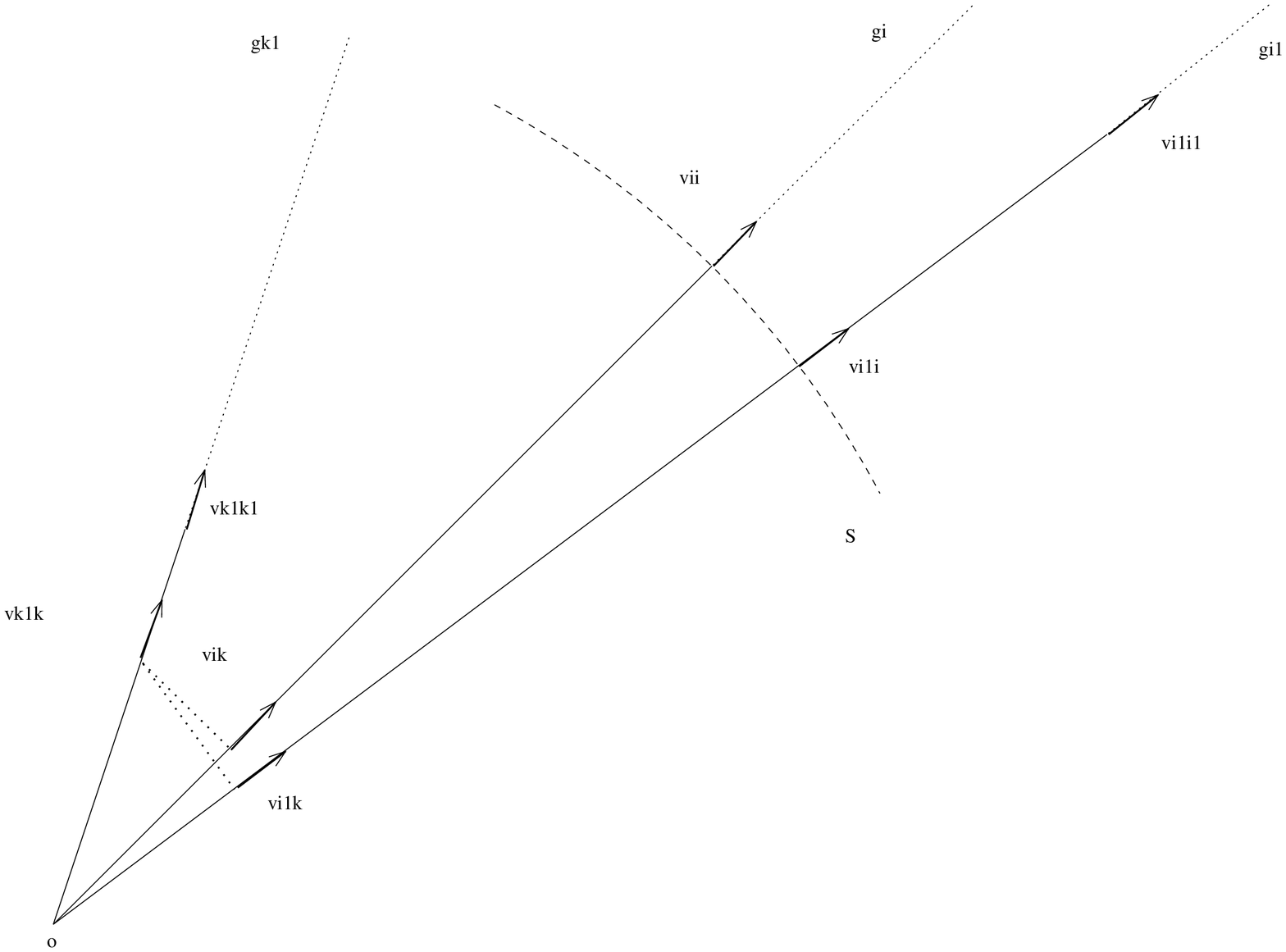}
\end{figure}

Starting with a continuous map $v_0:Y\rightarrow U_oX$ from a compact
manifold to the unit tangent sphere at $o$ and the time $t_{-1}:=0$ we construct sequences of times
$t_i$ and maps $v_i:Y\rightarrow U_oX$ as illustrated in Figure~\ref{verschiebung}.\\

Suppose $v_i$ and $t_{i-1}$ are given. Define
$$t_{i}:=\min\{t\geq t_{i-1}+\tilde B\;|\;d_{UR}(\phi_t(v_i(Y)))\leq11\epsilon\}$$
and choose $v_{i+1}$ to be the deformation of $v_i$ at time $t_i$ as
described in Section~\ref{howto}. Notice that $t_i$ might equal
$\infty$. In this case define $v_j:=v_i$ and $t_{j}:=t_{i-1}+(j-i+1)\tilde B$
for all $j\geq i$. Our sequences satisfy for all $x\in Y$ and $j\in \mathbb
N$
$$\begin{array}{cccccc}
10\epsilon&\leq&d_{UR}(\phi_{t_j}(v_{j+1}))&&&(\frac\tau2\text{-good})\\ \\
&&d(\phi_{t_j}(v_j(x)),\phi_{t_j}(v_{j+1}(x)))&\leq&\delta\epsilon&(\text{by
  definition of }\delta)\\ \\
t_j+\tilde B&\leq&t_{j+1}&&&(\text{by definition of }t_{j+1})
\end{array}$$
By the first inequality 
\begin{equation}\phi_{t_j}(v_{j+1})\in K_\epsilon.\label{viink}\end{equation}

We will show that the geodesics $\gamma_{v_i(x)}$ are good on the interval
$[0,t_i]$ and so is the limit. To do so we need to estimate the distance of
the vector $\phi_{t_j}(v_{j+1}(x))$ to the geodesics $\gamma_{v_i(x)}|_{[0,t_i]}$
which appear later in the sequence (i.\ e.\ $i>j$).\\
We define the auxiliatory maps $t_{i,j}:Y\rightarrow\mathbb R_+$, which need
not be continuous but are $4\epsilon$-close to $t_j$ as we will see in
Lemma~\ref{fuenfundsechs}.\\
For $0\leq k<i-1$ define $t_{i,k}:=t_{i,k}(x)$
\index{$t i,k$@$t_{i,k}$}
\index{t i,k@ $t_{i,k}$}
by
$$d(\phi_{t_{i,k}}(v_i),\phi_{t_k}(v_{k+1}))=d(\dot\gamma_{v_i},\phi_{t_k}(v_{k+1})),$$
i.\ e.\ $\phi_{t_{i,k}(x)}(v_i(x))$ is a vector tangent
to the geodesic ray
$\gamma_{v_i(x)}$ which is closest to the vector $\phi_{t_k}(v_{k+1}(x))$.%
\footnote{Notice that $t_{i,k}(x)$ might not be unique}
To shorten notation we will usually write $v_{i,k}:=\phi_{t_{i,k}}(v_i)$%
\index{$v i,k$@$v_{i,k}$}\index{v i,k@$v_{i,k}$}
and omit $x$ where possible.\\
Sometimes we might write $t_{i,i}$ or $t_{i+1,i}$ for $t_i$.\\
So we know 
\begin{equation}\label{vik}
d(v_{i+1,k},v_{k+1,k})=d(\dot\gamma_{v_{i+1}},v_{k+1,k})\quad\text{for}\quad
0\leq k<i
\end{equation}
and (\ref{viink}) translates into
\begin{equation}\label{viinK}
v_{i+1,i}\in K_\epsilon.
\end{equation}

\mysubsection{The Construction Has the Desired Properties}{Properties}
So far we know

\begin{alignat}{3}
\beta\epsilon\;<&\;d(v_{i,i},\R_>)\label{nummer1}\\
&\;d(v_{i+1,i}(x),v_{i,i}(x))&&\;\leq\;\delta \epsilon&&\text{ for all }x\in Y\label{nummer2}\\
10\epsilon\;\leq&\;d(v_{i+1,i},\R_>)\label{nummer3}
\end{alignat}

By~(\ref{viinK}) and the definition of $\beta$ we know
\begin{equation}
\beta\epsilon\leq d_{\R_>}(\phi_t(v_i))\quad\text{ for }\quad t_{i-1}-\tilde
B\leq t\leq t_{i-1}+\tilde B \label{nummer4}
\end{equation}
and by the choice of $t_{i}$ we know that if $t_{i}-t_{i-1}>\tilde B$ then 
\begin{equation}\label{nummer5}
11\epsilon<d_{UR}(\phi_t(v_i))\qquad \text{ for }\qquad t_{i-1}+\tilde
B\leq t<t_{i}.
\end{equation}

\begin{lemm}{}{einsundzwei}{}
For $0\leq k<i$ and all $x\in Y$ 
\begin{eqnarray}
d(v_{i,k}(x),v_{k,k}(x))&<&(\delta+\lambda)\epsilon,\label{z1}\\
d(v_{i,k}(x),v_{k+1,k}(x))&<&2\epsilon.\label{z2}
\end{eqnarray}
This is illustrated in Figure~\ref{verschiebung}.
\end{lemm}

\begin{lemmp}{einsundzwei}{}
Fix $x\in Y$.
By~(\ref{nummer2}) holds~(\ref{z1}) for $i=k+1$, since $\lambda>0$. We will
prove it for all $k$ by induction on $j:=i-k$. Suppose~(\ref{z1}) holds for
$i-k=j$. Using the fact that $(i+1)-(k+1)=j$, we get
\begin{align*}
d(\gamma_{v_{k+1}}(t_{k+1}),\gamma_{v_{i+1}})&\leq
d(\dot\gamma_{v_{k+1}}(t_{k+1}),\dot\gamma_{v_{i+1}})\\
&=d(\dot\gamma_{v_{k+1}}(t_{k+1}),\dot\gamma_{v_{i+1}}(t_{i+1,k+1}))\\
&=d(v_{k+1,k+1},v_{i+1,k+1})\\
&\overset{\text{per ind.}%\ref{z1}
}<(\delta+\lambda)\epsilon<\Delta.
\end{align*}
And since $t_{k+1}-t_{k}\geq\tilde B>A$ we can apply Corollary~\ref{best
  hyper} which yields (notice that $v_{k+1,k}\in K_\epsilon$)
$$d(v_{k+1,k},v_{i+1,k})=d(\dot\gamma_{v_{k+1}}(t_{k+1,k}),\dot\gamma_{v_{i+1}})<\frac{(\delta+\lambda)}{N}\epsilon\leq
2\epsilon$$
which proves~(\ref{z2}) for $(i+1)-k=j+1$ and furthermore
\begin{align*}
d(v_{i+1,k},v_{k,k})&\leq d(v_{i+1,k},v_{k+1,k})+d(v_{k+1,k},v_{k,k})\\
&\leq\frac{\delta+\lambda}{N}\epsilon+\delta \epsilon\\
&=(\delta+\frac{\delta+\lambda}{N})\epsilon<(\delta+\lambda)\epsilon
\end{align*}
thus we have proved~(\ref{z1}) for $(i+1)-k=j+1$.
\end{lemmp}

\begin{lemm}{}{dreiundvier}{}
For $0\leq k<i$ 
\begin{equation}
2\epsilon<d_{UR}(v_{i,k})\label{z3}
\end{equation}
hence we know
\begin{equation}
v_{i,k}(x)\in K_{\epsilon} \qquad\text{ for all }x\in Y\text{ and }0\leq k< i.\label{z4}
\end{equation}
\end{lemm} 

\begin{lemmp}{dreiundvier}{}
Suppose $2\epsilon\geq d_{UR}(v_{i,k})$ for some $k<i$. Then
\begin{align*}
10\epsilon&\overset{(\ref{nummer3})}{\leq\;} d_{UR}(v_{k+1,k})\\
&\leq\; d(v_{k+1,k},v_{i,k})+d_{UR}(v_{i,k})\\
&\leq\; 2\epsilon +2\epsilon=4\epsilon
\end{align*}
This is a contradiction.
\end{lemmp}

\begin{lemm}{}{fuenfundsechs}{}
For $0\leq k\leq i$ we can estimate $t_{i,k}$ by $t_{k,k}$ as follows:
\begin{equation}\label{z5}
|t_{k}-t_{i,k}|\begin{cases} =0\qquad\qquad\;\;\text{ if } i=k \text{ or } i=k+1\\
                             \leq 4\epsilon \qquad\qquad\text{ if } i>k+1.
               \end{cases}
\end{equation}
\end{lemm}

\begin{lemmp}{fuenfundsechs}{}
We can find an estimate for the angle $\measuredangle(v_i,v_{k+1})$ as
illustrated in Figure~\ref{angleestimate} and thus

\begin{align*}
t_{i,k}&=d(v_i,\phi_{t_{i,k}}(v_i))=d(v_i,v_{i,k})\\
&\geq \drunter{=t_{k,k}}{d(v_{k+1},v_{k+1,k})}
-
\drunter{
\begin{array}{l}
=\measuredangle(v_{k+1},v_i)\\
\leq\arcsin\left(\frac{2\epsilon}{t_k}\right)\\
{<}\frac{4\epsilon}{r_0}<2\epsilon
\end{array}
}{d(v_{k+1},v_i)}-\drunter{<2\epsilon}{d(v_{k+1,k},v_{i,k})}\\
&\geq t_{k}-4\epsilon
\end{align*}

and analogously
\begin{align*}
t_{k}=t_{k+1,k}&=d(v_{k+1},v_{k+1,k}))\\
&\geq d(v_i,v_{i,k})-d(v_{k+1},v_{i})-d(v_{i,k},v_{k+1,k})\\
&\geq t_{i,k}-2\epsilon-2\epsilon\\
&=t_{i,k}-4\epsilon.
\end{align*}
\begin{figure}[h]
\caption{Estimation of $\measuredangle(v_i,v_{k+1})$}
\label{angleestimate}
  \psfrag{o}{$o$} \psfrag{gk1}{$\gamma_{v_{k+1}}$} \psfrag{gi}{$\gamma_{v_i}$}
  \psfrag{gk10}{$v_{k+1}$} \psfrag{vk1k}{$v_{k+1,k}$} \psfrag{vik}{$v_{i,k}$}
  \psfrag{tik}{$t_{i,k}$} \psfrag{tk1k}{$t_{k+1,k}=t_k$} \psfrag{2e}{$<2\epsilon$}
  \psfrag{alpha}{$\alpha<\arcsin\left(\frac{2\epsilon}{t_{k}}\right)$}
  \psfrag{gi0}{$v_i$}
\setlength{\meinbuffer}{\baselineskip} \multiply \meinbuffer by 14
  \includegraphics
[height=\meinbuffer]
{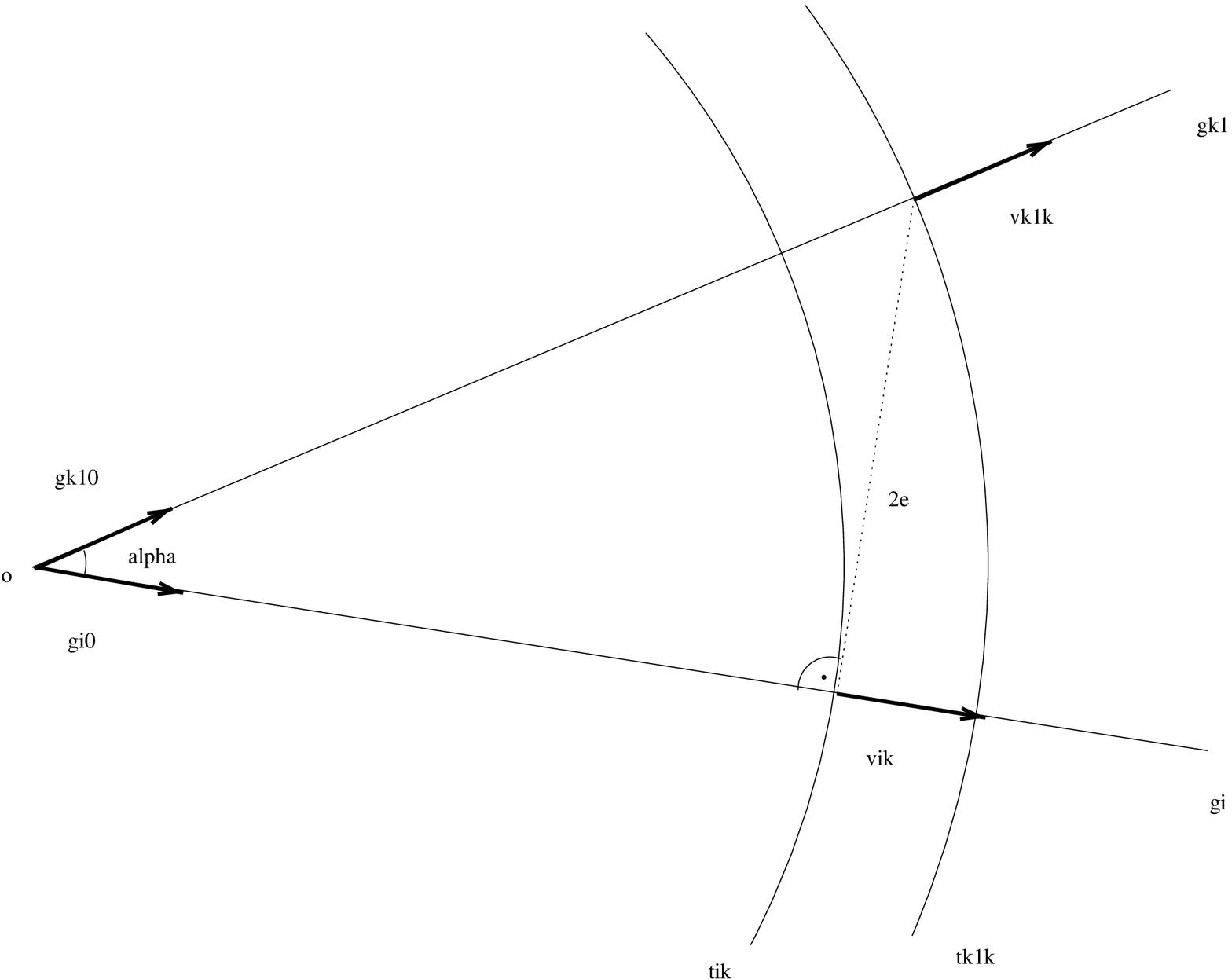}
\end{figure}
\end{lemmp}

\begin{lemm}{}{finale}{}
For all $i\geq 1$ and all $0\leq t\leq t_{i}$ 
\begin{equation*}
d_{\R_>}(\phi_t(v_i))>\beta\epsilon.
\end{equation*}
\end{lemm}

\begin{lemmp}{finale}{}
Suppose $d_{\R_>}(\phi_t(v_i))<\beta\epsilon$ for some $0<t<t_{i,i}$.\\
First of all $t>t_{i,1}$, by the definition of $\beta$ since $v_{i,1}\in
K_\epsilon$ by~(\ref{z4}). Hence there
is a $0\leq k<i$ with
$$t_{i,k}\leq t<t_{i,k+1}$$
and by~(\ref{nummer5}) we even know that $k+1<i$ since $\beta<11$. Since
all $v_{i,k}$ are in $K_{\epsilon}$ by~(\ref{z4}), the definition of
$\beta$ yields
\begin{align*}
t_{i,k}+\tilde B<t&<\quad t_{i,k+1}- \tilde B\\
       &\overset{(\ref{fuenfundsechs})}{\leq\quad} t_{k+1,k+1}+4\epsilon-\tilde B\\
       &<\quad t_{k+1,k+1}=t_{k+1,k}.
\end{align*}
Now
\begin{align*}
d(\phi_t(v_i),\phi_t(v_{k+1}))&\quad\leq\quad
d(\phi_{t_{k+1,k}}(v_i),\phi_{t_{k+1,k}}(v_{k+1}))\\
&\quad\leq\quad
\drunter{=|t_{k+1,k}-t_{i,k}|}{d(\phi_{t_{k+1,k}}(v_i),\phi_{t_{i,k}}(v_i))}+\drunter{=d(v_{i,k},v_{k+1,k})}{d(\phi_{t_{i,k}}(v_i),\phi_{t_{k+1,k}}(v_{k+1}))}\\
&\overset{\text{L. }\ref{fuenfundsechs}}{\quad\leq\quad}4\epsilon+d(v_{i,k},v_{k+1,k})\\
&\overset{(\ref{z2})}{\quad<\quad}6\epsilon
\end{align*}
and hence for $t_+:=t+4\epsilon$
\begin{align*}
d_{\R_>}(\phi_{t_+}(v_{k+1}))
&\leq\drunter{<\beta\epsilon}{d_{\R_>}(\phi_t(v_i))}
+\drunter{<6\epsilon}{d(\phi_t(v_i),\phi_t(v_{k+1}))}
+\drunter{=4\epsilon}{d(\phi_t(v_{k+1}),\phi_{t_+}(v_{k+1}))}\\
&<(10+\beta)\epsilon<11\epsilon.
\end{align*}

This is a contradiction to~(\ref{nummer5}) since
\begin{align*}
t_{k+1,k}+\tilde B\leq t_{i,k}+4\epsilon+\tilde
B&<t_+\\
&<t_{i,k+1}+4\epsilon-\tilde B\\
&\leq
t_{k+1,k+1}+4\epsilon+4\epsilon-\tilde B<t_{k+1,k+1}
\end{align*}
and hence
$$t_{k+1,k}+\tilde B<t_+<t_{k+1,k+1}.$$
\end{lemmp}

\begin{lemm}{}{convergence1}{}
For any $0\leq k<i$ 
$$d(v_{i+1,i-k}(x),\dot\gamma_{v_i(x)})\leq\frac\delta{N^k}\epsilon\text{
  for all }x\in Y.$$
\end{lemm}
\begin{lemmp}{convergence1}{}
We prove this by induction over $k$.\\
For $k=0$ the inequality becomes
$d(v_{i+1,i},\dot\gamma_{v_i})\leq\delta\epsilon$ which holds
by~(\ref{nummer2}).\\
Now suppose for some $0\leq k<i-1$ that
$$d(v_{i+1,i-k},\dot\gamma_{v_i})\leq\frac{\delta}{N^k}\epsilon.$$
We want to apply Corollary~\ref{best hyper}. To this end we have to
estimate $d(\gamma_{v_{i+1,i-k-1}}(A),\gamma_{v_i})$ to see that
$d(v_{i+1,i-k-1},\dot\gamma_{v_i})$ is less than $\frac1N$-times this
distance.\\
Recall that 
$$t_{i+1,i-k}-t_{i+1,i-k-1}\overset{\text{L. }\ref{fuenfundsechs}}{\quad\geq\quad}\tilde
B-8\epsilon\overset{(\ref{defBtilde})}{\quad>\quad}A$$
to get 
\begin{align*}
d(\gamma_{v_{i+1,i-k-1}}(A),\gamma_{v_i})
&\leq
d(\drunter{=p(v_{i+1,i-k})}{\gamma_{v_{i+1,i-k-1}}(t_{i+1,i-k}-t_{i+1,i-k-1})},\gamma_{v_i})\\
&\leq d(v_{i+1,i-k},\dot\gamma_{v_i})\\
&\leq\frac\delta{N^k}\epsilon<\Delta.
\end{align*}
We know by Lemma~\ref{dreiundvier} that $v_{i+1,i-k-1}\in K_\epsilon$ and
hence Corollary~\ref{best hyper} can be applied:
$$d(v_{i+1,i-k-1},\dot\gamma_{v_i})<\frac1N\frac\delta{N^k}\epsilon=\frac\delta{N^{k+1}}\epsilon.$$
\end{lemmp}

\begin{prop}{}{convergence2}{}
The sequence $(v_i)_i$ of maps $Y\rightarrow U_oX$ converges uniformly to a
map $v_\infty$ which is $\frac{\delta}{N(N-1)}\epsilon$-close to the
original map $v_0$ and with the property
\begin{equation}\label{abstand}
d_{\R_>}(\phi_t(v_\infty(x)))\geq\beta\epsilon
\end{equation}
for all $t\geq 0$ and $x\in Y$.
\end{prop}
\begin{propp}{convergence2}{}
By Lemma~\ref{convergence1} we know that (set $k:=i-1$)
$$d(v_{i+1,i},\dot\gamma_{v_i})\leq\frac\delta{N^{i-1}}\epsilon$$
We can therefore estimate $d(v_{i+1},v_i)$ as illustrated in Figure~\ref{konvergenz}. 
\begin{figure}[h]\caption{Estimation of $d(v_{i+1},v_i)$}\label{konvergenz}
\psfrag{o}{$o$}
\psfrag{angle}{$\measuredangle$}
\psfrag{vektori}{$v_i$}
\psfrag{vektori1}{$v_{i+1}$}
\psfrag{vektor}{$v_{i+1,i}$}
\psfrag{gammai}{$\gamma_i$}
\psfrag{gammai1}{$\gamma_{i+1}$}
\psfrag{distance}{$<\frac\delta{N^{i-1}}\epsilon$}
\psfrag{laenge}{$>1$}
  \includegraphics
[width=\textwidth]
{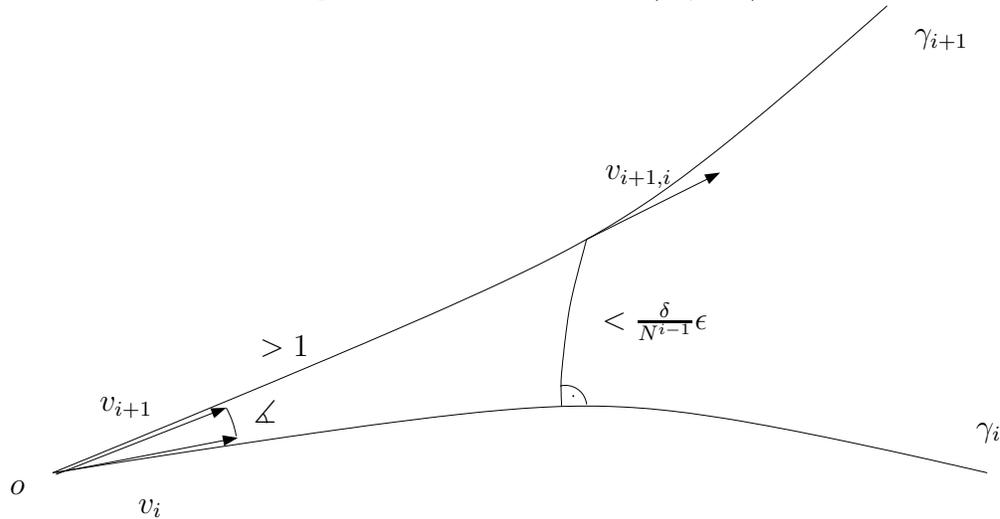}
\end{figure}
We need the estimations
\begin{xalignat*}{2}
d(p(v_{i+1,i},\gamma_{v_i})&\leq d(v_{i+1,i},\dot\gamma_{v_i})\qquad&t_{i+1,1}&\geq t_1-4\epsilon\\
&\leq\frac\delta{N^{i-1}}\epsilon&
&\geq\tilde B-4\epsilon\geq r_0>1
\end{xalignat*}
to conclude by comparison to a flat triangle
\begin{align*}
d(v_i,v_{i+1})&=\measuredangle(v_i,v_{i+1})\\
&\leq\arctan\left(\frac{\frac\delta{N^{i-1}}\epsilon}{1}\right)\\
&\leq\frac\delta{N^{i-1}}\epsilon.
\end{align*}
The sequence of distances $d(v_i,v_{i+1})$ is summable since $N>1$ and hence the $v_i$
converge to some $v_\infty$ which is $\sum
d(v_i,v_{i+1})=\frac{N\delta}{N-1}\epsilon$-close to $v_0$. Since the
estimation of $d(v_i,v_{i+1})$ is independent of $x$ the convergence is
uniform and $v_\infty$ is continuous. Property~(\ref{abstand}) is evident
from Lemma~\ref{finale}.
\end{propp}

Summing up the results of Section~\ref{rankr5.tex} we get the following
Proposition on the existence of rays avoiding $\R_>$.

\begin{prop}{Summary}{summary}{}
Let $X$ denote a rank one Hadamard manifold with compact quotient
$M=X/\Gamma$. Suppose $\sdim(\R_>)<\dim X-1$ then there are constants $\epsilon$ and $c$ such that for any $\eta<\epsilon$
there is an $\eta'<\eta$ for which the following holds: 
\begin{quote} 
For any compact manifold $Y$ with 
$$\dim Y<\dim X -\sdim \R_>$$ and any continuous map
$v_0:Y\rightarrow U_oX$ from $Y$ to the unit tangent sphere at a point
$o\in X$ we can find a continuous map $v_\infty:Y\rightarrow U_oX$ which is
$c\eta$-close to $v_0$ and satisfies $$d({\R_>},\phi_{\mathbbR_+}(v_{\infty}))\geq\eta'.$$
\end{quote}
\end{prop}

\begin{propp}{summary}{}
Take all the constants we had before. Note that by definition
$\delta$ and $N$ are independent of
$\epsilon$. So we can define a global constant $c:=\frac{\delta}{N(N-1)}$
(compare Proposition~\ref{convergence2}).
However, $\beta$, $\tilde B$, $A$ and $\Delta$ depend on $\epsilon$. So
given $\eta<\epsilon$ we have to go through the definitions of these
variables again, where we replace $\epsilon$ by $\eta$. Denote the
resulting constants by a hat and write $\eta':=\hat\beta\eta$ to get the
desired result. 
\end{propp}

If we start with just one vector and want to find a single ray avoiding
$\R_>$ we can apply Proposition~\ref{summary} for the case $\dim Y=0$.

\begin{coro}{}{Strahl}{}
Let $M$ denote a compact manifold of nonpositive curvature with
$$\sdim\R_><\dim M.$$ %-1 war's 
Then there are constants $\epsilon$ and $c$ such that for any $\eta<\epsilon$
there is an $\eta'<\eta$ for which the following holds: \\
For every vector $v\in UM$ there is a geodesic ray $\gamma$ in $M$
starting in the same base point
with initial direction $c\eta$-close to $v$ and
avoiding an $\eta'$-neighbourhood of
$\R_>$. Or short
$$\gamma(0)=p(v)\qquad\quad d(v,\dot\gamma(0))<c\eta\qquad\quad d(\R_>,\dot\gamma(\mathbbR_+))>\eta'.$$
\end{coro}

Notice that Corollary~\ref{Strahl} works for surfaces if and only if there
are only finitely many geodesics of higher rank and these are all
closed. In this case these geodesics are the totally geodesic s-support of $\R_>$.

\mysection{From Rays to Geodesics}{}\label{topologie.tex}
We have seen that any manifold of directions in one point can be deformed
into a manifold which consists of directions of geodesic rays that do not
come close to any vector of higher rank. In this section we will use a
topological argument to find a whole geodesic through that point. In fact,
for every small $\tau$ we will find a small $\epsilon$ such that in every
$\epsilon$-neigh\-bour\-hood of a given vector $v$ there is geodesic
through the 
base point of $v$ that avoids a $\tau$-neigh\-bour\-hood of all vectors of
higher rank. Since we are using a dimensional argument, the proof works
only for manifolds where the set $\R_>$ of all vectors of higher rank has
small s-dimension.

\begin{theo}{}{maintheorem}{}
Let $X$ denote a rank one Hadamard manifold with compact quotient $M$,
and suppose that
$$\sdim(\R_>)<\frac{\dim X}2.$$
Then there are constants $\epsilon$ and $c$ such that for any
$\eta<\epsilon$ there is an $\eta'<\eta$ with the following property:\\
For every vector $v_0\in
UX$ there is a $c\eta$-close vector $v$ with the same base point such that
$\gamma_v$ avoids an $\eta'$-neighbourhood of all vectors of higher rank.
\end{theo}

\begin{theop}{maintheorem}{}
The main ingredients for this proof are Proposition~\ref{summary} and some
intersection theory (e.\ g.~\cite{DiffTop}).\\
Start with the constants from Proposition~\ref{summary}. Given any vector $v_0\in U_pX$ we can find two great spheres
$Y_1$, $Y_2$ 
in $U_pX$ of dimensions%
\footnote{Recall the definition of the Gauss bracket: For
$x\in\mathbbR$ the largest integer less or equal to $x$ is denoted by $[x]$.
\index{x@$[x]$}\index{$x$@$[x]$}\index{$[.]$@$[.]$}\index{Gauss bracket}}
$\dim Y_1=\left[\frac{\dim X-1}2\right]$ and 
$\dim Y_2 =\left[\frac{\dim X}2\right]$ 
which intersect transversally in $v_0$ and $-v_0$. By the
condition on the s-dimension of $\R_>$ we have $\dim(Y_i)<\dim X -\dim R$
for $i=1,2$, where $R$ is an s-support of $\R_>$. Applying the proposition we can homotope $Y_1$ and
$Y_2$ into
$c\eta$-close submanifolds $\hat Y_i$ which consist only of initial vectors of
geodesic rays avoiding an $\eta'$-neighbourhood of $\R_>$. The
intersection number of simplices is invariant under homotopies. So the
intersection number of $\hat Y_1$ and $-\hat Y_2$ is the same as that of
$Y_1$ and $-Y_2$. But since
the points of intersection are (for $\eta$ small enough) much farther apart
then the displacement by 
the homotopy the two points of intersection $v_0$ and $-v_0$ can not
cancel out. Hence there are at least two points of intersection of $\hat
Y_1$ and $-\hat Y_2$. Pick the one $c\eta$-close to $v_0$ and call it $\hat v$.
Since $\hat v\in Y_1$ we know that $\dot\gamma_{\hat v}(\mathbbR_+)$ avoids an
$\eta'$-neighbourhood of $\R_>$. Since $-\hat v\in Y_2$ we know that the same
applies for $\dot\gamma_{-\hat v}(\mathbbR_+)=-\dot\gamma_{\hat v}(\mathbbR_-)$. Now notice that the set
$\R_>$ is symmetric - if $w\in U_qX$ is of higher rank, then so is the
vector $-w\in U_qX$. Therefore $\dot\gamma_{\hat v}(\mathbbR_-)$ avoids
an $\eta'$-neighbourhood of $\R_>$, too. 
\end{theop}

An easy consequence of~~Theorem~\ref{maintheorem} is

\begin{theo}{}{THE Theorem}{}
Let $M$ be a compact manifold of nonpositive curvature. Suppose the
s-dimension of the set of vectors of higher rank $\R_>$ is bounded by 
$$\sdim(\R_>)<\frac{\dim M}2.$$
  Then for every
  $\epsilon>0$ there is a closed, flow invariant, full, $\epsilon$-dense
  subset $Z_\epsilon$ of the unit tangent bundle $UM$ consisting only of
  vectors of rank one.
\end{theo}

\begin{rema}{}{total3}{}
\begin{enumerate}
\item $\R_>$ is flow invariant, since the
  rank is a property of the geodesic. Therefore whenever $p(\R_>)$ is a
  stratified subset of $M$ then $$\sdim R_>=\dim p(\R_>).$$
\item %\label{dimension ist fuenf}
  If $\R_>$ is not empty then
  the s-dimension of $\R_>$
  is at least one, this being the case when all the higher rank vectors are
  tangent 
  to a finite number of periodic geodesics of higher rank. 
  So the first dimension of a manifold where our theorems yield a
  result is in dimension three iff $p(\R_>)$ is a finite union of closed
  geodesics. 
\item The set constructed in Theorem~\ref{THE Theorem} is not only
  $\epsilon$-dense, but $\epsilon$-dense on every fibre $U_oM$.

\end{enumerate}
\end{rema}

\begin{appendix}

\markboth{Appendix}{}

{\LARGE \bf Appendix}
\addcontentsline{toc}{section}{Appendix}

\cleardoublepage
\setcounter{section}{1}

\section*{\thesection{}~Basics}
\addcontentsline{toc}{subsection}{\protect\numberline{\thesection}Basics}
\thispagestyle{empty}
\markboth{Appendix}{\thesection{} Basics}
\label{manifold.tex}
The aim of this appendix is to give an overview of the theory of manifolds
and Riemannian manifolds without going into detail. Most of the notation
and the theory introduced here is taken from the book by Sakai
\cite{Sakai1996}. Please
take a look at this very nice and concise book if you are interested in the
theory in more detail.

\subsection*{Manifold}
Suppose $M$ is a second countable, Hausdorff, topological space.\\
A pair $(U,\phi)$ is called an $n$-dimensional \emph{chart}\index{chart} of $M$ if
$\phi:U\rightarrow\mathbbR^n$ is a homeomorphism from an open subset
$U\subset M$ onto the image $\phi(U)$. In
this case, $U$ is called a \emph{coordinate neighbourhood}\index{coordinate
  neighbourhood}\index{neighbourhood!coordinate}. If $v_i$ denotes the
$i$-th coordinate of $v\in \mathbbR^n$ then the maps $x_i:u\mapsto\phi(u)_i$ are
called \emph{local coordinates}\index{local coordinates}\index{coordinates!local} for the coordinate neighbourhood $U$. \\
A collection $\cal A$
\index{A@$\cal A$}\index{$a$@$\cal A$}
of charts is called an \emph{atlas}\index{atlas} of $M$ if the coordinate
neighbourhoods of the charts cover all of $M$.\\

\begin{defi}{}{}{}
An $n$-dimensional, \emph{topological manifold}
\index{manifold!topological}\index{topological manifold}
$M$ is a second countable,
Hausdorff, topological space equipped with an atlas. 
\end{defi}

Suppose we have two charts $C_1:=(U_1,\phi_1)$ and $C_2:=(U_2,\phi_2)$ in
the atlas $\cal A$. If $U_1\cap U_2\neq\emptyset$ then we can define the
\emph{coordinate transformation}\index{coordinate transformation}\index{transformation!coordinate} from $C_1$ to $C_2$ to be
$\phi_2\circ\phi_1\inv:\phi_1(U_1\cap U_2)\rightarrow\phi_2(U_1\cap U_2)$. 
A topological manifold is called a \emph{smooth} (resp.\ \emph{real analytic})
\emph{manifold}
\index{manifold!smooth}\index{manifold!real analytic}
\index{smooth manifold}\index{real analytic manifold}\index{analytic manifold}
if all
coordinate transformations of its atlas are smooth (resp.\ real analytic).
For every atlas we can define the \emph{maximal atlas}\index{maximal atlas}\index{atlas!maximal}, that is the set of all
charts for which the coordinate transformations to and from the charts of
the atlas are smooth (resp.\ real analytic).
 
\begin{rema}{}{}{}
A special case are open subsets $M$ of $\mathbbR^n$ or of finite
dimensional vector spaces. The
inclusion into $\mathbbR^n$ is a chart which is an atlas since it covers
all of $M$. 
\end{rema}

\subsection*{Universal Covering}
For every manifold $M$ there is a
covering\index{covering}\index{covering!simply connected} $\pi:X\rightarrow M$ by a
simply connected manifold $X$. This manifold is unique up to equivalence of
coverings. Furthermore any covering $\psi:Y\rightarrow M$ can be covered by
the universal cover $X$ via a covering map\index{covering map}\index{map!covering} $\Psi:X\rightarrow Y$ such that
$\pi=\psi\circ\Psi$. This is due to the fact that every covering of $M$
corresponds to a subgroup of the fundamental group\index{fundamental
  group}\index{group!fundamental} $\pi_1(M)$
\index{$pi 1(M)$@$\pi_1(M)$}\index{pi 1(M)@$\pi_1(M)$}
of $M$. (rf. \cite[III.8]{Bredon}) The unique covering $\phi:X\rightarrow M$ is
called the \emph{universal covering} of $M$.
\index{universal covering}\index{covering!universal}
Very often, the universal covering is denoted by
$\tilde M$.\index{M@$\tilde M$}\index{$m$@$\tilde M$}

\subsection*{Tangent Bundle}
A map $f:M\rightarrow N$ between two smooth (resp.\ real analytic) manifolds
is called \emph{smooth}
\index{smooth map}\index{smooth map!in $p$}
\index{map!smooth}\index{map!analytic}
(resp.\ real \emph{analytic}) \emph{in $p\in M$} if for all charts
$(U,\phi)\in{\cal A}(M)$ with $p\in U$ and $(V,\psi)\in{\cal A}(N)$ with
$f(p)\in V$ the map $\psi\circ f\circ\phi\inv$ is smooth (resp.\ real analytic)
in $\phi(p)$. \\
$f$ is called \emph{smooth} (resp.\ real analytic) if it is smooth (resp.\ real
analytic) in every point.\\

Write ${\frak F}(p)$
\index{F(p)@${\frak F}(p)$}\index{$f(p)$@${\frak F}(p)$}
for the set of all \emph{germs}\index{germs}%
\footnote{We get the set of germs by considering all functions which are
  defined on a neighbourhood of $p$ and identifying two of them if
  they are identical on a (maybe smaller) neighbourhood of $p$.}
of smooth real valued functions $f:U\rightarrow\mathbbR$ defined on an open
neighbourhood $U$ of $p$ in the
smooth manifold $M$. This is an $\mathbbR$-algebra.
\\
A map $X:{\frak F}(p)\rightarrow\mathbbR$ is called a
\emph{derivation}
\index{derivation!of ${\frak F}(p)$}
of ${\frak F}(p)$, if it satisfies for all $a,b\in\mathbbR$ and
$f,g\in{\frak F}(p)$ 
\begin{align}
X(af+bg)&=aX(f)+bX(g)\\
X(fg)&=f(p)X(g)+g(p)X(f).
\end{align}
The set of all derivations of ${\frak F}(p)$ is called the \emph{tangent
  space}\index{tangent space}
 of
$M$ at $p$, is denoted by $T_pM$\index{T p M@$T_pM$}\index{$t p M$@$T_pM$}
and carries the structure of a vector
space in a natural way. Notice that by the above properties $Xf$ depends
only on the behaviour of $f$ on a neighbourhood of $p$%
\footnote{First of all notice that $Xf=0$ if $f|_U\equiv 0$. This is easy
  to see, since we can find $g\in\frak F(p)$ with $g(p)=0$ and $gf\equiv
  f$. Hence $X(f)=X(gf)=f(p)X(g)+g(p)X(f)=0$. By linearity we see that
  $Xf=X \tilde f$ whenever $f$ and $\tilde f$ agree on a neighbourhood of
  $p$. 
}. Therefore we only
considered germs of functions.\\
We introduce the \emph{tangent bundle}
\index{tangent bundle}\index{bundle!tangent}
$TM:=\underset{p\in M}\bigcup T_pM$\index{TM@$TM$}\index{$tm$@$TM$}
and the \emph{base point projection}
\index{projection!base point}\index{base point projection}
\index{p@$p:TM\rightarrow M$}\index{$p:TM\rightarrow M$@$p:TM\rightarrow M$}
$p:TM\rightarrow M$ which is defined by $p(X)=p$ whenever
$X\in T_pM$. The tangent bundle carries the structure of a $2n$-dimensional
smooth manifold. \\
For any smooth curve $c:I\rightarrow M$ we define the \emph{tangent
  vector}\index{tangent vector}\index{vector!tangent} $\dot c(0)$
\index{c(0)@$\dot c(0)$}\index{$c(0)$@$\dot c(0)$}
to $c$ in
$p:=c(0)$ to be the derivation of ${\frak F}(p)$ defined by
$$\dot c(0):f\longmapsto \frac{\partial}{\partial t}\big|_{t=0}f\circ c(t).$$
Notice that every element of $T_pM$ can be realised as tangent vector to a
smooth curve.\\

For a smooth map $f:M\rightarrow N$ we define $d_pf:T_pM\rightarrow
T_{f(p)}N$,\index{d p f@$d_pf$}\index{$d p f$@$d_pf$} the
\emph{differential in $p$},\index{differential in $p$} by 
$$d_pf(X)g:= X(g\circ f)\quad\text{for }X\in T_pM\;\;g\in{\frak F}(f(p))$$
and the \emph{differential}\index{differential} $df:TM\rightarrow TN$
\index{df@$df$}\index{$df$} to be the collection of the
$d_pf$. Each of the 
$d_pf$ is a linear map between vector spaces, $df$ is a smooth map between
smooth manifolds. Notice that the chain rule applies: Given smooth maps
$L\overset f\rightarrow M\overset g\rightarrow N$ and $p\in L$ we have
$d_p(g\circ f)=d_{f(p)}g\circ d_pf$.

\begin{rema}{}{}{}
Consider again the case where $M$ is an open subset of
$V\equiv\mathbbR^n$. For every $p\in M$ we can identify $V$ with $T_pM$ by
assigning to $v\in V$ the derivation
\begin{align*}
v:\frak F(p)&\longrightarrow \mathbbR\\
f&\longmapsto \frac d{dt}\Big|_{t=0}f(p+tv)
\end{align*}
In the special case where we have an open interval of the real line, $I:=]a,b[\subset\mathbbR$, a smooth map
$c:I\rightarrow M$ is called a curve in $M$. The differential map
$d_tc:T_tI\equiv\mathbbR\rightarrow TM$ is easily understood if we note
that $d_tc(\lambda)= \lambda d_tc(1)$ and $d_tc(1)=:\dot c(t)$ is given by
$$
\begin{array}[t]{cccc}
\dot c(t):&\frak F(c(t))&\rightarrow&\mathbbR\\
         &  f&\mapsto&\frac{\partial}{\partial s}\big|_{s=t}f\circ c(s).
\end{array}$$

\end{rema}

\subsection*{Vector Fields}
A smooth map $X:M\rightarrow TM$ satisfying $X_p\in T_pM$ for all $p\in M$
is called a smooth \emph{vector field} \index{vector field}on $M$. We denote the space of all vector
fields on $M$ by $\frak X(M)$.
\index{X(M)@${\frak X}(M)$}\index{$x(m)$@${\frak X}(M)$}
This is a vector space (even ${\frak
  F}(M)$-module), the space of sections of the vector bundle
$p:TM\rightarrow M$.\\
A vector field $X$ acts as a \emph{derivation%
\footnote{I.\ e.\ for all $a,b\in\mathbbR$ and $f,g\in{\frak F}(M)$ 
$$X(af+bg)=aXf+bXg\qquad\text{and}\qquad X(fg)=fXg+gXf.$$
}
 of ${\frak F}(M)$},\index{F(M)@${\frak F}(M)$}\index{$f(m)$@${\frak F}(M)$}
the \emph{vector space of differentiable functions on $M$}
\index{vector space!of differentiable functions on $M$}
by assigning to $f\in{\frak F}(M)$ the function $Xf: p\mapsto X_pf$. \\

We define the \emph{bracket}\index{bracket} on ${\frak X}(M)$ by 
$$[X,Y]f:=X(Yf)-Y(Xf)\quad\text{for all }X,Y\in {\frak X}(M)\text{ and
  }f\in{\frak F}(M)\index{[.,.]@$[.,.]$}\index{$[.,.]$}$$
which has the following properties:

\begin{align*}
[X,Y]&=-[X,Y]\\
[fX,Y]&=f[X,Y]-(Yf)X\\
[X+Y,Z]&=[X,Z]+[Y,Z]\\
[X,[Y,Z]]&=[[X,Y],Z]+[Y,[X,Z]]\tag{Jacobi Identity}\index{Jacobi identity}
\end{align*}
and hence $({\frak X}(M),[.,.])$ is a Lie algebra.\\
A \emph{linear connection}\index{linear connection}\index{connection!linear} $\nabla$\index{nabla@$\nabla$}\index{$nabla$@$\nabla$} on the tangent bundle $p:TM\rightarrow M$ is a map
\begin{align*}
\nabla:{\frak X}(M)\times{\frak X}(M)&\longrightarrow{\frak X}(M)\\
(X,Y)&\longmapsto \nabla_XY
\end{align*}
such that for all $X,Y,Z\in{\frak X}(M)$ and $f,g\in{\frak F}(M)$ 
\begin{align*}
\nabla_{fX+gY}Z&=f\nabla_XZ+g\nabla_YZ\\
\nabla_X(Y+Z)&=\nabla_XY+\nabla_XZ\\
\nabla_X(fY)&=(Xf)Y+f\nabla_XY.
\end{align*}
$\nabla_XY$\index{nabla X Y@$\nabla_XY$}\index{$nabla x y$@$\nabla_XY$} is called the \emph{covariant derivative}\index{covariant derivative}\index{derivative!covariant} of $Y$ via $X$. Notice that
$(\nabla_XY)_p$ depends only on $X_p$ and on the behaviour of $Y$ in a
small neighbourhood of $p$.%
\footnote{Fix a $\frak F(M)$-basis $E_i$ of $\frak X(M)$. Then there are
  functions $x_i\in\frak F(M)$ with $X=\sum x_iE_i$. We see
$$(\nabla_XY)_p=(\sum x_i\nabla_{E_i}Y)_p=\sum x_i(p)(\nabla_{E_i}Y)_p$$
and since the $E_i$ and $Y$ are independent of $X$ the term depends only on
the coeffients of $X_p=\sum x_i(p)E_i(p)$.\\
Suppose $Y$ vanishes on a neighbourhood $U$ of $p$. In this case we can
find $f\in \frak F(p)$ such that $f|_U=0$ and $Y=fY$. Notice that
$(X(f))_p=0$ in this case. Now
$$(\nabla_XY)_p=(\nabla_X(fY))_p=(X(f))_pY+\drunter{=0}{f(p)}(\nabla_XY)_p=0$$
and by linearity we see that $(\nabla_XY)_p$ depends only on the local
behaviour of $Y$.
}

\subsection*{Curvature}
For any linear connection $\nabla$ we define the \emph{curvature
  tensor}\index{curvature tensor}\index{tensor!curvature}
\index{R@$R$}\index{$r$@$R$}
\begin{gather*}
R:{\frak X}(M)\times{\frak X}(M)\times{\frak
  X}(M)\rightarrow{\frak X}(M)\\ 
\text{ by }\quad\quad\quad R(X,Y)Z:=\nabla_X\nabla_YZ-\nabla_Y\nabla_XZ-\nabla_{[X,Y]}Z\qquad\qquad\\
\text{ or }\quad\quad(R(X,Y)=\nabla_X\nabla_Y-\nabla_Y\nabla_X-\nabla_{[X,Y]}:{\frak
  X}(M)\rightarrow{\frak X}(M)).
\end{gather*}
The \emph{torsion tensor}
\index{torsion tensor}\index{tensor!torsion}\index{T@$T$}\index{$t$@$T$}
$T:{\frak X}(M)\times{\frak X}(M)\rightarrow{\frak X}(M)$
of the linear connection $\nabla$ is defined by
$T(X,Y):=\nabla_XY-\nabla_YX-[X,Y]$. The connection $\nabla$ is called
\emph{torsion free}\index{torsion free}, if $T\equiv0$.\\
The \emph{Ricci tensor}
\index{Ricci!tensor}\index{tensor!Ricci}
\index{Ric(X,Y)@$\Ric(X,Y)$}\index{$ric(x,y)$@$\Ric(X,Y)$}
$\Ric(X,Y)\in\mathbbR$ is defined for $X,Y\in T_pM$ by
$$\Ric(X,Y):=\operatorname{trace}\left(Z\rightarrow R(X,Y)Z\right),$$
sometimes with the factor $\frac 1{n-1}$. The \emph{Ricci
  curvature}\index{curvature!Ricci}\index{Ricci curvature} is easily
defined by \index{Ric(X)@$\Ric(X)$}\index{$ric(X)$@$\Ric(X)$}
$$\Ric(X):=\Ric(X,X)$$
but usually only used for vectors of norm one in the Riemannian case
(see below).

\subsection*{Parallel Transport}
Consider a smooth 
curve $c:I\rightarrow M$. A \emph{vector field along $c$}
\index{vector field!along $c$}
is a smooth map $X:I\rightarrow TM$ with $X(t)\in T_{c(t)}M$. We often write
$X_{c(t)}$
\index{X c(t)@$X_{c(t)}$}\index{$x c(t)$@$X_{c(t)}$}
\index{X t@$X_t$}\index{$x t$@$X_t$}
or $X_t$ for $X(t)$. The \emph{vector space of vector fields along $c$}
\index{vector space!of vector fields along $c$}
is denoted by $\frak X(c)$.
\index{X(c)@$\frak X(c)$}\index{$x(c)$@$\frak X(c)$}
Notice that locally any vector field along a curve with $\dot
c\neq 0$ can be extended to a vector field on a neighbourhood of $c$.\\
Every vector field $Y\in\frak X(M)$ defines a vector field along $c$,
namely $Y\circ c$. It is easy to show that the vector field $t\rightarrow
(\nabla_{\dot c}Y)_{c(t)}$ depends only on $Y\circ c$. Hence we can define
the \emph{covariant derivative along $c$}\index{covariant derivative along
  a curve}\index{derivative!covariant, along a curve} as a map
\index{nabla/dt@$\frac\nabla{dt}$} \index{$nabla/dt$@$\frac\nabla{dt}$}
$$\frac\nabla{dt}:\frak X(c)\rightarrow\frak X(c)$$
defined by $\frac\nabla{dt}Y:=(\nabla_{\dot c}\tilde Y)\circ c$, where
$\tilde Y\in\frak X(M)$ is a vector field with $Y=\tilde Y\circ c$. A
shorter way to denote $\frac\nabla{dt}Y$ is to write $Y'$,
\index{X'@$X'$}\index{$x'$@$X'$}
if there is no danger of confusion. \\
A vector field $Y$ along a curve $c$ is called
\emph{parallel}\index{parallel vector field}\index{vector field!parallel} if
$\frac\nabla{dt}Y\equiv0$. For every vector $v\in T_{c(t)}M$ there is
exactly one parallel vector field $P$ along $c$ with $P_{c(t)}=v$. We
define the \emph{parallel transport}\index{parallel transport}\index{transport!parallel} of $v$ along $c|_{[t,s]}$ by
$_c\|_t^sv:=P_{c(s)}$.
\index{$"\"|$gamma@$_\gamma"\"|_t^s$}
%\index{"\"|c@$_c"\"|_t^s$}
This is an isomorphism of vector spaces
$_c\|^s_t:T_{c(t)}M\rightarrow T_{c(s)}M$. We write $\frak P(c)$
\index{P(c)@$\frak P(c)$}\index{$p(c)$@$\frak P(c)$}
for the
vector space of parallel vector fields along $c$.\\
A curve $c:I\rightarrow M$ is called a \emph{geodesic}\index{geodesic} if its velocity field $\dot
c:I\rightarrow TM$ is a parallel vector field along $c$, i.\ e.\ if
$\frac\nabla{dt}\dot c\equiv 0$ (or $\nabla_{\dot c}\dot c\equiv 0$)%
\footnote{ 
 Notice that there is a different definition for geodesics in metric spaces
 (see below). However for Riemannian manifolds equipped with the
 Levi-Civita connection both definition coincide.}.
By parallel transport we can see that for every $v\in TM$ there is a unique
maximal
geodesic $c$ with $\dot c(0)=v$. We denote this curve by
$\gamma_v$.
\index{gamma v@$\gamma_v$}\index{$gamma v$@$\gamma_v$}
Obviously
\begin{xalignat*}{2}
\gamma_v(0)&=p(v)&\dot \gamma_v(0)&=v\\
\gamma_v(t)&=p(_c\|_0^tv)&\dot \gamma_v(t)&={}_c\|_0^tv.
\end{xalignat*}

With parallel transport and geodesics we can define the \emph{geodesic
  flow}
\index{geodesic flow}\index{flow!geodesic}
\index{phi t@$\phi_t$}\index{$phi t$@$\phi_t$}
$$\begin{array}{cccc}
\phi:&\mathbbR\times TM&\longrightarrow&TM\\
&(t,v)&\longmapsto&\phi_t(v):=\dot\gamma_v(t).
\end{array}$$

\subsection*{Riemannian Metric}
A smooth \emph{Riemannian metric}
\index{metric!Riemannian}\index{Riemannian metric}
$g$ on a smooth manifold assigns to every point
$p\in M$ an inner product
\index{inner product}\index{product!inner}\index{g p@$g_p$}\index{$g p$@$g_p$}
$g_p$ on the tangent space $T_pM$ at $p$, such
that for all smooth vector fields $X,Y\in{\frak X}(M)$ the map $g(X,Y):p\mapsto
g_p(X_p,Y_p)$ is smooth. Instead of $g$ we will usually denote the
Riemannian metric by $\langle.,.\rangle$
\index{<.,.>@$\langle.,.\rangle$}\index{$<.,.>$@$\langle.,.\rangle$}
and the inner product on $T_pM$ by
$\langle.,.\rangle_p$.
\index{<.,.>p@$\langle.,.\rangle_p$}\index{$<.,.>p$@$\langle.,.\rangle_p$}
\begin{theo}{Levi/Civita}{}{}\index{Levi/Civita}\index{Theorem!Levi/Civita}
There is a unique linear connections $\nabla$ on the tangent bundle
$p:TM\rightarrow M$, such that for all $X,Y,Z\in{\frak X}(M)$ 
\begin{align*}
[X,Y]&=\nabla_XY-\nabla_YX\tag{torsion free}
\index{torsion free connection}\index{connection!torsion free}\\
X\langle Y,Z\rangle&=\langle\nabla_XY,Z\rangle+\langle
Y,\nabla_XZ\rangle.\tag{metric}
\index{metric connection}\index{connection!metric}
\end{align*}
This connection is called the \emph{Levi-Civita connection}
\index{connection!Levi-Civita}\index{Levi-Civita connection}
and is defined by
\begin{align*}
\langle\nabla_XY,Z\rangle=\frac12\Big(&  X\langle Y,Z\rangle
                                        +Y\langle Z,X\rangle
                                        -Z\langle X,Y\rangle \\
                                      & +\langle[X,Y],Z\rangle
                                        -\langle[Y,Z],X\rangle
                                        +\langle[Z,X],Y\rangle
                                 \Big)
\end{align*}
\end{theo}
In general we will only consider the Levi-Civita connection and so e.\ g.\ if
we talk about the curvature tensor it will always be the curvature tensor
\index{curvature!Levi-Civita}\index{Levi-Civita curvature}
obtained from the Levi-Civita connection.\\
We define the \emph{gradient}
\index{gradient}\index{nabla f@$\nabla f$}\index{$nabla f$@$\nabla f$}
$\nabla f\in{\frak X}(M)$ of a function $f\in{\frak F}(M)$ by 
$$\langle\nabla f,X\rangle=Xf\qquad\text{ for all } X\in{\frak X}(M).$$
The \emph{Hessian}\index{Hessian}\index{H f@${\cal H}f$}\index{$h f$@${\cal H}f$}
${\cal H}f$ is defined by
$${\cal H}f(X,Y):=\langle\nabla_X\nabla f,Y\rangle=XYf-(\nabla_XY)f.$$
The \emph{divergence}\index{divergence}
\index{div X@$\operatorname{div} X$}
\index{$div x$@$\operatorname{div} X$}
$\operatorname{div} X$ of a vector field $X\in{\frak X}(M)$ is defined by
$$\operatorname{div} X:=\operatorname{trace}(Y\mapsto\nabla_YX)$$
The \emph{Laplacian}\index{Laplacian}
\index{D f@$\Delta f$}\index{$d f$@$\Delta f$}
$\Delta f\in{\frak F}(M)$ of a function $f\in {\frak F}(M)$ is
defined by
$$\Delta f=-\operatorname{trace}({\cal H}f)=-\operatorname{div}(\nabla f)$$

\subsection*{Riemannian Curvature}
\index{Riemannian curvature}\index{curvature!Riemannian}
For manifolds with a linear connection we explained the curvature tensor $R$
before. Now with the scalar product on $T_pM$ we can introduce the
sectional curvature.\\
For a two-dimensional subspace $\sigma$ of $T_pM$ choose a orthonormal
basis $\{u,v\}$. The \emph{sectional curvature}
\index{sectional curvature}\index{curvature!sectional}
$K_\sigma$\index{K sigma@$K_\sigma$}\index{$k sigma$@$K_\sigma$}
of $\sigma$ is defined
by
$$K_\sigma:=\langle R(u,v)v,u\rangle$$ 
which does not depend on the choice of our orthonormal basis. In fact for
any two non-collinear vectors $u,v\in T_pM$ we have 
$$K(u,v):=K_\sigma=\frac{\langle R(u,v)v,u\rangle}{\|u\|^2\|v\|^2-\langle
  u,v\rangle^2}\index{K(u,v)@$K(u,v)$}\index{$k(u,v)$@$K(u,v)$}$$
where $\sigma=\operatorname{span}(u,v)$. In the case of a two-dimensional
manifold the tangent space in every point is two-dimensional and hence
we can assign the curvature to the point. This curvature is usually called
the \emph{Gaussian curvature}\index{curvature!Gaussian}\index{Gaussian curvature}. For higher dimensional manifolds the curvature
tensor can be expressed by the sectional curvature.\\
If we fix an orthonormal basis $\{e_i\}_{i=1..n}$ of $T_pM$ we get nice
expressions for the Ricci curvature.
\begin{align}
\tag{Ricci Tensor}\index{Ricci tensor}\index{tensor!Ricci} \Ric(X,Y)&=\sum\langle R(e_i,X)Y,e_i\rangle\\
\tag{Ricci Curvature}\index{Ricci curvature}\index{curvature!Ricci} \Ric(X)&=\underset{e_i\neq X}\sum K(X,e_i)
\end{align}
Furthermore we can define the \emph{scalar curvature}
\index{scalar curvature}\index{curvature!scalar}
\index{K(p)@$K(p)$}\index{$k(p)$@$K(p)$}
$K:M\rightarrow\mathbbR$ by
$$K(p):=\sum \Ric(e_i)$$
which is the trace of the symmetric bilinear form $\Ric$ and hence
independent of the choice of the orthonormal basis.\\
Sometimes the Ricci tensor differs from this definition by a factor of
$\frac1{n-1}$. In this case the scalar curvature is usually defined by
$\frac 1n\sum\Ric(e_i)$ and hence differs by the factor $\frac1{n(n-1)}$.\\

Throughout this thesis we will only consider manifolds of nonpositive
curvature\index{curvature!nonpositive}\index{nonpositive curvature}. That means for any two-dimensional subspace $\sigma$ of any
tangent space $T_pM$ we have $K_\sigma\leq 0$. 

\subsection*{Geodesics}
Given the Riemannian metric we can define the \emph{length}\index{length} of any piecewise
smooth curve $c:[a,b]\rightarrow M$ by
$$
\index{L(c)@$L(c)$}\index{$l(c)$@$L(c)$}
L(c):=\overset b{\underset a\int}\|\dot c(t)\|dt$$
where $\|u\|:=\sqrt{\langle u,u\rangle}$ denotes the \index{norm}\emph{norm}\index{$"\"|."\"|$@$"\"|."\"|$} associated to the
scalar product $\langle.,.\rangle$. A curve $c$ is said to be \emph{parametrised
by arc length}\index{arc length}\index{parametrisation!by arc length}, if $\|\dot c(t)\|\equiv 1$. In this case $L(c)=|b-a|$. \\
We define a \emph{metric} on $M$ by introducing the \emph{distance}\index{metric}\index{distance}\index{d(p,q)@$d(p,q)$}\index{$d(p,q)$@$d(p,q)$} $d(p,q)$ between two
points $p,q\in M$:
$$d(p,q):=\inf\{L(c)\;|\;c:[a,b]\rightarrow M\text{ smooth, } c(a)=p\text{
  and }c(b)=q\}.$$
$(M,d)$ is a \emph{metric space}\index{metric space}\index{space!metric},
i.\ e.\ for all $p,q,r\in M$ 
\begin{align}
&\tag{positivity}\index{positivity} d(p,q)\geq 0 \text{ and }(d(p,q)=0 \text{ implies }p=q)\\
&\tag{symmetry}\index{symmetry} d(p,q)=d(q,p)\\
&\tag{triangle inequality}\index{triangle inequality} d(p,r)\leq d(p,q)+d(q,r).
\end{align}
Now we want to find curves that are the shortest connections between two
given points. To this end we need to consider variations.\\
A \emph{smooth variation}\index{variation} is a smooth map $h:I\times [a,b]\rightarrow M,$
 $(s,t)\mapsto h(s,t)=:h_s(t)$. The curves $h_s:[a,b]\rightarrow M$ are
 considered neighbouring curves of the curve $h_0$. The variation is called
 a variation of $h_0$.\\
Every variation defines a vector field $X_t:=\frac\partial{\partial
  s}\big|_{s=0}h(s,t)$ along $h_0$. This vector field is called \emph{variational
vector field}\index{vector field!variational}\index{variational vector field} of $h$. For every vector field $X$ along a curve $c$ we can find a
variation $h$ such that $h_0=c$ and $X$ is the variational vector field of
$h$.\\
A variation $h$ of a curve $c:[a,b]\rightarrow M$ is called
\emph{proper}\index{proper variation} if $h_s(a)=c(a)$ and $h_s(b)=c(b)$ for
all $s\in I$. In this case the
variational vector field vanishes in $c(a)$ and $c(b)$. On the other hand
for every vector field that vanishes at the end points we can
find a proper variation such that the vector field is its variational
vector field.\\
The variational vector fields gives us informations about the behaviour of
the lengths $L_s:=L(h_s)$. Namely if $c=h_0$ there is the \emph{first
  variation formula}
\index{first variation formula}\index{variation formula!first}
$$L'(s)=-\int_a^b\left\langle X(t),\frac\nabla{dt}\frac{\dot c(t)}{\|\dot
  c(t)\|}\right\rangle\;dt+\left\langle X(b),\frac{\dot c(b)}{\|\dot c(b)\|}\right\rangle-
 \left\langle X(a),\frac{\dot c(a)}{\|\dot c(a)\|}\right\rangle.$$
This formula holds for the Levi-Civita connection. Since this connection
is torsion free we have
$$\frac\nabla{\partial s}\frac\partial{\partial t}h=\frac\nabla{\partial
  t}\frac\partial{\partial s}h\quad\text{ and hence }\quad \left[\frac\partial{\partial t}h,\frac\partial{\partial s}h\right]=0.$$

\begin{rema}{}{}{}
It is possible to consider piecewise smooth
variations\index{variation!piecewise smooth}. We then get special
terms belonging to points of indifferentiability. For our purposes,
however, it will be sufficient to consider smooth variations.
\end{rema}

Considering only proper variations we see that a necessary
condition for a curve to be the shortest connection from $c(a)$ to $c(b)$
is $\left\langle X(t),\frac\nabla{dt}\frac{\dot c(t)}{\|\dot
  c(t)\|}\right\rangle\equiv 0$ for all $t\in[a,b]$. Since this must hold for
all variational vector fields we find $\frac\nabla{dt}\frac{\dot c(t)}{\|\dot
  c(t)\|}\equiv 0$.\\
I.\ e.\ every shortest curve connecting two points (if
such a curve exists) is a geodesic up to parametrisation.\\
In general we will only consider manifolds where the shortest connecting
curve between any two points exist. By the following theorem this holds
whenever the underlying metric space is complete.

\begin{theo}{Hopf/Rinow}{Hopf-Rinow}{}
\index{Hopf/Rinow}\index{Theorem!Hopf/Rinow}
The following conditions are equivalent for a connected Riemannian manifold
$(M,g)$.
\begin{enumerate}
\item \label{complete} Every Cauchy-sequence in $(M,d)$ converges.
\item\label{geodcompp} There is a point $p\in M$ such that every geodesic emanating from $p$
  can be extended to $\mathbbR$.
\item\label{geodcomp} Every geodesic in $M$ can be extended to $\mathbbR$.
\item Any closed ball $\bar B_r(p):=\{q\in M\;|\;d(p,q)\leq r\}$ is compact.
\end{enumerate}
Any of these implies
\begin{enumerate}\setcounter{enumi}{4}
\item\label{geodesicspace} Any two points in $M$ can be joined by a geodesic realising the distance.
\end{enumerate}
\end{theo}
\begin{defi}{}{}{}
The properties mentioned in Theorem~\ref{Hopf-Rinow} are named:
\begin{enumerate}
\item$(M,d)$ is a \emph{complete} metric space\index{complete metric
    space}\index{metric space!complete}.
\item $(M,g)$ is \emph{geodesically complete in $p$}.\index{geodesically
    complete Riemannian manifold}
\item $(M,g)$ is \emph{geodesically complete}\index{complete!geodesically}.
\item $(M,d)$ is \emph{proper}.\index{proper space}\index{space!proper}
\item $(M,d)$ is a \emph{geodesic space}\index{geodesic space}\index{space!geodesic}.
\end{enumerate}
\end{defi}

We will only be interested in \emph{complete manifolds}.\index{complete
  manifold}\index{manifold!complete} % (in any of the
%above senses).
On a complete manifold we can define the \emph{exponential
  map}\index{exponential map}
$\exp:TM\rightarrow M$ by $\exp(v):=\gamma_v(1)\index{exp@$\exp$}$\index{$exp$@$\exp$}. 
Write $\exp_p$ for the restriction of $\exp$ to the tangent space $T_pM$. 
The \emph{injectivity radius}\index{injectivity radius}\index{radius of
  injectivity} in a point $p\in M$ is the largest number $\iota>0$ such
that $\exp_p$ is injective on the ball $$\{u\in T_pM\;|\;\|u\|\leq\iota\}.$$

\subsection*{Jacobi Fields}
Consider a variation $h$ of a geodesic $\gamma=h_0$. If all $h_s$ are
geodesics then the variational vector field $X$ will satisfy the \emph{Jacobi
equation}\index{Jacobi equation}\index{equation!Jacobi}
$$X''+R(X,\dot\gamma)\dot\gamma\equiv 0$$
as the following calculation shows.
{\begin{align*}
X''&=\frac\nabla{\partial t}\frac\nabla{\partial t}\frac{\partial}{\partial
  s}h\\
&=\frac\nabla{\partial t}\frac\nabla{\partial s}\frac{\partial}{\partial
  t}h\\
&=R(\frac\partial{\partial t}h,\frac\partial{\partial
  s}h)\frac{\partial}{\partial t}h+\nabla_{\left[\frac\partial{\partial
      t}h,\frac\partial{\partial s}h\right]}\frac\partial{\partial t}h
+\frac\nabla{\partial s}\drunter{=\frac\nabla{\partial t}\dot{h_s}}{\frac\nabla{\partial t}\frac{\partial}{\partial t}h}\\
&=-R(\frac\partial{\partial
  s}h,\frac\partial{\partial t}h)\frac{\partial}{\partial t}h\\
&=-R(X,\dot\gamma)\dot\gamma
\end{align*}
Notice that $\frac\nabla{\partial t}=\nabla_{\frac{\partial}{\partial t}h}$
and $\frac\nabla{\partial s}=\nabla_{\frac{\partial}{\partial s}h}$ and
that $\frac\nabla{\partial t}\dot{h_s}=0$ since the $h_s$ are geodesics.
}

Any vector field along $\gamma$ satisfying the Jacobi equation will be called a \emph{Jacobi
field}\index{Jacobi field}\index{vector field!Jacobi} and in fact every Jacobi field gives rise to a variation of $\gamma$
consisting only of geodesics.\\
Denote by $\frak J(\gamma)$
\index{J(gamma)@$\frak J(\gamma)$}\index{$j(gamma)$@$\frak J(\gamma)$}
the vector space of all Jacobi fields along the
geodesic $\gamma$. This is a $2n$-dimensional vector space and can be
identified with $T_{\gamma(0)}M\times T_{\gamma(0)}M$ since for every
choice of $v,w\in T_{\gamma(0)}M$ there is a unique Jacobi field $X$ along
$\gamma$ with $X_0=v$ and $X'_0=w$. \\
$\frak J(\gamma)$ is the direct sum of the vector spaces of orthogonal and
parallel Jacobi fields $\frak J^\perp(\gamma)$
\index{J(gamma)orth@$\frak J^\perp(\gamma)$}\index{$j(gamma)orth$@$\frak J^\perp(\gamma)$}
\index{J(gamma)par@$\frak J^{"\"|}(\gamma)$}\index{$j(gamma)par$@$\frak J^{"\"|}(\gamma)$}
and $\frak J^\|(\gamma)$. For any Jacobi field $X$ the longitudinal component
$X^\|:=\frac{\langle X,\dot\gamma\rangle}{\|\dot\gamma\|^2}\dot\gamma$ and
the orthogonal component $X^\perp:=X-X^\|$ are both Jacobi fields. It is
easy to see that
$J^\|(\gamma)=\{t\rightarrow(x+vt)\dot\gamma(t)\;|\;x,v\in\mathbbR\}$. 
Suppose we have a variation $h$ consisting entirely of geodesics of the
same speed ($\|\dot h_s\|=const$.). In this case $X^\|=x\dot\gamma$ since 
\begin{align*}
\langle X',\dot\gamma\rangle&=\langle\frac{\nabla}{\partial
  t}\frac\partial{\partial s}h,\frac\partial{\partial t}h\rangle\\
&=\langle\frac{\nabla}{\partial
  s}\frac\partial{\partial t}h,\frac\partial{\partial t}h\rangle\\
&=\frac12\frac d{ds}\|\dot{h_s}\|^2=0.
\end{align*}
Given this variation we can define a variation by $f(s,t):=h(s,t-xs)$ with
variational vector field in $\frak J^\perp(\gamma)$ and $f_s\simeq h_s$ up to
parametrisation. Since we are only interested in geodesic variations with
geodesics of constant speed we will only have to deal with orthogonal
Jacobi fields. Notice further that for a geodesic variation (of constant
speed) with fixed
starting point the variational vector field is always an orthogonal Jacobi field.

\end{appendix}
\cleardoublepage
\pagestyle{myheadings}
\markboth{Index}{Index}
\thispagestyle{empty}
\addcontentsline{toc}{section}{Index}
\printindex
\cleardoublepage
\bibliographystyle{alpha} 
\addcontentsline{toc}{section}{Bibliography}

\markboth{Bibliography}{Bibliography}
\thispagestyle{empty}
\bibliography{archiv}
\cleardoublepage
\section*{Notes}
\markboth{}{}
\thispagestyle{empty}
\cleardoublepage
\section*{}
\cleardoublepage
\section*{}
\cleardoublepage
\end{document}